\documentclass[10pt]{article}
\usepackage{amsmath,amssymb,amsfonts}     
\date{}
\newtheorem{lemma}{Lemma}[section]
\newtheorem{theorem}{Theorem}[section]
\newtheorem{proposition}{Proposition}[section]
\newcommand{\eproof}{\makebox[1cm]{}\hfill{\framebox[2.5mm]{}}}
\begin{document}
\title{\bf Designs for 24-vertex snarks\\}
\author{Anthony~D.~Forbes\\
        Department of Mathematics and Statistics\\
        The Open University\\
        Walton Hall\\
        Milton Keynes MK7 6AA\\
        UNITED KINGDOM}
\maketitle
\begin{abstract}
The design spectrum problem is solved for the thirty-eight 24-vertex non-trivial snarks.
\end{abstract}

\newcommand{\adfsplit}{\\ \makebox{~~~~~~~~~~}}
\newcommand{\adfSgap}{\vskip 0.5mm}
\newcommand{\adfLgap}{\vskip 1.0mm}

\newcommand{\adfVfy}[1]{}

\section{Introduction}
\label{sec:introduction}
Let $G$ be a simple graph without isolated vertices.
If the edge set of a simple graph $K$ can be partitioned into edge sets
of graphs each isomorphic to $G$, we say that there exists a {\em decomposition} of $K$ into $G$.
In the case where $K$ is the complete graph $K_n$ we refer to the decomposition as a
$G$ {\em design} of order $n$. The {\em spectrum} of $G$ is the set of positive integers $n$ for
which there exists a $G$ design of order $n$.
We refer the reader to the survey article of Adams, Bryant and Buchanan, \cite{ABB} and,
for more up to date results, the Web site maintained by Bryant and McCourt, \cite{BM}.

A {\em snark} is a connected, bridgeless 3-regular graph with chromatic index 4.
However, a snark is usually regarded as trivial (or reducible) if it has girth less than 5 or if it has
three edges the deletion of which results in a disconnected graph each of whose components
is non-trivial (as a graph).

The spectra for non-trivial snarks of up to 22 vertices have been successfully determined, as follows.
\begin{enumerate}
\item{The smallest non-trivial snark is the 10-vertex Petersen graph for which
designs of order $n$ exist if and only if
$n \equiv 1$ or 10 (mod 15), \cite{AB-pet}.}
\item{There are no non-trivial snarks on 12, 14 or 16 vertices.}
\item{There are two non-trivial snarks on 18 vertices, namely the $(1,1)$- and $(1,2)$-Blanu\v{s}a snarks.
For each of them, designs of order $n$ exist if and only if
$n \equiv 1$ (mod 27), \cite{For-snark}.}
\item{For each of the six 20-vertex non-trivial snarks (including the flower snark J5),
designs of order $n$ exist if and only if
$n \equiv 1,$ 16, 25, 40 (mod 60), $n \neq 16$, \cite{For-snark}.}
\item{For each of the twenty 22-vertex non-trivial snarks (including the two Loupekine snarks),
designs of order $n$ exist if and only if
$n \equiv 1$ or 22 (mod 33), \cite{For-snark}.}
\end{enumerate}
The purpose of this paper is to extend these results to 24-vertex non-trivial snarks of which there are precisely 38, \cite{BGHM}.
We prove the following.

{\theorem \label{thm:snark24}
Designs of order $n$ exist for each of the thirty-eight non-trivial snarks on $24$ vertices if and only if
$n \equiv 1$ or $64~(\mathrm{mod}~72)$.}\\

\noindent The methods use to obtain results like Theorems~\ref{thm:snark24} are explained
in \cite{For-snark}, to which the reader should regard this paper as a sequel.
The main tool is \cite[Proposition 1.8]{For-snark},
repeated here for convenience as Proposition~\ref{prop:d=3, v=24}.

{\proposition \label{prop:d=3, v=24}
Let $G$ be a 3-regular graph on $24$ vertices.
Suppose there exist $G$ designs of order $64$, $73$, $136$ and $145$.
Suppose also that there exist decompositions into $G$ of the complete multipartite graphs
$K_{12,12,12}$, $K_{24,24,15}$, $K_{72,72,63}$, $K_{24,24,24,24}$ and $K_{24,24,24,21}$.
Then there exists a $G$ design of order $n$ if and only if
$n \equiv 1$ or $64~(\textup{mod}~ 72)$.}\\

\noindent\textbf{Proof.} See \cite{For-snark}.\eproof\\

To prove Theorem~\ref{thm:snark24} it is sufficient, therefore, to construct the
four designs and five multipartite graph decompositions required by Proposition~\ref{prop:d=3, v=24}
for each of the 38 non-trivial snarks.
To save space the details of the constructions are given only for designs of order 136,
which, in the author's opinion, were by far the most difficult to obtain. The details of the
remaining constructions are given in the Appendix,
which is present only in the full version of this paper,
available by request from the author.
The decompositions were obtained and checked by computer programs in the same manner as those of \cite{For-snark}.




\section{The 38 non-trivial snarks on 24 vertices}
\label{sec:Snark24}
\newcommand{\adfGa}{\mathrm{G1}}
\newcommand{\adfGb}{\mathrm{G2}}
\newcommand{\adfGc}{\mathrm{G3}}
\newcommand{\adfGd}{\mathrm{G4}}
\newcommand{\adfGe}{\mathrm{G5}}
\newcommand{\adfGf}{\mathrm{G6}}
\newcommand{\adfGg}{\mathrm{G7}}
\newcommand{\adfGh}{\mathrm{G8}}
\newcommand{\adfGi}{\mathrm{G9}}
\newcommand{\adfGj}{\mathrm{G10}}
\newcommand{\adfGk}{\mathrm{G11}}
\newcommand{\adfGl}{\mathrm{G12}}
\newcommand{\adfGm}{\mathrm{G13}}
\newcommand{\adfGn}{\mathrm{G14}}
\newcommand{\adfGo}{\mathrm{G15}}
\newcommand{\adfGp}{\mathrm{G16}}
\newcommand{\adfGq}{\mathrm{G17}}
\newcommand{\adfGr}{\mathrm{G18}}
\newcommand{\adfGs}{\mathrm{G19}}
\newcommand{\adfGt}{\mathrm{G20}}
\newcommand{\adfGu}{\mathrm{G21}}
\newcommand{\adfGv}{\mathrm{G22}}
\newcommand{\adfGw}{\mathrm{G23}}
\newcommand{\adfGx}{\mathrm{G24}}
\newcommand{\adfGy}{\mathrm{G25}}
\newcommand{\adfGz}{\mathrm{G26}}
\newcommand{\adfGA}{\mathrm{G27}}
\newcommand{\adfGB}{\mathrm{G28}}
\newcommand{\adfGC}{\mathrm{G29}}
\newcommand{\adfGD}{\mathrm{G30}}
\newcommand{\adfGE}{\mathrm{G31}}
\newcommand{\adfGF}{\mathrm{G32}}
\newcommand{\adfGG}{\mathrm{G33}}
\newcommand{\adfGH}{\mathrm{G34}}
\newcommand{\adfGI}{\mathrm{G35}}
\newcommand{\adfGJ}{\mathrm{G36}}
\newcommand{\adfGK}{\mathrm{G37}}
\newcommand{\adfGL}{\mathrm{G38}}
The thirty-eight non-trivial snarks on 24 vertices are represented by ordered 24-tuples of vertices:
(1, 2, \dots, 24)${}_{\adfGa}$, (1, 2, \dots, 24)${}_{\adfGb}$, \dots, (1, 2, \dots, 24)${}_{\adfGL}$.
Their edge sets, as appearing in Royale's list, \cite{Roy}, are respectively

$\adfGa$:
\{$\{1,2\}$, $\{1,3\}$, $\{1,4\}$, $\{2,5\}$, $\{2,6\}$, $\{3,7\}$, $\{3,8\}$, $\{4,9\}$, $\{4,10\}$, $\{5,7\}$, $\{5,9\}$, $\{6,8\}$, $\{6,11\}$, $\{7,12\}$, $\{8,13\}$, $\{9,14\}$, $\{10,11\}$, $\{10,12\}$, $\{11,15\}$, $\{12,16\}$, $\{13,14\}$, $\{13,17\}$, $\{14,18\}$, $\{15,16\}$, $\{15,19\}$, $\{16,20\}$, $\{17,21\}$, $\{17,22\}$, $\{18,23\}$, $\{18,24\}$, $\{19,21\}$, $\{19,23\}$, $\{20,22\}$, $\{20,24\}$, $\{21,24\}$, $\{22,23\}$\},

$\adfGb$:
\{$\{1,2\}$, $\{1,3\}$, $\{1,4\}$, $\{2,5\}$, $\{2,6\}$, $\{3,7\}$, $\{3,8\}$, $\{4,9\}$, $\{4,10\}$, $\{5,7\}$, $\{5,9\}$, $\{6,8\}$, $\{6,11\}$, $\{7,12\}$, $\{8,13\}$, $\{9,14\}$, $\{10,11\}$, $\{10,15\}$, $\{11,16\}$, $\{12,15\}$, $\{12,16\}$, $\{13,14\}$, $\{13,17\}$, $\{14,18\}$, $\{15,19\}$, $\{16,20\}$, $\{17,21\}$, $\{17,22\}$, $\{18,23\}$, $\{18,24\}$, $\{19,21\}$, $\{19,23\}$, $\{20,22\}$, $\{20,24\}$, $\{21,24\}$, $\{22,23\}$\},

$\adfGc$:
\{$\{1,2\}$, $\{1,3\}$, $\{1,4\}$, $\{2,5\}$, $\{2,6\}$, $\{3,7\}$, $\{3,8\}$, $\{4,9\}$, $\{4,10\}$, $\{5,7\}$, $\{5,9\}$, $\{6,8\}$, $\{6,11\}$, $\{7,12\}$, $\{8,13\}$, $\{9,14\}$, $\{10,11\}$, $\{10,15\}$, $\{11,16\}$, $\{12,15\}$, $\{12,16\}$, $\{13,17\}$, $\{13,18\}$, $\{14,19\}$, $\{14,20\}$, $\{15,21\}$, $\{16,22\}$, $\{17,19\}$, $\{17,23\}$, $\{18,20\}$, $\{18,24\}$, $\{19,24\}$, $\{20,23\}$, $\{21,22\}$, $\{21,23\}$, $\{22,24\}$\},

$\adfGd$:
\{$\{1,2\}$, $\{1,3\}$, $\{1,4\}$, $\{2,5\}$, $\{2,6\}$, $\{3,7\}$, $\{3,8\}$, $\{4,9\}$, $\{4,10\}$, $\{5,7\}$, $\{5,9\}$, $\{6,8\}$, $\{6,11\}$, $\{7,12\}$, $\{8,13\}$, $\{9,14\}$, $\{10,11\}$, $\{10,15\}$, $\{11,16\}$, $\{12,15\}$, $\{12,17\}$, $\{13,18\}$, $\{13,19\}$, $\{14,20\}$, $\{14,21\}$, $\{15,22\}$, $\{16,17\}$, $\{16,22\}$, $\{17,23\}$, $\{18,20\}$, $\{18,23\}$, $\{19,21\}$, $\{19,24\}$, $\{20,24\}$, $\{21,23\}$, $\{22,24\}$\},

$\adfGe$:
\{$\{1,2\}$, $\{1,3\}$, $\{1,4\}$, $\{2,5\}$, $\{2,6\}$, $\{3,7\}$, $\{3,8\}$, $\{4,9\}$, $\{4,10\}$, $\{5,7\}$, $\{5,9\}$, $\{6,8\}$, $\{6,11\}$, $\{7,12\}$, $\{8,13\}$, $\{9,14\}$, $\{10,12\}$, $\{10,15\}$, $\{11,12\}$, $\{11,16\}$, $\{13,14\}$, $\{13,17\}$, $\{14,18\}$, $\{15,16\}$, $\{15,19\}$, $\{16,20\}$, $\{17,21\}$, $\{17,22\}$, $\{18,23\}$, $\{18,24\}$, $\{19,21\}$, $\{19,23\}$, $\{20,22\}$, $\{20,24\}$, $\{21,24\}$, $\{22,23\}$\},

$\adfGf$:
\{$\{1,2\}$, $\{1,3\}$, $\{1,4\}$, $\{2,5\}$, $\{2,6\}$, $\{3,7\}$, $\{3,8\}$, $\{4,9\}$, $\{4,10\}$, $\{5,7\}$, $\{5,9\}$, $\{6,8\}$, $\{6,11\}$, $\{7,12\}$, $\{8,13\}$, $\{9,14\}$, $\{10,12\}$, $\{10,15\}$, $\{11,15\}$, $\{11,16\}$, $\{12,16\}$, $\{13,14\}$, $\{13,17\}$, $\{14,18\}$, $\{15,19\}$, $\{16,20\}$, $\{17,21\}$, $\{17,22\}$, $\{18,23\}$, $\{18,24\}$, $\{19,21\}$, $\{19,23\}$, $\{20,22\}$, $\{20,24\}$, $\{21,24\}$, $\{22,23\}$\},

$\adfGg$:
\{$\{1,2\}$, $\{1,3\}$, $\{1,4\}$, $\{2,5\}$, $\{2,6\}$, $\{3,7\}$, $\{3,8\}$, $\{4,9\}$, $\{4,10\}$, $\{5,7\}$, $\{5,9\}$, $\{6,8\}$, $\{6,11\}$, $\{7,12\}$, $\{8,13\}$, $\{9,14\}$, $\{10,12\}$, $\{10,15\}$, $\{11,15\}$, $\{11,16\}$, $\{12,16\}$, $\{13,17\}$, $\{13,18\}$, $\{14,19\}$, $\{14,20\}$, $\{15,21\}$, $\{16,22\}$, $\{17,19\}$, $\{17,23\}$, $\{18,20\}$, $\{18,24\}$, $\{19,24\}$, $\{20,23\}$, $\{21,22\}$, $\{21,23\}$, $\{22,24\}$\},

$\adfGh$:
\{$\{1,2\}$, $\{1,3\}$, $\{1,4\}$, $\{2,5\}$, $\{2,6\}$, $\{3,7\}$, $\{3,8\}$, $\{4,9\}$, $\{4,10\}$, $\{5,7\}$, $\{5,9\}$, $\{6,8\}$, $\{6,11\}$, $\{7,12\}$, $\{8,13\}$, $\{9,14\}$, $\{10,12\}$, $\{10,15\}$, $\{11,15\}$, $\{11,16\}$, $\{12,17\}$, $\{13,18\}$, $\{13,19\}$, $\{14,20\}$, $\{14,21\}$, $\{15,22\}$, $\{16,17\}$, $\{16,23\}$, $\{17,22\}$, $\{18,20\}$, $\{18,23\}$, $\{19,21\}$, $\{19,24\}$, $\{20,24\}$, $\{21,23\}$, $\{22,24\}$\},

$\adfGi$:
\{$\{1,2\}$, $\{1,3\}$, $\{1,4\}$, $\{2,5\}$, $\{2,6\}$, $\{3,7\}$, $\{3,8\}$, $\{4,9\}$, $\{4,10\}$, $\{5,7\}$, $\{5,9\}$, $\{6,8\}$, $\{6,11\}$, $\{7,12\}$, $\{8,13\}$, $\{9,14\}$, $\{10,12\}$, $\{10,15\}$, $\{11,16\}$, $\{11,17\}$, $\{12,16\}$, $\{13,18\}$, $\{13,19\}$, $\{14,20\}$, $\{14,21\}$, $\{15,17\}$, $\{15,22\}$, $\{16,22\}$, $\{17,23\}$, $\{18,20\}$, $\{18,23\}$, $\{19,21\}$, $\{19,24\}$, $\{20,24\}$, $\{21,23\}$, $\{22,24\}$\},

$\adfGj$:
\{$\{1,2\}$, $\{1,3\}$, $\{1,4\}$, $\{2,5\}$, $\{2,6\}$, $\{3,7\}$, $\{3,8\}$, $\{4,9\}$, $\{4,10\}$, $\{5,7\}$, $\{5,9\}$, $\{6,8\}$, $\{6,11\}$, $\{7,12\}$, $\{8,13\}$, $\{9,14\}$, $\{10,15\}$, $\{10,16\}$, $\{11,15\}$, $\{11,17\}$, $\{12,15\}$, $\{12,18\}$, $\{13,19\}$, $\{13,20\}$, $\{14,21\}$, $\{14,22\}$, $\{16,17\}$, $\{16,18\}$, $\{17,23\}$, $\{18,24\}$, $\{19,21\}$, $\{19,23\}$, $\{20,22\}$, $\{20,24\}$, $\{21,24\}$, $\{22,23\}$\},

$\adfGk$:
\{$\{1,2\}$, $\{1,3\}$, $\{1,4\}$, $\{2,5\}$, $\{2,6\}$, $\{3,7\}$, $\{3,8\}$, $\{4,9\}$, $\{4,10\}$, $\{5,7\}$, $\{5,9\}$, $\{6,8\}$, $\{6,11\}$, $\{7,12\}$, $\{8,13\}$, $\{9,14\}$, $\{10,15\}$, $\{10,16\}$, $\{11,15\}$, $\{11,17\}$, $\{12,15\}$, $\{12,18\}$, $\{13,19\}$, $\{13,20\}$, $\{14,21\}$, $\{14,22\}$, $\{16,18\}$, $\{16,23\}$, $\{17,18\}$, $\{17,24\}$, $\{19,21\}$, $\{19,23\}$, $\{20,22\}$, $\{20,24\}$, $\{21,24\}$, $\{22,23\}$\},

$\adfGl$:
\{$\{1,2\}$, $\{1,3\}$, $\{1,4\}$, $\{2,5\}$, $\{2,6\}$, $\{3,7\}$, $\{3,8\}$, $\{4,9\}$, $\{4,10\}$, $\{5,7\}$, $\{5,9\}$, $\{6,8\}$, $\{6,11\}$, $\{7,12\}$, $\{8,13\}$, $\{9,14\}$, $\{10,15\}$, $\{10,16\}$, $\{11,15\}$, $\{11,17\}$, $\{12,16\}$, $\{12,18\}$, $\{13,19\}$, $\{13,20\}$, $\{14,21\}$, $\{14,22\}$, $\{15,18\}$, $\{16,17\}$, $\{17,23\}$, $\{18,24\}$, $\{19,21\}$, $\{19,23\}$, $\{20,22\}$, $\{20,24\}$, $\{21,24\}$, $\{22,23\}$\},

$\adfGm$:
\{$\{1,2\}$, $\{1,3\}$, $\{1,4\}$, $\{2,5\}$, $\{2,6\}$, $\{3,7\}$, $\{3,8\}$, $\{4,9\}$, $\{4,10\}$, $\{5,7\}$, $\{5,9\}$, $\{6,8\}$, $\{6,11\}$, $\{7,12\}$, $\{8,13\}$, $\{9,14\}$, $\{10,15\}$, $\{10,16\}$, $\{11,17\}$, $\{11,18\}$, $\{12,15\}$, $\{12,17\}$, $\{13,19\}$, $\{13,20\}$, $\{14,21\}$, $\{14,22\}$, $\{15,18\}$, $\{16,17\}$, $\{16,23\}$, $\{18,24\}$, $\{19,21\}$, $\{19,23\}$, $\{20,22\}$, $\{20,24\}$, $\{21,24\}$, $\{22,23\}$\},

$\adfGn$:
\{$\{1,2\}$, $\{1,3\}$, $\{1,4\}$, $\{2,5\}$, $\{2,6\}$, $\{3,7\}$, $\{3,8\}$, $\{4,9\}$, $\{4,10\}$, $\{5,7\}$, $\{5,9\}$, $\{6,11\}$, $\{6,12\}$, $\{7,10\}$, $\{8,11\}$, $\{8,13\}$, $\{9,13\}$, $\{10,14\}$, $\{11,15\}$, $\{12,14\}$, $\{12,16\}$, $\{13,17\}$, $\{14,18\}$, $\{15,17\}$, $\{15,19\}$, $\{16,20\}$, $\{16,21\}$, $\{17,22\}$, $\{18,23\}$, $\{18,24\}$, $\{19,20\}$, $\{19,23\}$, $\{20,24\}$, $\{21,22\}$, $\{21,23\}$, $\{22,24\}$\},

$\adfGo$:
\{$\{1,2\}$, $\{1,3\}$, $\{1,4\}$, $\{2,5\}$, $\{2,6\}$, $\{3,7\}$, $\{3,8\}$, $\{4,9\}$, $\{4,10\}$, $\{5,7\}$, $\{5,9\}$, $\{6,11\}$, $\{6,12\}$, $\{7,10\}$, $\{8,11\}$, $\{8,13\}$, $\{9,14\}$, $\{10,15\}$, $\{11,14\}$, $\{12,15\}$, $\{12,16\}$, $\{13,17\}$, $\{13,18\}$, $\{14,17\}$, $\{15,19\}$, $\{16,20\}$, $\{16,21\}$, $\{17,22\}$, $\{18,20\}$, $\{18,23\}$, $\{19,23\}$, $\{19,24\}$, $\{20,24\}$, $\{21,22\}$, $\{21,23\}$, $\{22,24\}$\},

$\adfGp$:
\{$\{1,2\}$, $\{1,3\}$, $\{1,4\}$, $\{2,5\}$, $\{2,6\}$, $\{3,7\}$, $\{3,8\}$, $\{4,9\}$, $\{4,10\}$, $\{5,7\}$, $\{5,9\}$, $\{6,11\}$, $\{6,12\}$, $\{7,10\}$, $\{8,11\}$, $\{8,13\}$, $\{9,14\}$, $\{10,15\}$, $\{11,16\}$, $\{12,15\}$, $\{12,17\}$, $\{13,14\}$, $\{13,18\}$, $\{14,16\}$, $\{15,19\}$, $\{16,20\}$, $\{17,21\}$, $\{17,22\}$, $\{18,21\}$, $\{18,23\}$, $\{19,23\}$, $\{19,24\}$, $\{20,22\}$, $\{20,24\}$, $\{21,24\}$, $\{22,23\}$\},

$\adfGq$:
\{$\{1,2\}$, $\{1,3\}$, $\{1,4\}$, $\{2,5\}$, $\{2,6\}$, $\{3,7\}$, $\{3,8\}$, $\{4,9\}$, $\{4,10\}$, $\{5,7\}$, $\{5,9\}$, $\{6,11\}$, $\{6,12\}$, $\{7,10\}$, $\{8,11\}$, $\{8,13\}$, $\{9,14\}$, $\{10,15\}$, $\{11,16\}$, $\{12,15\}$, $\{12,17\}$, $\{13,18\}$, $\{13,19\}$, $\{14,20\}$, $\{14,21\}$, $\{15,16\}$, $\{16,22\}$, $\{17,22\}$, $\{17,23\}$, $\{18,20\}$, $\{18,23\}$, $\{19,21\}$, $\{19,24\}$, $\{20,24\}$, $\{21,23\}$, $\{22,24\}$\},

$\adfGr$:
\{$\{1,2\}$, $\{1,3\}$, $\{1,4\}$, $\{2,5\}$, $\{2,6\}$, $\{3,7\}$, $\{3,8\}$, $\{4,9\}$, $\{4,10\}$, $\{5,7\}$, $\{5,9\}$, $\{6,11\}$, $\{6,12\}$, $\{7,10\}$, $\{8,11\}$, $\{8,13\}$, $\{9,14\}$, $\{10,15\}$, $\{11,16\}$, $\{12,15\}$, $\{12,17\}$, $\{13,18\}$, $\{13,19\}$, $\{14,20\}$, $\{14,21\}$, $\{15,22\}$, $\{16,17\}$, $\{16,22\}$, $\{17,23\}$, $\{18,20\}$, $\{18,23\}$, $\{19,21\}$, $\{19,24\}$, $\{20,24\}$, $\{21,23\}$, $\{22,24\}$\},

$\adfGs$:
\{$\{1,2\}$, $\{1,3\}$, $\{1,4\}$, $\{2,5\}$, $\{2,6\}$, $\{3,7\}$, $\{3,8\}$, $\{4,9\}$, $\{4,10\}$, $\{5,7\}$, $\{5,9\}$, $\{6,11\}$, $\{6,12\}$, $\{7,10\}$, $\{8,11\}$, $\{8,13\}$, $\{9,14\}$, $\{10,15\}$, $\{11,16\}$, $\{12,17\}$, $\{12,18\}$, $\{13,19\}$, $\{13,20\}$, $\{14,16\}$, $\{14,19\}$, $\{15,21\}$, $\{15,22\}$, $\{16,20\}$, $\{17,21\}$, $\{17,23\}$, $\{18,22\}$, $\{18,24\}$, $\{19,23\}$, $\{20,24\}$, $\{21,24\}$, $\{22,23\}$\},

$\adfGt$:
\{$\{1,2\}$, $\{1,3\}$, $\{1,4\}$, $\{2,5\}$, $\{2,6\}$, $\{3,7\}$, $\{3,8\}$, $\{4,9\}$, $\{4,10\}$, $\{5,7\}$, $\{5,9\}$, $\{6,11\}$, $\{6,12\}$, $\{7,10\}$, $\{8,13\}$, $\{8,14\}$, $\{9,13\}$, $\{10,15\}$, $\{11,14\}$, $\{11,15\}$, $\{12,16\}$, $\{12,17\}$, $\{13,18\}$, $\{14,19\}$, $\{15,16\}$, $\{16,20\}$, $\{17,21\}$, $\{17,22\}$, $\{18,21\}$, $\{18,23\}$, $\{19,22\}$, $\{19,24\}$, $\{20,23\}$, $\{20,24\}$, $\{21,24\}$, $\{22,23\}$\},

$\adfGu$:
\{$\{1,2\}$, $\{1,3\}$, $\{1,4\}$, $\{2,5\}$, $\{2,6\}$, $\{3,7\}$, $\{3,8\}$, $\{4,9\}$, $\{4,10\}$, $\{5,7\}$, $\{5,9\}$, $\{6,11\}$, $\{6,12\}$, $\{7,10\}$, $\{8,13\}$, $\{8,14\}$, $\{9,13\}$, $\{10,15\}$, $\{11,14\}$, $\{11,15\}$, $\{12,16\}$, $\{12,17\}$, $\{13,18\}$, $\{14,19\}$, $\{15,20\}$, $\{16,21\}$, $\{16,22\}$, $\{17,23\}$, $\{17,24\}$, $\{18,19\}$, $\{18,21\}$, $\{19,23\}$, $\{20,22\}$, $\{20,24\}$, $\{21,24\}$, $\{22,23\}$\},

$\adfGv$:
\{$\{1,2\}$, $\{1,3\}$, $\{1,4\}$, $\{2,5\}$, $\{2,6\}$, $\{3,7\}$, $\{3,8\}$, $\{4,9\}$, $\{4,10\}$, $\{5,7\}$, $\{5,9\}$, $\{6,11\}$, $\{6,12\}$, $\{7,10\}$, $\{8,13\}$, $\{8,14\}$, $\{9,13\}$, $\{10,15\}$, $\{11,14\}$, $\{11,16\}$, $\{12,15\}$, $\{12,17\}$, $\{13,18\}$, $\{14,19\}$, $\{15,20\}$, $\{16,21\}$, $\{16,22\}$, $\{17,21\}$, $\{17,23\}$, $\{18,19\}$, $\{18,23\}$, $\{19,24\}$, $\{20,22\}$, $\{20,24\}$, $\{21,24\}$, $\{22,23\}$\},

$\adfGw$:
\{$\{1,2\}$, $\{1,3\}$, $\{1,4\}$, $\{2,5\}$, $\{2,6\}$, $\{3,7\}$, $\{3,8\}$, $\{4,9\}$, $\{4,10\}$, $\{5,7\}$, $\{5,9\}$, $\{6,11\}$, $\{6,12\}$, $\{7,10\}$, $\{8,13\}$, $\{8,14\}$, $\{9,13\}$, $\{10,15\}$, $\{11,14\}$, $\{11,16\}$, $\{12,15\}$, $\{12,17\}$, $\{13,18\}$, $\{14,19\}$, $\{15,20\}$, $\{16,21\}$, $\{16,22\}$, $\{17,21\}$, $\{17,23\}$, $\{18,19\}$, $\{18,24\}$, $\{19,23\}$, $\{20,22\}$, $\{20,24\}$, $\{21,24\}$, $\{22,23\}$\},

$\adfGx$:
\{$\{1,2\}$, $\{1,3\}$, $\{1,4\}$, $\{2,5\}$, $\{2,6\}$, $\{3,7\}$, $\{3,8\}$, $\{4,9\}$, $\{4,10\}$, $\{5,7\}$, $\{5,9\}$, $\{6,11\}$, $\{6,12\}$, $\{7,10\}$, $\{8,13\}$, $\{8,14\}$, $\{9,13\}$, $\{10,15\}$, $\{11,14\}$, $\{11,16\}$, $\{12,17\}$, $\{12,18\}$, $\{13,16\}$, $\{14,19\}$, $\{15,20\}$, $\{15,21\}$, $\{16,22\}$, $\{17,20\}$, $\{17,23\}$, $\{18,21\}$, $\{18,24\}$, $\{19,22\}$, $\{19,23\}$, $\{20,24\}$, $\{21,23\}$, $\{22,24\}$\},

$\adfGy$:
\{$\{1,2\}$, $\{1,3\}$, $\{1,4\}$, $\{2,5\}$, $\{2,6\}$, $\{3,7\}$, $\{3,8\}$, $\{4,9\}$, $\{4,10\}$, $\{5,7\}$, $\{5,9\}$, $\{6,11\}$, $\{6,12\}$, $\{7,10\}$, $\{8,13\}$, $\{8,14\}$, $\{9,13\}$, $\{10,15\}$, $\{11,15\}$, $\{11,16\}$, $\{12,17\}$, $\{12,18\}$, $\{13,19\}$, $\{14,16\}$, $\{14,20\}$, $\{15,21\}$, $\{16,19\}$, $\{17,20\}$, $\{17,22\}$, $\{18,23\}$, $\{18,24\}$, $\{19,23\}$, $\{20,24\}$, $\{21,22\}$, $\{21,24\}$, $\{22,23\}$\},

$\adfGz$:
\{$\{1,2\}$, $\{1,3\}$, $\{1,4\}$, $\{2,5\}$, $\{2,6\}$, $\{3,7\}$, $\{3,8\}$, $\{4,9\}$, $\{4,10\}$, $\{5,7\}$, $\{5,9\}$, $\{6,11\}$, $\{6,12\}$, $\{7,10\}$, $\{8,13\}$, $\{8,14\}$, $\{9,13\}$, $\{10,15\}$, $\{11,15\}$, $\{11,16\}$, $\{12,17\}$, $\{12,18\}$, $\{13,19\}$, $\{14,16\}$, $\{14,20\}$, $\{15,21\}$, $\{16,22\}$, $\{17,19\}$, $\{17,23\}$, $\{18,22\}$, $\{18,24\}$, $\{19,20\}$, $\{20,24\}$, $\{21,23\}$, $\{21,24\}$, $\{22,23\}$\},

$\adfGA$:
\{$\{1,2\}$, $\{1,3\}$, $\{1,4\}$, $\{2,5\}$, $\{2,6\}$, $\{3,7\}$, $\{3,8\}$, $\{4,9\}$, $\{4,10\}$, $\{5,7\}$, $\{5,9\}$, $\{6,11\}$, $\{6,12\}$, $\{7,10\}$, $\{8,13\}$, $\{8,14\}$, $\{9,13\}$, $\{10,15\}$, $\{11,15\}$, $\{11,16\}$, $\{12,17\}$, $\{12,18\}$, $\{13,19\}$, $\{14,17\}$, $\{14,20\}$, $\{15,21\}$, $\{16,20\}$, $\{16,22\}$, $\{17,19\}$, $\{18,22\}$, $\{18,23\}$, $\{19,24\}$, $\{20,23\}$, $\{21,23\}$, $\{21,24\}$, $\{22,24\}$\},

$\adfGB$:
\{$\{1,2\}$, $\{1,3\}$, $\{1,4\}$, $\{2,5\}$, $\{2,6\}$, $\{3,7\}$, $\{3,8\}$, $\{4,9\}$, $\{4,10\}$, $\{5,7\}$, $\{5,9\}$, $\{6,11\}$, $\{6,12\}$, $\{7,10\}$, $\{8,13\}$, $\{8,14\}$, $\{9,13\}$, $\{10,15\}$, $\{11,15\}$, $\{11,16\}$, $\{12,17\}$, $\{12,18\}$, $\{13,19\}$, $\{14,17\}$, $\{14,20\}$, $\{15,21\}$, $\{16,20\}$, $\{16,22\}$, $\{17,23\}$, $\{18,22\}$, $\{18,24\}$, $\{19,20\}$, $\{19,24\}$, $\{21,23\}$, $\{21,24\}$, $\{22,23\}$\},

$\adfGC$:
\{$\{1,2\}$, $\{1,3\}$, $\{1,4\}$, $\{2,5\}$, $\{2,6\}$, $\{3,7\}$, $\{3,8\}$, $\{4,9\}$, $\{4,10\}$, $\{5,7\}$, $\{5,9\}$, $\{6,11\}$, $\{6,12\}$, $\{7,10\}$, $\{8,13\}$, $\{8,14\}$, $\{9,15\}$, $\{10,16\}$, $\{11,13\}$, $\{11,16\}$, $\{12,17\}$, $\{12,18\}$, $\{13,17\}$, $\{14,19\}$, $\{14,20\}$, $\{15,21\}$, $\{15,22\}$, $\{16,23\}$, $\{17,23\}$, $\{18,19\}$, $\{18,21\}$, $\{19,22\}$, $\{20,21\}$, $\{20,24\}$, $\{22,24\}$, $\{23,24\}$\},

$\adfGD$:
\{$\{1,2\}$, $\{1,3\}$, $\{1,4\}$, $\{2,5\}$, $\{2,6\}$, $\{3,7\}$, $\{3,8\}$, $\{4,9\}$, $\{4,10\}$, $\{5,7\}$, $\{5,9\}$, $\{6,11\}$, $\{6,12\}$, $\{7,10\}$, $\{8,13\}$, $\{8,14\}$, $\{9,15\}$, $\{10,16\}$, $\{11,13\}$, $\{11,16\}$, $\{12,17\}$, $\{12,18\}$, $\{13,19\}$, $\{14,15\}$, $\{14,20\}$, $\{15,19\}$, $\{16,21\}$, $\{17,20\}$, $\{17,22\}$, $\{18,23\}$, $\{18,24\}$, $\{19,23\}$, $\{20,24\}$, $\{21,22\}$, $\{21,24\}$, $\{22,23\}$\},

$\adfGE$:
\{$\{1,2\}$, $\{1,3\}$, $\{1,4\}$, $\{2,5\}$, $\{2,6\}$, $\{3,7\}$, $\{3,8\}$, $\{4,9\}$, $\{4,10\}$, $\{5,7\}$, $\{5,9\}$, $\{6,11\}$, $\{6,12\}$, $\{7,10\}$, $\{8,13\}$, $\{8,14\}$, $\{9,15\}$, $\{10,16\}$, $\{11,13\}$, $\{11,16\}$, $\{12,17\}$, $\{12,18\}$, $\{13,19\}$, $\{14,15\}$, $\{14,20\}$, $\{15,21\}$, $\{16,17\}$, $\{17,22\}$, $\{18,20\}$, $\{18,23\}$, $\{19,23\}$, $\{19,24\}$, $\{20,24\}$, $\{21,22\}$, $\{21,23\}$, $\{22,24\}$\},

$\adfGF$:
\{$\{1,2\}$, $\{1,3\}$, $\{1,4\}$, $\{2,5\}$, $\{2,6\}$, $\{3,7\}$, $\{3,8\}$, $\{4,9\}$, $\{4,10\}$, $\{5,7\}$, $\{5,9\}$, $\{6,11\}$, $\{6,12\}$, $\{7,10\}$, $\{8,13\}$, $\{8,14\}$, $\{9,15\}$, $\{10,16\}$, $\{11,13\}$, $\{11,17\}$, $\{12,16\}$, $\{12,18\}$, $\{13,19\}$, $\{14,15\}$, $\{14,20\}$, $\{15,19\}$, $\{16,21\}$, $\{17,22\}$, $\{17,23\}$, $\{18,20\}$, $\{18,22\}$, $\{19,24\}$, $\{20,23\}$, $\{21,23\}$, $\{21,24\}$, $\{22,24\}$\},

$\adfGG$:
\{$\{1,2\}$, $\{1,3\}$, $\{1,4\}$, $\{2,5\}$, $\{2,6\}$, $\{3,7\}$, $\{3,8\}$, $\{4,9\}$, $\{4,10\}$, $\{5,7\}$, $\{5,9\}$, $\{6,11\}$, $\{6,12\}$, $\{7,10\}$, $\{8,13\}$, $\{8,14\}$, $\{9,15\}$, $\{10,16\}$, $\{11,13\}$, $\{11,17\}$, $\{12,16\}$, $\{12,18\}$, $\{13,19\}$, $\{14,15\}$, $\{14,20\}$, $\{15,19\}$, $\{16,21\}$, $\{17,22\}$, $\{17,23\}$, $\{18,22\}$, $\{18,24\}$, $\{19,24\}$, $\{20,21\}$, $\{20,22\}$, $\{21,23\}$, $\{23,24\}$\},

$\adfGH$:
\{$\{1,2\}$, $\{1,3\}$, $\{1,4\}$, $\{2,5\}$, $\{2,6\}$, $\{3,7\}$, $\{3,8\}$, $\{4,9\}$, $\{4,10\}$, $\{5,7\}$, $\{5,9\}$, $\{6,11\}$, $\{6,12\}$, $\{7,10\}$, $\{8,13\}$, $\{8,14\}$, $\{9,15\}$, $\{10,16\}$, $\{11,13\}$, $\{11,17\}$, $\{12,18\}$, $\{12,19\}$, $\{13,18\}$, $\{14,20\}$, $\{14,21\}$, $\{15,22\}$, $\{15,23\}$, $\{16,17\}$, $\{16,18\}$, $\{17,24\}$, $\{19,20\}$, $\{19,22\}$, $\{20,23\}$, $\{21,22\}$, $\{21,24\}$, $\{23,24\}$\},

$\adfGI$:
\{$\{1,2\}$, $\{1,3\}$, $\{1,4\}$, $\{2,5\}$, $\{2,6\}$, $\{3,7\}$, $\{3,8\}$, $\{4,9\}$, $\{4,10\}$, $\{5,7\}$, $\{5,9\}$, $\{6,11\}$, $\{6,12\}$, $\{7,10\}$, $\{8,13\}$, $\{8,14\}$, $\{9,15\}$, $\{10,16\}$, $\{11,16\}$, $\{11,17\}$, $\{12,18\}$, $\{12,19\}$, $\{13,15\}$, $\{13,17\}$, $\{14,18\}$, $\{14,20\}$, $\{15,20\}$, $\{16,21\}$, $\{17,22\}$, $\{18,23\}$, $\{19,22\}$, $\{19,24\}$, $\{20,24\}$, $\{21,23\}$, $\{21,24\}$, $\{22,23\}$\},

$\adfGJ$:
\{$\{1,2\}$, $\{1,3\}$, $\{1,4\}$, $\{2,5\}$, $\{2,6\}$, $\{3,7\}$, $\{3,8\}$, $\{4,9\}$, $\{4,10\}$, $\{5,7\}$, $\{5,9\}$, $\{6,11\}$, $\{6,12\}$, $\{7,10\}$, $\{8,13\}$, $\{8,14\}$, $\{9,15\}$, $\{10,16\}$, $\{11,16\}$, $\{11,17\}$, $\{12,18\}$, $\{12,19\}$, $\{13,15\}$, $\{13,18\}$, $\{14,17\}$, $\{14,20\}$, $\{15,20\}$, $\{16,21\}$, $\{17,22\}$, $\{18,23\}$, $\{19,22\}$, $\{19,24\}$, $\{20,24\}$, $\{21,23\}$, $\{21,24\}$, $\{22,23\}$\},

$\adfGK$:
\{$\{1,2\}$, $\{1,3\}$, $\{1,4\}$, $\{2,5\}$, $\{2,6\}$, $\{3,7\}$, $\{3,8\}$, $\{4,9\}$, $\{4,10\}$, $\{5,7\}$, $\{5,9\}$, $\{6,11\}$, $\{6,12\}$, $\{7,13\}$, $\{8,14\}$, $\{8,15\}$, $\{9,14\}$, $\{10,16\}$, $\{10,17\}$, $\{11,16\}$, $\{11,18\}$, $\{12,19\}$, $\{12,20\}$, $\{13,18\}$, $\{13,21\}$, $\{14,22\}$, $\{15,19\}$, $\{15,23\}$, $\{16,21\}$, $\{17,18\}$, $\{17,23\}$, $\{19,24\}$, $\{20,22\}$, $\{20,23\}$, $\{21,24\}$, $\{22,24\}$\}
and

$\adfGL$:
\{$\{1,2\}$, $\{1,3\}$, $\{1,4\}$, $\{2,5\}$, $\{2,6\}$, $\{3,7\}$, $\{3,8\}$, $\{4,9\}$, $\{4,10\}$, $\{5,7\}$, $\{5,9\}$, $\{6,11\}$, $\{6,12\}$, $\{7,13\}$, $\{8,14\}$, $\{8,15\}$, $\{9,14\}$, $\{10,16\}$, $\{10,17\}$, $\{11,16\}$, $\{11,18\}$, $\{12,19\}$, $\{12,20\}$, $\{13,18\}$, $\{13,21\}$, $\{14,22\}$, $\{15,19\}$, $\{15,23\}$, $\{16,21\}$, $\{17,18\}$, $\{17,24\}$, $\{19,24\}$, $\{20,22\}$, $\{20,23\}$, $\{21,23\}$, $\{22,24\}$\}.


{\lemma \label{lem:Snark24 136}
There exist a design of order $136$ for each of the thirty-eight non-trivial 24-vertex snarks.}\\


\noindent\textbf{Proof.}
Let the vertex set of $K_{136}$ be $Z_{135} \cup \{\infty\}$. The decompositions consist of

$(\infty,1,96,44,95,99,21,54,102,3,90,11,\adfsplit 130,94,5,133,120,31,2,25,77,116,98,9)_{\adfGa}$,

$(33,35,88,18,101,9,132,34,58,126,42,114,\adfsplit 120,39,72,119,107,28,70,105,12,133,56,106)_{\adfGa}$,

$(29,13,28,12,27,84,31,37,56,95,46,133,\adfsplit 122,10,114,121,130,60,3,4,85,14,41,119)_{\adfGa}$,

$(120,47,86,95,25,14,80,29,104,58,26,7,\adfsplit 76,28,97,121,41,33,18,46,70,63,133,4)_{\adfGa}$,

$(132,70,3,96,55,5,24,60,59,130,26,56,\adfsplit 37,76,98,16,111,102,109,22,79,47,63,74)_{\adfGa}$,

\adfSgap
$(97,84,54,95,83,35,119,77,41,9,12,65,\adfsplit 57,120,55,42,126,50,53,62,123,63,98,51)_{\adfGa}$,

$(2,102,26,48,105,20,56,50,59,114,22,45,\adfsplit 80,24,9,98,116,77,29,11,74,62,113,122)_{\adfGa}$,

\adfLgap
$(\infty,71,112,54,34,92,74,91,27,13,99,46,\adfsplit 107,115,123,121,97,25,19,47,31,51,106,14)_{\adfGb}$,

$(116,84,89,26,70,52,113,75,125,123,122,82,\adfsplit 132,81,95,43,78,32,97,13,62,20,83,115)_{\adfGb}$,

$(30,93,125,4,59,85,79,15,78,23,0,81,\adfsplit 126,106,46,36,82,72,99,103,100,77,43,41)_{\adfGb}$,

$(89,126,98,31,28,105,95,10,22,16,110,77,\adfsplit 61,86,107,106,133,60,14,96,101,3,68,79)_{\adfGb}$,

$(98,49,63,78,25,76,67,60,117,30,19,70,\adfsplit 22,58,20,6,105,120,59,126,111,83,3,38)_{\adfGb}$,

\adfSgap
$(131,57,125,42,87,59,96,119,54,29,123,21,\adfsplit 51,84,12,62,120,38,45,128,56,53,36,50)_{\adfGb}$,

$(44,90,56,122,26,77,41,27,111,62,60,117,\adfsplit 87,98,50,48,2,83,21,108,80,78,9,14)_{\adfGb}$,

\adfLgap
$(\infty,31,33,86,82,70,32,55,90,95,1,21,\adfsplit 115,68,16,2,13,91,19,73,39,111,45,134)_{\adfGc}$,

$(107,50,31,72,78,116,33,48,42,12,91,30,\adfsplit 90,112,59,18,13,80,122,32,17,98,11,58)_{\adfGc}$,

$(38,82,85,67,95,128,78,112,50,97,116,96,\adfsplit 132,13,101,55,36,16,93,83,104,12,76,66)_{\adfGc}$,

$(92,68,119,3,76,104,105,60,44,85,54,32,\adfsplit 125,78,72,45,8,86,71,12,34,65,33,102)_{\adfGc}$,

$(29,59,126,91,44,94,128,42,77,80,103,70,\adfsplit 119,24,9,114,111,33,10,0,92,130,32,58)_{\adfGc}$,

\adfSgap
$(35,88,118,1,21,24,97,73,69,100,112,93,\adfsplit 36,117,22,45,82,85,70,127,34,60,124,106)_{\adfGc}$,

$(1,60,36,79,97,64,40,10,131,25,111,30,\adfsplit 31,49,15,76,67,90,101,75,115,22,112,19)_{\adfGc}$,

\adfLgap
$(\infty,75,70,68,60,122,6,33,35,109,71,48,\adfsplit 91,111,13,65,22,63,121,118,42,7,47,126)_{\adfGd}$,

$(7,83,59,128,76,43,22,102,14,71,21,124,\adfsplit 8,129,72,98,5,132,52,125,39,110,95,97)_{\adfGd}$,

$(79,81,129,16,52,40,29,109,69,98,118,56,\adfsplit 67,47,37,87,48,63,116,44,21,26,50,131)_{\adfGd}$,

$(108,81,122,131,105,11,91,50,95,67,20,82,\adfsplit 117,26,17,75,128,73,89,107,106,66,96,45)_{\adfGd}$,

$(62,91,38,64,95,83,20,73,126,65,117,3,\adfsplit 17,90,70,30,33,85,130,109,96,81,58,69)_{\adfGd}$,

\adfSgap
$(127,95,52,91,24,121,6,70,94,79,4,66,\adfsplit 38,51,115,19,1,102,64,84,18,67,0,21)_{\adfGd}$,

$(7,4,121,39,124,64,109,85,127,99,13,25,\adfsplit 97,78,24,45,117,118,71,15,0,102,70,103)_{\adfGd}$,

\adfLgap
$(\infty,93,53,76,43,0,113,49,24,2,34,54,\adfsplit 64,77,6,82,28,21,125,117,96,122,87,71)_{\adfGe}$,

$(2,105,34,106,112,78,114,68,10,24,79,63,\adfsplit 23,72,52,82,45,32,64,39,116,1,18,43)_{\adfGe}$,

$(125,24,118,77,86,16,92,6,84,78,131,35,\adfsplit 76,23,18,4,117,52,95,60,130,80,29,128)_{\adfGe}$,

$(46,112,67,127,26,106,114,131,128,110,61,101,\adfsplit 59,43,91,14,16,19,40,84,111,7,86,93)_{\adfGe}$,

$(120,57,101,16,0,52,106,115,109,133,107,8,\adfsplit 55,87,76,104,82,75,39,80,93,38,119,72)_{\adfGe}$,

\adfSgap
$(76,98,51,38,71,83,45,68,50,65,104,93,\adfsplit 106,81,44,82,128,117,21,57,111,105,6,87)_{\adfGe}$,

$(5,17,33,99,32,29,59,123,9,71,11,41,\adfsplit 129,20,53,42,18,125,36,90,51,120,21,84)_{\adfGe}$,

\adfLgap
$(\infty,64,33,83,116,34,85,45,111,82,28,109,\adfsplit 37,0,50,35,86,110,70,68,115,40,69,56)_{\adfGf}$,

$(11,89,129,50,90,60,94,3,115,122,39,40,\adfsplit 88,92,102,133,33,101,66,32,46,100,112,99)_{\adfGf}$,

$(55,11,50,92,35,24,53,101,132,118,64,47,\adfsplit 43,122,67,96,108,23,111,103,134,4,37,15)_{\adfGf}$,

$(80,125,14,101,50,36,69,64,132,8,13,129,\adfsplit 130,115,15,68,112,106,74,58,34,71,44,96)_{\adfGf}$,

$(8,25,72,60,97,39,1,115,63,76,74,27,\adfsplit 23,7,18,9,96,55,14,31,51,128,130,36)_{\adfGf}$,

\adfSgap
$(67,46,92,32,38,25,126,4,123,80,50,30,\adfsplit 131,84,65,0,39,78,21,66,69,72,105,117)_{\adfGf}$,

$(5,52,92,129,17,44,127,0,60,95,33,80,\adfsplit 108,54,117,114,78,96,30,45,111,56,8,3)_{\adfGf}$,

\adfLgap
$(\infty,123,16,92,100,30,29,69,132,57,27,67,\adfsplit 121,18,95,3,77,60,90,23,61,2,66,12)_{\adfGg}$,

$(131,5,95,72,13,10,56,49,57,3,119,48,\adfsplit 66,11,89,32,86,16,46,74,22,114,88,4)_{\adfGg}$,

$(80,6,101,23,62,35,91,82,0,134,68,92,\adfsplit 48,26,119,123,61,76,4,53,50,130,75,20)_{\adfGg}$,

$(34,95,67,83,91,51,85,88,42,89,65,39,\adfsplit 115,61,82,5,118,53,133,57,71,69,19,60)_{\adfGg}$,

$(89,24,101,121,103,30,83,11,15,91,62,1,\adfsplit 65,98,31,25,93,37,28,21,22,102,33,132)_{\adfGg}$,

\adfSgap
$(98,101,95,4,84,104,48,10,3,91,87,132,\adfsplit 61,111,69,51,109,115,9,93,36,0,58,117)_{\adfGg}$,

$(1,55,27,3,54,75,121,7,130,105,49,70,\adfsplit 66,80,103,69,25,117,33,63,16,42,110,109)_{\adfGg}$,

\adfLgap
$(\infty,40,32,105,117,59,76,93,77,68,104,37,\adfsplit 5,85,128,18,4,34,79,110,65,64,51,126)_{\adfGh}$,

$(57,36,114,81,4,55,19,66,112,10,46,17,\adfsplit 40,31,113,84,125,134,86,117,129,95,9,7)_{\adfGh}$,

$(20,37,120,112,122,129,56,3,74,0,6,28,\adfsplit 75,30,52,73,72,79,71,10,83,53,63,110)_{\adfGh}$,

$(54,18,123,88,74,128,81,69,4,8,51,110,\adfsplit 17,40,62,5,28,103,120,85,50,57,27,119)_{\adfGh}$,

$(20,31,82,26,105,7,103,87,52,111,0,3,\adfsplit 36,39,6,70,48,121,28,42,41,86,118,50)_{\adfGh}$,

\adfSgap
$(90,129,51,77,116,33,119,20,74,25,59,82,\adfsplit 83,89,19,22,94,67,88,73,104,124,62,31)_{\adfGh}$,

$(2,85,95,11,28,89,100,98,32,7,110,130,\adfsplit 35,41,70,26,124,9,66,113,92,82,77,134)_{\adfGh}$,

\adfLgap
$(\infty,23,123,79,87,12,80,58,38,90,134,109,\adfsplit 42,15,50,122,29,37,132,121,82,104,99,1)_{\adfGi}$,

$(85,50,125,73,47,78,16,96,133,120,25,48,\adfsplit 59,71,60,83,12,23,84,69,75,22,40,91)_{\adfGi}$,

$(3,128,52,36,83,27,58,7,50,80,61,60,\adfsplit 53,18,6,59,120,133,11,35,69,126,12,107)_{\adfGi}$,

$(84,48,90,134,87,45,16,21,127,46,53,7,\adfsplit 59,64,56,31,129,119,70,8,43,0,6,52)_{\adfGi}$,

$(36,83,90,94,78,35,133,34,37,3,117,108,\adfsplit 24,40,45,109,86,107,104,116,110,22,43,38)_{\adfGi}$,

\adfSgap
$(86,67,60,58,119,71,34,44,22,38,19,53,\adfsplit 18,56,65,104,14,2,127,109,107,70,10,28)_{\adfGi}$,

$(1,103,34,31,50,80,64,133,73,57,52,83,\adfsplit 125,107,41,106,56,40,130,121,122,68,70,88)_{\adfGi}$,

\adfLgap
$(\infty,17,37,72,38,19,104,108,74,117,122,48,\adfsplit 119,93,106,92,64,42,68,15,107,88,44,84)_{\adfGj}$,

$(31,55,121,113,26,51,22,108,27,117,48,53,\adfsplit 95,10,88,41,103,40,46,47,24,3,36,101)_{\adfGj}$,

$(132,71,12,108,9,85,106,120,98,134,93,45,\adfsplit 61,89,63,6,59,94,56,127,119,100,126,104)_{\adfGj}$,

$(16,101,92,23,72,68,35,60,103,66,65,57,\adfsplit 32,69,117,78,30,95,58,122,106,4,111,76)_{\adfGj}$,

$(133,21,48,16,0,93,9,49,36,112,53,25,\adfsplit 105,20,28,114,59,119,125,124,113,91,37,72)_{\adfGj}$,

\adfSgap
$(3,106,28,46,11,81,19,49,103,58,124,94,\adfsplit 112,109,115,67,61,29,97,52,37,14,131,89)_{\adfGj}$,

$(1,43,16,44,35,31,73,68,62,4,83,52,\adfsplit 25,127,46,38,56,94,104,37,112,71,131,69)_{\adfGj}$,

\adfLgap
$(\infty,8,99,25,127,54,124,24,9,38,74,122,\adfsplit 2,62,53,66,101,47,79,44,5,109,51,4)_{\adfGk}$,

$(20,54,93,44,43,48,125,32,21,124,100,33,\adfsplit 7,112,88,39,131,83,92,41,1,130,109,82)_{\adfGk}$,

$(69,65,31,28,129,128,106,38,57,97,13,132,\adfsplit 29,113,35,107,55,119,58,34,110,4,60,102)_{\adfGk}$,

$(45,91,72,62,19,105,64,13,76,48,66,8,\adfsplit 89,113,102,36,67,100,97,70,1,9,46,63)_{\adfGk}$,

$(97,57,49,86,62,52,34,58,80,82,79,69,\adfsplit 78,105,48,63,132,87,26,127,131,114,66,95)_{\adfGk}$,

\adfSgap
$(130,27,95,62,29,60,105,28,97,33,11,116,\adfsplit 128,129,117,84,44,77,93,63,51,74,35,120)_{\adfGk}$,

$(0,59,77,78,126,108,26,66,71,80,39,96,\adfsplit 131,81,29,12,6,41,47,92,48,15,14,99)_{\adfGk}$,

\adfLgap
$(\infty,6,19,47,82,37,55,86,93,98,40,88,\adfsplit 26,10,85,73,56,35,42,90,71,87,22,21)_{\adfGl}$,

$(130,134,50,119,3,103,76,98,7,28,84,133,\adfsplit 71,104,60,58,116,99,36,111,23,124,17,57)_{\adfGl}$,

$(58,29,80,131,59,88,47,72,87,107,24,120,\adfsplit 90,11,122,67,81,130,127,8,56,43,55,14)_{\adfGl}$,

$(119,23,109,42,120,95,84,25,126,35,113,125,\adfsplit 64,38,96,127,0,116,102,16,15,39,22,105)_{\adfGl}$,

$(52,98,87,59,101,32,44,115,22,106,42,96,\adfsplit 27,51,15,65,9,41,99,57,111,72,74,62)_{\adfGl}$,

\adfSgap
$(8,75,7,22,84,70,16,9,133,64,56,15,\adfsplit 18,19,55,82,123,76,61,86,21,85,67,100)_{\adfGl}$,

$(1,6,129,96,132,13,130,79,46,4,34,49,\adfsplit 25,115,90,9,16,40,66,67,109,73,127,108)_{\adfGl}$,

\adfLgap
$(\infty,125,103,69,70,64,82,106,24,20,81,123,\adfsplit 47,66,77,33,15,49,115,40,45,16,61,75)_{\adfGm}$,

$(109,91,70,124,77,55,118,117,71,20,129,49,\adfsplit 45,110,112,121,98,35,99,68,114,130,0,43)_{\adfGm}$,

$(8,48,10,64,80,128,14,78,53,116,37,90,\adfsplit 85,38,51,106,66,39,79,84,60,41,1,32)_{\adfGm}$,

$(14,31,13,51,116,54,69,122,42,113,68,103,\adfsplit 58,119,12,96,93,23,118,7,89,129,134,87)_{\adfGm}$,

$(54,102,48,1,122,127,52,105,103,93,29,3,\adfsplit 104,114,34,73,38,7,60,51,76,35,86,84)_{\adfGm}$,

\adfSgap
$(58,21,104,49,72,96,123,125,132,11,12,25,\adfsplit 117,5,113,92,128,20,34,98,43,122,80,38)_{\adfGm}$,

$(128,14,41,40,116,62,53,122,35,8,20,17,\adfsplit 59,83,119,46,134,95,23,91,11,54,129,100)_{\adfGm}$,

\adfLgap
$(\infty,124,62,102,110,112,34,70,15,93,86,103,\adfsplit 89,12,7,10,98,90,8,75,37,117,104,97)_{\adfGn}$,

$(54,13,19,133,108,124,93,88,85,118,38,28,\adfsplit 8,128,34,5,33,94,4,77,114,16,134,42)_{\adfGn}$,

$(54,120,1,72,56,18,76,22,45,78,85,77,\adfsplit 53,64,49,87,36,67,113,106,74,28,24,117)_{\adfGn}$,

$(85,44,52,75,127,132,95,74,24,13,9,134,\adfsplit 96,64,62,1,59,26,33,91,48,89,69,116)_{\adfGn}$,

$(56,67,0,24,13,96,132,4,120,101,10,127,\adfsplit 65,63,28,70,66,11,134,53,36,126,98,102)_{\adfGn}$,

\adfSgap
$(10,98,47,15,131,29,71,68,8,38,45,24,\adfsplit 125,20,3,66,110,134,26,32,59,22,41,27)_{\adfGn}$,

$(5,29,26,52,41,104,83,77,57,89,20,108,\adfsplit 53,101,71,78,11,68,66,36,62,131,96,2)_{\adfGn}$,

\adfLgap
$(\infty,112,90,29,49,24,116,74,127,134,5,25,\adfsplit 34,124,1,77,2,119,56,125,89,13,78,132)_{\adfGo}$,

$(11,19,4,118,5,131,119,69,60,21,38,1,\adfsplit 9,51,54,92,71,61,66,13,126,48,134,94)_{\adfGo}$,

$(47,4,121,20,10,12,78,17,35,24,88,1,\adfsplit 8,114,98,18,108,96,28,118,90,80,40,77)_{\adfGo}$,

$(28,134,24,132,88,102,63,34,127,12,79,120,\adfsplit 49,2,76,78,133,16,117,58,112,103,43,77)_{\adfGo}$,

$(121,18,15,68,91,110,92,13,125,54,127,88,\adfsplit 130,108,53,51,83,9,48,96,81,11,12,75)_{\adfGo}$,

\adfSgap
$(97,61,48,42,121,1,72,88,66,125,81,76,\adfsplit 75,110,36,40,26,29,80,27,33,62,119,113)_{\adfGo}$,

$(1,100,123,96,93,51,41,77,11,5,27,53,\adfsplit 23,110,3,107,56,47,32,131,83,116,122,71)_{\adfGo}$,

\adfLgap
$(\infty,77,30,61,108,51,14,37,58,60,27,102,\adfsplit 0,117,46,126,91,28,103,74,17,70,15,118)_{\adfGp}$,

$(40,73,110,11,61,4,16,88,92,116,134,59,\adfsplit 69,21,106,71,120,36,35,103,14,28,52,83)_{\adfGp}$,

$(68,11,16,17,53,8,98,12,112,2,9,63,\adfsplit 124,110,88,15,93,30,61,77,23,95,62,105)_{\adfGp}$,

$(83,87,107,39,82,66,34,30,95,51,93,80,\adfsplit 17,27,118,44,74,8,126,28,16,124,106,111)_{\adfGp}$,

$(105,60,62,116,0,41,54,24,57,15,29,20,\adfsplit 106,80,76,117,120,36,53,48,82,16,78,46)_{\adfGp}$,

\adfSgap
$(29,128,56,1,85,68,36,116,27,7,41,32,\adfsplit 73,127,100,21,13,92,109,55,67,126,25,4)_{\adfGp}$,

$(5,97,112,115,48,43,106,121,22,134,21,73,\adfsplit 16,51,98,85,15,7,125,127,79,64,58,105)_{\adfGp}$,

\adfLgap
$(\infty,87,65,55,68,108,114,25,7,132,34,110,\adfsplit 120,100,128,75,31,89,117,78,2,99,4,48)_{\adfGq}$,

$(12,1,111,20,129,57,31,79,3,65,13,61,\adfsplit 15,52,81,97,23,91,124,37,47,74,120,99)_{\adfGq}$,

$(57,98,91,14,105,59,63,34,35,0,83,92,\adfsplit 87,93,120,97,18,37,32,39,1,8,17,72)_{\adfGq}$,

$(117,17,90,50,16,97,76,15,110,86,116,96,\adfsplit 85,78,13,38,109,120,66,83,122,80,130,92)_{\adfGq}$,

$(129,6,25,90,84,89,130,63,113,42,41,35,\adfsplit 18,46,48,107,120,13,65,1,109,103,112,71)_{\adfGq}$,

\adfSgap
$(23,16,40,106,59,127,122,58,130,95,110,68,\adfsplit 52,28,98,125,115,134,26,10,56,128,82,4)_{\adfGq}$,

$(5,112,31,20,77,1,74,107,47,133,44,101,\adfsplit 7,110,94,61,37,35,124,22,28,125,50,118)_{\adfGq}$,

\adfLgap
$(\infty,41,28,45,96,58,76,32,134,65,21,13,\adfsplit 48,0,57,61,62,25,67,125,6,15,101,16)_{\adfGr}$,

$(32,0,9,119,17,84,53,107,12,131,70,69,\adfsplit 40,79,56,104,42,35,47,50,112,59,124,126)_{\adfGr}$,

$(29,100,87,20,115,24,52,128,76,62,93,121,\adfsplit 81,66,132,48,79,83,16,126,45,1,57,123)_{\adfGr}$,

$(77,91,126,110,50,8,122,85,97,130,88,35,\adfsplit 117,114,13,24,58,7,108,106,52,96,118,74)_{\adfGr}$,

$(74,8,20,0,32,112,2,73,57,109,82,27,\adfsplit 16,96,51,39,75,37,9,10,63,69,56,64)_{\adfGr}$,

\adfSgap
$(18,34,32,71,77,36,16,7,115,10,65,119,\adfsplit 47,53,35,96,94,49,109,55,59,112,110,17)_{\adfGr}$,

$(1,41,3,17,118,20,74,56,102,68,6,26,\adfsplit 13,71,132,59,90,7,86,92,57,80,5,130)_{\adfGr}$,

\adfLgap
$(\infty,121,122,120,7,30,84,123,21,129,2,133,\adfsplit 124,44,11,67,82,116,113,58,48,91,4,63)_{\adfGs}$,

$(79,18,124,51,86,47,36,89,8,43,7,91,\adfsplit 16,60,132,2,134,117,10,115,67,53,108,68)_{\adfGs}$,

$(33,99,82,15,75,92,117,79,5,123,96,71,\adfsplit 89,7,98,101,104,34,30,27,63,62,46,47)_{\adfGs}$,

$(99,27,107,112,6,8,121,101,53,102,37,46,\adfsplit 22,14,49,97,61,7,103,34,29,31,9,63)_{\adfGs}$,

$(49,56,120,6,41,12,122,16,101,17,109,69,\adfsplit 124,11,89,8,47,73,20,50,71,77,28,23)_{\adfGs}$,

\adfSgap
$(17,68,30,63,81,119,0,126,33,32,87,132,\adfsplit 51,88,78,48,121,120,21,102,16,45,96,18)_{\adfGs}$,

$(4,99,106,72,40,96,73,81,10,75,93,15,\adfsplit 129,78,123,47,48,55,19,98,121,0,124,66)_{\adfGs}$,

\adfLgap
$(\infty,2,52,0,66,95,87,17,63,118,80,48,\adfsplit 53,31,78,89,21,115,75,132,125,74,85,42)_{\adfGt}$,

$(78,90,27,37,84,95,106,3,9,49,50,36,\adfsplit 17,29,0,10,61,30,75,121,94,113,8,2)_{\adfGt}$,

$(0,56,124,32,110,90,44,65,93,102,81,25,\adfsplit 40,26,1,128,91,88,6,99,119,41,97,92)_{\adfGt}$,

$(43,0,63,62,116,16,107,47,25,49,75,122,\adfsplit 88,52,67,13,58,65,15,64,103,132,4,45)_{\adfGt}$,

$(54,121,0,123,24,36,131,87,81,91,75,21,\adfsplit 82,50,39,25,44,109,49,23,92,111,124,95)_{\adfGt}$,

\adfSgap
$(96,67,56,7,88,103,49,74,95,28,83,100,\adfsplit 62,23,8,32,106,102,19,20,13,112,73,16)_{\adfGt}$,

$(5,45,134,62,91,16,52,116,55,119,86,115,\adfsplit 49,92,35,95,73,131,22,77,8,26,20,110)_{\adfGt}$,

\adfLgap
$(\infty,126,41,121,109,56,78,82,13,2,114,34,\adfsplit 127,124,113,54,49,33,3,9,91,55,42,95)_{\adfGu}$,

$(43,56,119,122,55,19,11,85,117,79,52,101,\adfsplit 60,133,32,23,115,36,45,75,77,89,72,123)_{\adfGu}$,

$(127,108,132,42,2,55,30,15,13,62,37,51,\adfsplit 85,69,75,112,74,20,9,14,35,76,24,56)_{\adfGu}$,

$(48,127,131,64,76,44,87,56,9,90,14,73,\adfsplit 37,78,46,35,75,15,68,20,60,107,12,39)_{\adfGu}$,

$(53,69,55,108,27,18,15,5,35,49,24,11,\adfsplit 44,93,73,14,45,13,70,43,88,114,47,83)_{\adfGu}$,

\adfSgap
$(90,86,55,116,61,22,67,65,71,133,29,57,\adfsplit 16,37,33,53,83,19,25,59,2,4,32,47)_{\adfGu}$,

$(2,112,35,131,109,95,43,10,29,31,46,50,\adfsplit 17,110,37,62,86,88,116,103,97,113,80,94)_{\adfGu}$,

\adfLgap
$(\infty,61,119,105,104,129,46,118,24,10,1,60,\adfsplit 42,123,84,124,50,93,125,109,26,4,17,39)_{\adfGv}$,

$(33,67,62,43,61,36,69,12,118,131,11,7,\adfsplit 89,80,90,37,55,6,52,101,59,87,88,98)_{\adfGv}$,

$(54,48,34,109,126,26,30,52,29,9,67,74,\adfsplit 133,93,32,43,125,61,105,47,127,68,108,41)_{\adfGv}$,

$(60,33,63,37,38,80,128,61,116,133,73,126,\adfsplit 72,70,88,52,0,121,42,24,20,69,99,50)_{\adfGv}$,

$(37,17,79,51,21,28,5,63,20,104,128,60,\adfsplit 121,30,123,93,133,55,131,3,50,45,119,40)_{\adfGv}$,

\adfSgap
$(48,108,104,91,16,33,132,131,83,53,93,76,\adfsplit 56,112,116,65,49,106,89,134,119,74,11,7)_{\adfGv}$,

$(2,14,119,125,32,22,100,1,44,50,38,104,\adfsplit 71,117,88,91,41,4,89,61,64,131,122,80)_{\adfGv}$,

\adfLgap
$(\infty,94,105,59,27,134,123,19,114,117,46,118,\adfsplit 47,8,121,64,22,131,14,65,36,29,6,1)_{\adfGw}$,

$(25,35,88,54,2,94,107,56,91,120,48,86,\adfsplit 129,39,30,123,57,106,93,83,82,79,60,112)_{\adfGw}$,

$(63,84,127,13,86,104,101,66,90,88,67,8,\adfsplit 73,118,112,69,68,72,98,22,111,103,129,64)_{\adfGw}$,

$(112,95,25,121,29,71,105,116,78,1,39,22,\adfsplit 73,104,134,77,49,40,52,51,119,98,108,61)_{\adfGw}$,

$(55,66,121,116,46,13,80,81,71,105,14,49,\adfsplit 62,86,127,82,56,112,134,113,29,114,84,32)_{\adfGw}$,

\adfSgap
$(127,57,54,65,114,120,3,66,30,111,84,58,\adfsplit 44,50,123,20,131,101,87,107,42,36,60,6)_{\adfGw}$,

$(3,101,47,131,60,0,104,54,53,90,57,15,\adfsplit 18,17,27,81,99,96,111,97,59,35,126,24)_{\adfGw}$,

\adfLgap
$(\infty,95,16,90,68,19,79,22,85,66,10,106,\adfsplit 62,53,49,101,32,15,30,52,84,47,120,12)_{\adfGx}$,

$(31,81,0,97,103,14,4,26,101,29,76,40,\adfsplit 59,78,48,46,67,125,134,112,91,83,105,3)_{\adfGx}$,

$(38,44,24,101,16,37,91,74,100,115,25,78,\adfsplit 108,99,59,7,90,29,127,132,87,39,57,41)_{\adfGx}$,

$(124,28,0,114,110,48,20,128,105,25,38,127,\adfsplit 106,98,6,81,70,2,130,5,43,97,123,26)_{\adfGx}$,

$(113,110,129,94,13,64,67,65,42,18,43,35,\adfsplit 100,101,56,85,19,83,48,74,87,2,11,17)_{\adfGx}$,

\adfSgap
$(13,62,36,101,53,102,21,120,123,60,81,79,\adfsplit 131,33,114,122,128,32,111,96,54,57,119,0)_{\adfGx}$,

$(0,6,11,108,80,57,72,51,102,126,9,30,\adfsplit 123,24,105,48,90,93,1,27,111,89,50,78)_{\adfGx}$,

\adfLgap
$(\infty,68,108,127,24,55,134,11,115,126,0,96,\adfsplit 45,33,114,118,121,28,73,2,72,21,97,6)_{\adfGy}$,

$(103,106,68,54,93,47,117,108,95,35,19,24,\adfsplit 3,126,100,5,44,30,133,73,134,97,4,88)_{\adfGy}$,

$(35,0,6,99,112,9,60,51,131,96,129,133,\adfsplit 12,115,83,81,72,53,11,5,120,44,50,62)_{\adfGy}$,

$(91,20,19,64,3,27,124,94,59,113,54,58,\adfsplit 133,47,134,101,51,17,18,108,49,67,16,6)_{\adfGy}$,

$(54,133,129,116,8,27,4,34,76,25,17,90,\adfsplit 30,86,134,2,43,109,22,81,26,74,31,131)_{\adfGy}$,

\adfSgap
$(96,34,64,133,131,66,28,97,71,101,4,103,\adfsplit 95,10,77,2,86,1,62,128,11,44,134,49)_{\adfGy}$,

$(2,79,92,77,31,71,134,80,133,122,47,13,\adfsplit 97,44,59,110,45,49,74,82,95,118,41,26)_{\adfGy}$,

\adfLgap
$(\infty,58,111,29,43,22,106,26,54,80,66,61,\adfsplit 1,74,28,118,30,57,119,50,16,59,46,65)_{\adfGz}$,

$(96,63,55,4,107,104,121,32,26,94,62,48,\adfsplit 67,87,17,36,13,38,69,60,28,39,23,75)_{\adfGz}$,

$(55,132,45,63,115,8,93,58,40,122,100,96,\adfsplit 128,86,47,15,90,114,60,106,43,80,61,49)_{\adfGz}$,

$(46,108,45,100,56,84,33,0,50,26,4,98,\adfsplit 49,63,130,88,36,3,70,103,101,95,122,23)_{\adfGz}$,

$(124,3,53,82,118,78,115,13,50,6,10,7,\adfsplit 59,61,34,42,41,44,43,5,108,80,113,9)_{\adfGz}$,

\adfSgap
$(106,69,99,98,38,122,41,71,66,59,18,130,\adfsplit 83,101,72,131,93,123,80,62,114,57,0,45)_{\adfGz}$,

$(3,84,99,42,26,63,131,38,23,119,44,24,\adfsplit 6,77,66,105,78,102,120,89,73,83,65,36)_{\adfGz}$,

\adfLgap
$(\infty,68,0,109,87,60,20,118,48,112,65,66,\adfsplit 127,25,37,9,92,98,126,21,26,114,14,32)_{\adfGA}$,

$(95,53,38,92,19,83,90,129,54,22,25,118,\adfsplit 79,80,86,9,66,21,126,56,23,111,49,10)_{\adfGA}$,

$(86,107,90,91,67,61,112,65,130,102,32,17,\adfsplit 13,28,45,64,76,40,29,123,21,50,69,70)_{\adfGA}$,

$(48,66,113,16,39,83,111,122,70,69,56,34,\adfsplit 24,103,25,1,126,107,61,9,52,78,22,85)_{\adfGA}$,

$(23,22,45,76,16,43,95,5,96,133,94,112,\adfsplit 105,129,116,130,124,97,126,128,14,122,24,83)_{\adfGA}$,

\adfSgap
$(115,119,6,129,101,60,21,57,53,83,2,122,\adfsplit 126,31,47,63,45,32,42,35,62,78,17,14)_{\adfGA}$,

$(2,56,66,30,33,41,132,81,107,129,128,23,\adfsplit 84,35,68,53,99,77,0,134,123,113,54,109)_{\adfGA}$,

\adfLgap
$(\infty,54,7,8,112,66,103,72,134,131,96,53,\adfsplit 64,57,117,1,61,113,118,13,35,30,125,16)_{\adfGB}$,

$(27,121,75,106,125,12,68,9,13,6,72,59,\adfsplit 74,126,26,8,84,112,25,10,117,88,5,20)_{\adfGB}$,

$(109,66,21,10,69,120,19,56,94,26,97,82,\adfsplit 53,67,102,50,32,6,51,112,35,86,5,61)_{\adfGB}$,

$(38,57,39,42,41,88,85,44,69,3,20,68,\adfsplit 33,81,35,52,10,58,27,49,74,79,96,40)_{\adfGB}$,

$(124,14,87,64,108,99,88,19,75,101,115,9,\adfsplit 82,67,59,121,33,130,30,107,110,52,24,131)_{\adfGB}$,

\adfSgap
$(91,122,62,74,132,46,54,31,48,117,92,108,\adfsplit 77,93,86,61,30,125,8,32,20,107,81,51)_{\adfGB}$,

$(2,71,129,61,77,38,51,98,44,123,62,56,\adfsplit 92,26,72,14,85,39,5,132,101,131,68,29)_{\adfGB}$,

\adfLgap
$(\infty,103,120,11,101,56,111,46,67,93,134,70,\adfsplit 57,115,35,118,28,77,50,7,131,36,119,126)_{\adfGC}$,

$(2,76,77,78,31,47,93,80,50,104,64,95,\adfsplit 99,111,18,28,86,58,38,17,35,133,37,123)_{\adfGC}$,

$(4,84,25,104,47,76,68,51,109,46,42,70,\adfsplit 134,118,38,98,121,44,132,36,129,127,14,69)_{\adfGC}$,

$(117,77,73,20,118,57,88,5,42,90,105,11,\adfsplit 71,109,126,1,3,35,134,55,65,129,66,93)_{\adfGC}$,

$(113,71,11,28,79,116,26,74,40,32,45,104,\adfsplit 51,72,2,57,108,18,132,3,82,102,76,4)_{\adfGC}$,

\adfSgap
$(25,48,28,12,1,7,123,45,132,27,114,57,\adfsplit 97,84,19,3,40,109,91,112,52,94,43,0)_{\adfGC}$,

$(0,15,54,88,127,69,94,67,84,70,121,112,\adfsplit 34,45,3,10,27,130,73,49,42,25,51,1)_{\adfGC}$,

\adfLgap
$(\infty,116,40,57,7,90,43,53,59,73,42,130,\adfsplit 102,3,107,16,105,112,89,26,76,117,92,49)_{\adfGD}$,

$(101,29,18,120,91,97,107,19,12,115,1,95,\adfsplit 69,126,2,98,34,76,8,103,131,92,61,47)_{\adfGD}$,

$(64,50,113,107,78,84,83,112,119,73,102,21,\adfsplit 63,103,111,18,86,67,120,71,22,77,118,43)_{\adfGD}$,

$(127,83,122,14,47,39,58,46,71,96,86,121,\adfsplit 133,75,54,93,51,67,18,45,38,88,99,114)_{\adfGD}$,

$(34,59,93,22,27,32,49,112,12,16,25,28,\adfsplit 60,67,36,87,21,74,131,54,129,55,132,37)_{\adfGD}$,

\adfSgap
$(115,22,20,50,62,64,0,110,8,29,134,106,\adfsplit 96,56,98,66,120,71,23,116,44,30,92,87)_{\adfGD}$,

$(1,15,101,41,11,35,33,31,86,78,0,72,\adfsplit 45,131,17,113,68,5,83,39,75,110,126,59)_{\adfGD}$,

\adfLgap
$(\infty,119,43,93,115,17,83,26,55,88,128,56,\adfsplit 132,63,120,74,33,129,82,98,113,36,20,84)_{\adfGE}$,

$(27,117,61,49,83,50,128,45,93,82,132,40,\adfsplit 109,76,121,71,10,56,8,77,72,64,58,115)_{\adfGE}$,

$(16,34,74,96,129,15,131,55,67,8,2,95,\adfsplit 124,46,70,105,114,51,99,120,23,108,72,45)_{\adfGE}$,

$(122,100,28,66,85,126,81,60,26,105,80,108,\adfsplit 73,17,35,83,52,91,72,45,43,0,102,20)_{\adfGE}$,

$(134,50,75,52,45,120,33,108,129,104,21,0,\adfsplit 115,118,61,5,11,37,23,31,96,76,49,131)_{\adfGE}$,

\adfSgap
$(45,53,44,62,125,5,88,67,110,2,123,68,\adfsplit 131,4,46,122,97,20,71,118,109,70,86,91)_{\adfGE}$,

$(2,65,31,112,134,57,70,52,68,14,74,56,\adfsplit 94,124,88,44,100,131,121,101,46,71,82,22)_{\adfGE}$,

\adfLgap
$(\infty,70,51,113,45,67,106,81,101,82,97,85,\adfsplit 10,123,64,31,9,48,2,77,40,88,7,86)_{\adfGF}$,

$(38,21,11,66,103,56,128,7,44,25,47,79,\adfsplit 115,96,33,107,52,92,41,87,85,95,68,121)_{\adfGF}$,

$(28,67,70,38,129,102,27,1,63,133,81,12,\adfsplit 80,18,49,132,43,121,8,50,100,58,93,23)_{\adfGF}$,

$(12,129,70,0,59,30,20,3,104,119,10,105,\adfsplit 44,101,52,48,16,132,117,13,40,35,2,123)_{\adfGF}$,

$(31,121,43,88,72,9,23,103,67,75,37,19,\adfsplit 81,99,38,20,44,0,5,3,132,122,84,98)_{\adfGF}$,

\adfSgap
$(87,63,83,36,116,52,90,81,56,86,18,8,\adfsplit 95,14,35,38,62,2,89,68,82,69,48,93)_{\adfGF}$,

$(2,89,83,77,87,41,113,93,121,8,51,116,\adfsplit 95,11,132,66,101,126,74,102,15,75,35,68)_{\adfGF}$,

\adfLgap
$(\infty,96,44,49,13,35,120,45,27,122,4,21,\adfsplit 98,3,127,88,89,28,129,85,11,22,23,67)_{\adfGG}$,

$(126,120,92,15,53,54,8,99,111,55,11,75,\adfsplit 66,40,95,52,110,129,56,10,67,116,21,134)_{\adfGG}$,

$(15,105,78,132,41,14,114,7,129,27,128,11,\adfsplit 16,26,3,64,72,37,130,34,84,80,110,61)_{\adfGG}$,

$(22,123,3,21,31,43,64,7,112,53,75,27,\adfsplit 124,35,28,78,90,23,93,24,107,41,20,118)_{\adfGG}$,

$(113,133,86,18,15,24,121,45,28,13,112,132,\adfsplit 33,120,92,54,76,32,91,5,35,64,128,104)_{\adfGG}$,

\adfSgap
$(117,79,32,80,37,101,41,17,97,58,22,134,\adfsplit 8,50,52,118,47,124,68,28,76,14,31,128)_{\adfGG}$,

$(6,104,56,103,34,119,121,116,113,61,21,38,\adfsplit 118,37,32,47,71,55,70,127,62,17,67,65)_{\adfGG}$,

\adfLgap
$(\infty,103,69,134,98,133,46,114,16,104,38,83,\adfsplit 120,9,84,122,81,71,6,63,25,72,80,108)_{\adfGH}$,

$(89,92,13,128,101,102,30,127,100,112,22,80,\adfsplit 83,52,6,54,133,76,94,46,97,123,48,35)_{\adfGH}$,

$(66,126,68,86,112,94,8,53,115,117,109,28,\adfsplit 60,121,39,33,120,79,37,50,13,2,98,48)_{\adfGH}$,

$(74,94,41,9,106,40,131,102,124,104,1,2,\adfsplit 134,10,46,16,35,65,56,111,120,82,92,13)_{\adfGH}$,

$(109,68,7,44,76,26,119,129,34,8,60,106,\adfsplit 64,85,39,9,70,33,43,93,21,6,14,30)_{\adfGH}$,

\adfSgap
$(109,77,108,102,56,69,35,115,113,123,72,76,\adfsplit 83,114,119,23,111,75,44,0,128,90,96,134)_{\adfGH}$,

$(35,87,102,21,92,14,3,6,24,65,71,126,\adfsplit 42,17,68,117,50,131,27,74,9,30,39,44)_{\adfGH}$,

\adfLgap
$(\infty,11,133,99,111,59,74,16,125,95,66,5,\adfsplit 94,76,51,79,89,30,82,93,13,128,12,113)_{\adfGI}$,

$(57,134,10,133,132,126,118,4,31,131,100,28,\adfsplit 45,87,67,122,69,33,21,75,42,24,97,8)_{\adfGI}$,

$(107,16,79,32,23,65,76,29,20,131,6,97,\adfsplit 47,52,115,88,93,77,41,54,72,55,87,112)_{\adfGI}$,

$(10,78,55,98,58,13,97,17,93,35,22,63,\adfsplit 81,74,24,41,32,23,82,45,4,111,88,34)_{\adfGI}$,

$(101,116,36,7,82,55,61,115,19,51,109,40,\adfsplit 4,66,8,83,110,23,92,35,59,48,111,125)_{\adfGI}$,

\adfSgap
$(5,110,47,129,80,118,69,6,30,120,90,114,\adfsplit 45,105,29,81,96,75,7,111,54,134,92,3)_{\adfGI}$,

$(76,48,68,72,123,108,96,98,122,36,119,3,\adfsplit 114,45,105,69,12,44,9,23,99,50,66,116)_{\adfGI}$,

\adfLgap
$(\infty,72,73,113,28,115,61,130,3,117,129,66,\adfsplit 26,132,97,111,1,32,5,80,98,105,114,107)_{\adfGJ}$,

$(89,24,112,4,61,131,115,16,13,20,118,34,\adfsplit 26,31,80,3,54,14,33,32,85,73,63,71)_{\adfGJ}$,

$(84,60,19,74,87,126,25,62,48,12,50,5,\adfsplit 34,47,0,130,105,104,133,55,30,72,75,128)_{\adfGJ}$,

$(87,40,128,102,90,59,23,84,13,3,54,131,\adfsplit 112,63,71,126,25,110,106,117,2,1,76,34)_{\adfGJ}$,

$(40,124,22,133,41,108,110,44,98,106,49,100,\adfsplit 62,7,38,20,75,11,64,97,84,43,54,53)_{\adfGJ}$,

\adfSgap
$(52,118,131,121,92,49,36,38,95,68,54,128,\adfsplit 91,123,96,12,101,35,125,45,62,93,88,102)_{\adfGJ}$,

$(41,129,60,126,78,107,18,92,125,134,51,105,\adfsplit 99,69,122,44,111,48,65,77,86,110,108,84)_{\adfGJ}$,

\adfLgap
$(\infty,112,122,54,86,109,60,81,51,9,116,97,\adfsplit 13,90,26,19,36,103,78,98,47,37,87,80)_{\adfGK}$,

$(24,38,68,119,33,18,132,73,49,51,75,111,\adfsplit 127,43,94,83,76,32,112,92,34,74,19,71)_{\adfGK}$,

$(130,51,49,103,29,39,0,12,128,92,26,57,\adfsplit 63,78,133,40,20,5,64,11,127,13,90,72)_{\adfGK}$,

$(27,40,133,6,128,64,49,34,53,30,60,44,\adfsplit 125,57,96,0,70,108,124,69,119,37,62,15)_{\adfGK}$,

$(75,47,26,81,25,84,27,3,114,70,122,40,\adfsplit 112,97,64,80,18,77,55,17,131,34,99,119)_{\adfGK}$,

\adfSgap
$(24,82,39,9,29,48,97,54,134,67,38,79,\adfsplit 58,88,85,34,14,125,118,47,116,104,101,50)_{\adfGK}$,

$(96,130,127,86,134,41,31,8,56,110,32,118,\adfsplit 2,47,35,64,14,71,89,16,68,74,122,98)_{\adfGK}$

\noindent and

$(\infty,12,46,92,123,2,67,77,91,116,33,76,\adfsplit 3,58,113,70,125,7,100,114,63,20,131,22)_{\adfGL}$,

$(108,111,17,69,5,82,80,62,46,124,53,13,\adfsplit 58,96,76,133,114,19,8,38,105,3,121,132)_{\adfGL}$,

$(70,134,24,34,129,108,96,85,104,42,120,32,\adfsplit 83,74,12,109,23,76,88,113,55,133,90,118)_{\adfGL}$,

$(131,47,54,89,113,108,81,13,105,96,22,23,\adfsplit 27,100,17,130,116,36,29,31,7,74,119,78)_{\adfGL}$,

$(133,83,130,8,118,98,55,111,56,47,49,7,\adfsplit 72,19,65,31,124,102,117,1,73,123,3,40)_{\adfGL}$,

\adfSgap
$(57,93,56,128,5,23,83,134,126,36,21,80,\adfsplit 104,42,39,0,35,105,24,9,107,78,129,38)_{\adfGL}$,

$(131,39,75,60,126,90,47,123,38,113,45,42,\adfsplit 84,108,95,5,0,98,132,89,2,87,116,66)_{\adfGL}$

\noindent under the action of the mapping $\infty \mapsto \infty$, $x \mapsto x + 3$ (mod 135)
for the first five graphs in each design, and $x \mapsto x + 9$ (mod 135) for the last two. \eproof

\adfVfy{136, \{\{135,15,9\},\{1,1,1\}\}, 135, -1, \{5,\{\{135,3,3\},\{1,1,1\}\}\}} 


{\lemma \label{lem:Snark24 designs}
There exist designs of order $64$, $73$ and $145$
for each of the thirty-eight 24-vertex non-trivial snarks.}\\

\noindent\textbf{Proof.}
The decompositions are given in the Appendix. \eproof


{\lemma \label{lem:Snark24 multipartite}
There exist decompositions of the complete multipartite graphs
$K_{12,12,12}$, $K_{24,24,15}$, $K_{72,72,63}$, $K_{24,24,24,24}$ and $K_{24,24,24,21}$
for each of the thirty-eight 24-vertex non-trivial snarks.}\\

\noindent\textbf{Proof.}
The decompositions are given in the Appendix. \eproof


\vskip 2mm
Theorem~\ref{thm:snark24} follows from
Lemmas~\ref{lem:Snark24 136}, \ref{lem:Snark24 designs} and \ref{lem:Snark24 multipartite},
and Proposition~\ref{prop:d=3, v=24}.





\newpage
\appendix
\section{Proof of Lemma~\ref{lem:Snark24 designs}}
\label{app:Snark24 designs}


\noindent\textbf{Proof of Lemma~\ref{lem:Snark24 designs}.}
Recall that the lemma asserts the existence of designs of order $64$, $73$ and $145$
for each of the thirty-eight 24-vertex non-trivial snarks.

Let the vertex set of $K_{64}$ be $Z_{63} \cup \{\infty\}$. The decompositions consist of

$(\infty,53,7,45,29,21,19,59,46,47,33,24,\adfsplit 5,28,11,6,8,14,58,17,2,9,43,50)_{\adfGa}$,

$(58,17,51,53,13,60,7,43,9,52,32,20,\adfsplit 30,24,45,62,28,0,4,15,2,31,59,36)_{\adfGa}$,

\adfSgap
$(11,56,62,49,37,5,9,17,6,18,15,61,\adfsplit 54,16,41,23,0,48,22,30,33,40,13,10)_{\adfGa}$,

$(25,4,61,44,42,28,7,37,51,19,49,33,\adfsplit 10,21,16,30,43,47,40,27,3,1,54,57)_{\adfGa}$,

\adfLgap
$(\infty,20,51,16,38,50,22,35,12,46,15,44,\adfsplit 1,24,57,37,40,3,48,28,58,55,4,33)_{\adfGb}$,

$(26,27,29,15,55,19,38,33,47,21,45,24,\adfsplit 14,10,6,62,48,30,59,20,51,31,32,8)_{\adfGb}$,

\adfSgap
$(22,50,15,24,56,44,51,31,37,17,19,11,\adfsplit 8,47,49,25,43,23,46,4,41,53,30,61)_{\adfGb}$,

$(43,62,11,36,23,56,18,16,0,49,28,34,\adfsplit 19,61,47,21,59,58,4,48,10,46,14,53)_{\adfGb}$,

\adfLgap
$(\infty,8,48,16,49,9,18,10,55,1,3,28,\adfsplit 50,56,26,7,2,31,13,59,40,54,33,43)_{\adfGc}$,

$(2,51,52,58,48,26,39,24,4,13,16,10,\adfsplit 61,40,18,33,29,56,11,44,3,45,23,6)_{\adfGc}$,

\adfSgap
$(25,8,17,41,27,62,9,45,0,50,3,5,\adfsplit 29,23,13,44,36,55,47,33,59,15,26,38)_{\adfGc}$,

$(56,17,27,15,6,44,49,7,23,42,60,2,\adfsplit 48,50,38,32,29,21,57,39,46,3,62,11)_{\adfGc}$,

\adfLgap
$(\infty,43,17,24,45,42,29,30,0,27,52,34,\adfsplit 8,48,35,1,41,54,61,5,7,31,56,55)_{\adfGd}$,

$(3,57,38,26,40,62,19,29,51,10,41,17,\adfsplit 23,16,56,14,24,34,52,0,21,11,61,49)_{\adfGd}$,

\adfSgap
$(41,42,4,19,38,55,12,47,1,27,40,25,\adfsplit 18,60,6,36,62,29,17,30,52,32,7,51)_{\adfGd}$,

$(13,35,31,41,54,55,24,12,61,9,36,4,\adfsplit 6,57,52,30,59,17,33,21,27,37,48,2)_{\adfGd}$,

\adfLgap
$(\infty,46,36,41,56,11,9,42,2,38,33,22,\adfsplit 43,19,4,58,27,21,18,10,48,8,23,59)_{\adfGe}$,

$(4,44,36,56,15,57,12,31,60,35,9,29,\adfsplit 11,16,42,19,39,8,30,7,53,10,52,0)_{\adfGe}$,

\adfSgap
$(18,22,1,61,35,42,5,46,62,25,51,50,\adfsplit 54,0,60,43,17,28,23,2,49,58,59,53)_{\adfGe}$,

$(30,13,21,34,19,38,49,1,52,4,39,31,\adfsplit 8,59,61,2,53,46,29,32,16,20,47,10)_{\adfGe}$,

\adfLgap
$(\infty,24,52,8,38,33,0,16,39,60,27,41,\adfsplit 10,5,9,22,3,59,36,18,6,53,51,37)_{\adfGf}$,

$(24,44,14,58,52,31,13,54,60,57,22,16,\adfsplit 28,7,62,5,3,50,37,21,19,45,8,34)_{\adfGf}$,

\adfSgap
$(47,51,2,54,5,37,38,35,8,19,7,23,\adfsplit 4,32,56,27,24,59,50,13,52,17,44,20)_{\adfGf}$,

$(48,34,2,28,38,4,26,14,35,58,16,30,\adfsplit 32,29,44,23,62,31,7,37,17,47,19,5)_{\adfGf}$,

\adfLgap
$(\infty,15,38,13,43,25,17,42,57,35,53,46,\adfsplit 29,32,20,21,26,47,49,31,40,2,7,48)_{\adfGg}$,

$(28,25,21,6,30,45,39,58,22,33,15,11,\adfsplit 8,24,26,54,20,0,14,3,53,43,29,12)_{\adfGg}$,

\adfSgap
$(42,13,46,36,49,9,28,32,7,21,56,35,\adfsplit 18,25,27,51,23,33,16,48,4,2,19,10)_{\adfGg}$,

$(58,4,2,50,8,10,42,61,3,52,46,48,\adfsplit 29,26,25,14,34,43,40,31,17,19,49,28)_{\adfGg}$,

\adfLgap
$(\infty,60,41,7,26,30,3,58,15,54,20,9,\adfsplit 22,11,38,34,35,56,33,14,31,13,52,6)_{\adfGh}$,

$(46,31,55,2,19,59,25,53,58,21,4,42,\adfsplit 38,54,1,41,27,60,14,16,26,32,9,5)_{\adfGh}$,

\adfSgap
$(59,49,46,19,26,0,33,56,24,42,22,1,\adfsplit 23,21,51,35,15,7,3,12,52,54,2,25)_{\adfGh}$,

$(37,27,53,0,36,4,23,54,39,1,57,30,\adfsplit 22,13,59,33,7,21,38,8,18,60,20,61)_{\adfGh}$,

\adfLgap
$(\infty,49,24,17,40,46,1,9,6,7,27,34,\adfsplit 25,5,53,33,30,2,3,56,8,45,39,55)_{\adfGi}$,

$(54,41,15,8,11,34,46,22,18,28,32,17,\adfsplit 16,36,49,13,3,12,51,61,31,35,33,43)_{\adfGi}$,

\adfSgap
$(32,59,58,53,18,43,45,11,14,40,28,2,\adfsplit 9,1,55,44,23,56,47,38,51,17,29,3)_{\adfGi}$,

$(27,62,8,41,25,42,40,26,20,60,15,35,\adfsplit 11,44,32,29,17,61,50,48,3,14,23,12)_{\adfGi}$,

\adfLgap
$(\infty,53,36,55,14,24,37,41,59,5,9,29,\adfsplit 8,17,42,4,60,32,20,49,6,19,56,43)_{\adfGj}$,

$(44,60,28,53,29,18,6,7,54,58,42,62,\adfsplit 50,56,13,2,40,21,35,16,25,48,52,39)_{\adfGj}$,

\adfSgap
$(29,10,23,9,0,27,3,17,25,37,36,12,\adfsplit 11,16,28,61,33,6,21,39,40,19,46,49)_{\adfGj}$,

$(16,6,43,3,42,59,62,40,61,30,24,27,\adfsplit 12,15,33,22,49,7,18,58,31,10,54,55)_{\adfGj}$,

\adfLgap
$(\infty,44,43,15,55,41,50,18,45,26,23,37,\adfsplit 20,3,7,35,48,1,12,61,8,42,59,38)_{\adfGk}$,

$(26,47,38,39,15,22,55,32,43,21,14,23,\adfsplit 52,48,57,28,0,37,56,50,41,51,29,13)_{\adfGk}$,

\adfSgap
$(18,44,6,30,9,0,57,61,52,56,1,58,\adfsplit 19,55,12,21,49,27,37,40,31,16,15,34)_{\adfGk}$,

$(42,14,28,27,33,51,48,60,16,61,22,57,\adfsplit 58,4,43,10,19,25,9,13,55,24,18,21)_{\adfGk}$,

\adfLgap
$(\infty,49,23,15,36,29,37,3,50,51,26,54,\adfsplit 5,60,56,40,53,22,47,14,16,58,35,55)_{\adfGl}$,

$(43,28,35,26,21,22,51,48,31,16,50,11,\adfsplit 56,13,52,5,57,10,32,60,33,29,54,14)_{\adfGl}$,

\adfSgap
$(50,2,35,21,36,17,61,10,42,39,51,7,\adfsplit 47,24,4,15,3,58,22,40,30,28,25,18)_{\adfGl}$,

$(34,6,59,39,1,27,52,22,55,16,3,36,\adfsplit 51,33,12,45,37,54,13,0,10,24,9,48)_{\adfGl}$,

\adfLgap
$(\infty,44,39,34,26,61,59,4,47,19,24,23,\adfsplit 60,41,20,40,30,1,43,12,50,56,18,3)_{\adfGm}$,

$(4,31,56,55,51,47,22,16,0,54,10,26,\adfsplit 11,52,12,58,40,42,46,15,43,28,38,5)_{\adfGm}$,

\adfSgap
$(48,51,2,53,24,46,43,60,41,31,30,57,\adfsplit 26,42,29,36,54,6,28,39,50,58,9,38)_{\adfGm}$,

$(6,41,17,47,12,43,39,48,35,36,27,57,\adfsplit 50,0,2,44,5,26,45,11,32,24,42,9)_{\adfGm}$,

\adfLgap
$(\infty,36,8,22,13,57,62,3,26,5,56,12,\adfsplit 49,44,35,40,42,7,21,6,30,28,50,31)_{\adfGn}$,

$(1,0,21,60,24,27,59,62,54,34,35,39,\adfsplit 45,23,37,61,42,20,48,4,6,2,50,55)_{\adfGn}$,

\adfSgap
$(3,14,16,19,27,22,43,35,38,11,17,40,\adfsplit 50,24,2,47,49,55,7,29,32,13,37,28)_{\adfGn}$,

$(44,37,43,22,4,46,35,50,52,9,5,61,\adfsplit 34,40,42,7,55,31,53,41,39,26,2,10)_{\adfGn}$,

\adfLgap
$(\infty,47,16,24,38,42,9,5,26,25,29,37,\adfsplit 52,18,20,6,22,0,1,30,55,28,41,40)_{\adfGo}$,

$(21,38,33,19,4,41,29,10,26,1,23,9,\adfsplit 36,22,51,8,20,58,44,61,14,56,24,18)_{\adfGo}$,

\adfSgap
$(50,7,2,15,62,38,32,17,29,18,57,46,\adfsplit 43,49,27,20,1,24,33,60,48,10,42,25)_{\adfGo}$,

$(13,57,51,48,30,12,36,32,33,39,58,26,\adfsplit 40,16,54,41,25,31,19,10,18,0,62,27)_{\adfGo}$,

\adfLgap
$(\infty,8,54,49,25,35,27,20,47,58,12,51,\adfsplit 42,7,9,57,55,43,44,16,56,21,11,24)_{\adfGp}$,

$(49,61,23,9,54,47,29,27,28,14,21,56,\adfsplit 51,38,53,36,42,39,58,3,22,59,48,2)_{\adfGp}$,

\adfSgap
$(29,48,58,11,37,10,34,31,53,49,55,57,\adfsplit 59,8,60,19,54,25,13,43,7,16,46,32)_{\adfGp}$,

$(52,44,17,22,10,19,29,59,16,37,47,13,\adfsplit 41,20,58,31,34,49,61,6,53,50,46,9)_{\adfGp}$,

\adfLgap
$(\infty,37,42,23,5,55,39,35,18,24,34,36,\adfsplit 45,41,40,27,3,2,16,57,49,29,13,25)_{\adfGq}$,

$(18,39,13,6,22,20,38,15,0,61,24,55,\adfsplit 51,37,49,48,3,23,33,11,8,28,44,50)_{\adfGq}$,

\adfSgap
$(13,16,23,27,55,26,53,12,38,4,44,15,\adfsplit 61,6,1,11,43,41,17,46,20,29,50,32)_{\adfGq}$,

$(48,62,25,13,56,5,8,10,61,28,23,22,\adfsplit 45,1,35,26,29,14,59,20,40,31,16,11)_{\adfGq}$,

\adfLgap
$(\infty,51,10,23,60,20,7,53,39,24,57,49,\adfsplit 2,37,44,31,0,40,52,45,9,25,16,17)_{\adfGr}$,

$(48,23,7,51,56,29,47,36,32,39,42,18,\adfsplit 9,54,0,12,37,24,10,14,59,52,38,25)_{\adfGr}$,

\adfSgap
$(36,29,38,53,32,19,35,31,16,37,10,15,\adfsplit 18,3,44,40,8,22,47,59,20,61,4,49)_{\adfGr}$,

$(32,30,47,5,34,59,16,43,1,23,13,42,\adfsplit 41,52,55,22,44,26,19,11,61,10,8,56)_{\adfGr}$,

\adfLgap
$(\infty,0,26,40,31,51,22,54,33,21,16,10,\adfsplit 11,35,7,52,57,30,49,20,46,36,2,43)_{\adfGs}$,

$(43,53,31,8,46,1,25,23,59,30,60,19,\adfsplit 12,40,22,7,20,14,11,44,2,39,26,41)_{\adfGs}$,

\adfSgap
$(1,49,16,62,8,47,38,31,41,15,57,0,\adfsplit 54,19,24,42,39,21,45,9,6,5,12,2)_{\adfGs}$,

$(29,8,24,50,18,12,3,21,54,20,26,57,\adfsplit 44,14,6,9,33,48,0,51,23,45,60,25)_{\adfGs}$,

\adfLgap
$(\infty,36,11,43,44,16,22,32,41,40,50,58,\adfsplit 55,0,56,2,10,4,26,53,35,54,9,25)_{\adfGt}$,

$(28,58,34,41,3,27,24,12,26,23,17,18,\adfsplit 59,37,21,7,11,55,36,46,57,19,42,9)_{\adfGt}$,

\adfSgap
$(36,1,47,2,37,48,11,9,46,60,54,42,\adfsplit 29,21,25,20,12,5,33,44,41,40,49,30)_{\adfGt}$,

$(22,49,59,5,15,33,39,35,54,10,27,36,\adfsplit 24,30,3,41,43,53,8,0,26,52,16,51)_{\adfGt}$,

\adfLgap
$(\infty,10,35,42,32,61,26,7,44,46,18,48,\adfsplit 13,53,62,14,34,50,3,21,58,39,55,4)_{\adfGu}$,

$(59,60,29,17,11,55,3,19,0,27,32,40,\adfsplit 20,37,46,4,39,5,44,50,43,57,62,18)_{\adfGu}$,

\adfSgap
$(20,22,56,49,51,11,60,53,18,57,14,61,\adfsplit 40,50,1,59,21,0,54,3,7,25,27,33)_{\adfGu}$,

$(45,30,57,33,61,58,27,28,1,55,51,60,\adfsplit 62,54,15,3,52,35,46,9,39,12,49,24)_{\adfGu}$,

\adfLgap
$(\infty,11,6,58,30,61,36,17,16,57,2,14,\adfsplit 51,47,34,1,50,22,32,25,45,21,10,40)_{\adfGv}$,

$(33,60,55,32,10,58,18,19,7,6,12,23,\adfsplit 5,37,21,22,16,31,20,45,46,54,50,52)_{\adfGv}$,

\adfSgap
$(55,39,14,50,11,47,36,60,8,59,24,33,\adfsplit 57,32,2,40,51,17,45,5,35,62,38,29)_{\adfGv}$,

$(39,8,36,56,30,17,5,50,14,15,41,42,\adfsplit 57,53,35,29,45,12,2,27,38,9,62,43)_{\adfGv}$,

\adfLgap
$(\infty,42,52,5,28,25,30,9,60,32,59,38,\adfsplit 14,2,18,19,58,7,37,47,55,17,13,36)_{\adfGw}$,

$(8,52,53,34,10,37,9,31,20,57,48,62,\adfsplit 25,22,7,40,61,41,6,58,29,33,3,54)_{\adfGw}$,

\adfSgap
$(19,2,54,9,18,14,30,45,41,8,38,11,\adfsplit 0,39,59,60,53,42,51,35,29,20,57,26)_{\adfGw}$,

$(62,14,16,36,24,23,51,6,42,0,40,2,\adfsplit 5,30,53,12,11,26,57,50,21,54,15,3)_{\adfGw}$,

\adfLgap
$(\infty,33,53,43,20,28,46,11,30,59,42,32,\adfsplit 52,45,48,10,49,12,3,44,36,37,1,26)_{\adfGx}$,

$(30,55,31,22,58,26,7,49,52,32,28,14,\adfsplit 29,45,5,37,46,13,17,3,60,0,2,57)_{\adfGx}$,

\adfSgap
$(50,12,3,57,45,30,55,2,34,20,40,15,\adfsplit 27,33,59,0,25,60,18,53,35,41,11,44)_{\adfGx}$,

$(17,11,47,18,20,30,14,24,36,44,62,45,\adfsplit 5,13,21,6,15,26,41,51,48,32,55,29)_{\adfGx}$,

\adfLgap
$(\infty,0,23,46,6,60,20,27,40,10,25,1,\adfsplit 59,7,41,15,8,13,50,56,48,35,30,61)_{\adfGy}$,

$(9,24,49,45,3,15,50,54,4,6,48,33,\adfsplit 17,22,8,2,52,11,25,31,60,26,36,13)_{\adfGy}$,

\adfSgap
$(42,30,54,49,37,55,61,16,51,12,22,57,\adfsplit 50,46,11,5,45,2,56,43,44,19,23,40)_{\adfGy}$,

$(62,34,50,28,38,23,17,32,4,11,15,44,\adfsplit 26,7,52,40,45,20,37,31,35,53,2,59)_{\adfGy}$,

\adfLgap
$(\infty,39,26,4,54,43,44,13,10,21,12,2,\adfsplit 47,61,28,14,35,53,3,41,38,49,9,46)_{\adfGz}$,

$(11,6,20,32,5,61,27,34,43,1,10,25,\adfsplit 21,16,0,56,30,45,19,40,28,41,47,31)_{\adfGz}$,

\adfSgap
$(16,46,17,49,47,44,60,56,50,24,8,41,\adfsplit 42,20,57,0,14,2,3,9,51,30,27,18)_{\adfGz}$,

$(43,18,41,9,30,54,17,21,48,6,42,36,\adfsplit 57,46,51,13,50,8,15,12,45,11,3,24)_{\adfGz}$,

\adfLgap
$(\infty,41,22,48,9,38,30,18,11,46,59,10,\adfsplit 17,3,47,1,56,51,45,15,37,13,44,19)_{\adfGA}$,

$(23,15,59,43,45,49,42,9,40,17,2,34,\adfsplit 32,19,21,52,53,20,62,12,26,61,6,22)_{\adfGA}$,

\adfSgap
$(9,22,33,0,54,42,34,51,36,1,7,15,\adfsplit 12,19,43,24,6,4,25,8,11,2,56,50)_{\adfGA}$,

$(31,1,61,56,32,28,50,59,17,9,53,49,\adfsplit 35,12,55,22,21,48,57,58,27,3,16,46)_{\adfGA}$,

\adfLgap
$(\infty,40,24,20,55,15,28,1,48,31,38,36,\adfsplit 14,19,61,37,51,35,30,58,41,11,49,33)_{\adfGB}$,

$(30,24,57,44,62,17,40,37,58,34,1,0,\adfsplit 56,15,29,12,10,16,59,2,3,25,54,23)_{\adfGB}$,

\adfSgap
$(28,27,11,59,9,13,33,0,18,38,21,1,\adfsplit 55,5,53,3,35,52,6,36,42,40,23,51)_{\adfGB}$,

$(48,9,24,59,5,53,56,16,14,11,32,57,\adfsplit 62,42,2,61,38,58,21,12,35,29,17,26)_{\adfGB}$,

\adfLgap
$(\infty,20,54,7,2,14,38,23,4,9,51,36,\adfsplit 10,35,25,32,56,49,30,61,21,31,55,47)_{\adfGC}$,

$(25,35,56,9,24,28,57,2,59,0,16,39,\adfsplit 40,43,38,48,29,53,58,34,5,62,36,54)_{\adfGC}$,

\adfSgap
$(40,30,59,5,6,51,42,62,2,25,53,22,\adfsplit 34,45,37,15,60,52,7,41,3,12,36,39)_{\adfGC}$,

$(61,54,34,6,18,29,1,0,55,9,12,27,\adfsplit 49,58,10,19,3,48,51,22,40,46,62,2)_{\adfGC}$,

\adfLgap
$(\infty,3,20,43,38,5,46,59,29,39,18,36,\adfsplit 48,41,31,30,16,44,9,15,13,45,61,10)_{\adfGD}$,

$(45,30,3,38,58,49,40,41,17,7,34,50,\adfsplit 29,37,12,47,54,9,2,51,48,46,36,32)_{\adfGD}$,

\adfSgap
$(57,0,14,8,45,6,62,60,2,19,24,12,\adfsplit 37,17,50,38,30,52,53,47,49,43,59,40)_{\adfGD}$,

$(22,1,10,50,55,26,34,49,20,25,29,19,\adfsplit 13,11,14,44,7,46,0,31,16,59,40,52)_{\adfGD}$,

\adfLgap
$(\infty,43,38,18,44,56,52,51,36,29,30,17,\adfsplit 53,55,22,9,7,21,62,11,27,61,26,40)_{\adfGE}$,

$(62,17,10,56,22,59,57,58,15,31,43,52,\adfsplit 61,46,24,18,50,13,30,19,25,6,33,23)_{\adfGE}$,

\adfSgap
$(25,62,42,5,31,16,14,52,54,39,27,47,\adfsplit 12,17,37,11,59,50,0,3,2,28,35,51)_{\adfGE}$,

$(42,6,26,12,54,36,60,59,1,48,0,3,\adfsplit 29,62,47,24,13,27,15,2,4,41,21,53)_{\adfGE}$,

\adfLgap
$(\infty,58,41,21,47,54,43,57,35,12,34,19,\adfsplit 49,38,25,33,2,45,50,4,5,62,6,60)_{\adfGF}$,

$(25,4,6,36,53,31,30,11,13,5,59,34,\adfsplit 52,20,42,18,14,33,45,0,24,44,51,7)_{\adfGF}$,

\adfSgap
$(48,43,30,0,4,35,23,15,58,29,62,55,\adfsplit 33,10,1,50,47,61,46,49,21,52,31,34)_{\adfGF}$,

$(19,31,2,11,25,23,8,41,55,47,38,50,\adfsplit 17,35,16,4,45,7,32,40,37,60,27,53)_{\adfGF}$,

\adfLgap
$(\infty,17,15,46,27,1,52,53,24,39,0,60,\adfsplit 59,26,23,12,9,28,57,21,4,32,58,18)_{\adfGG}$,

$(4,26,51,33,9,18,2,7,3,11,6,0,\adfsplit 56,32,31,61,38,58,28,48,37,5,55,35)_{\adfGG}$,

\adfSgap
$(50,5,38,61,33,42,12,56,13,55,0,23,\adfsplit 21,32,34,44,43,31,2,58,29,10,14,26)_{\adfGG}$,

$(26,0,8,58,44,35,10,32,56,22,61,43,\adfsplit 17,13,19,16,31,28,38,25,37,20,1,59)_{\adfGG}$,

\adfLgap
$(\infty,43,0,5,20,14,50,46,57,11,52,7,\adfsplit 21,1,37,15,25,3,48,9,59,18,26,27)_{\adfGH}$,

$(51,22,58,42,61,32,49,9,21,16,30,31,\adfsplit 3,6,35,44,15,23,52,50,7,38,60,10)_{\adfGH}$,

\adfSgap
$(41,57,48,58,22,29,37,53,38,52,11,44,\adfsplit 7,18,15,36,20,1,17,26,23,43,4,8)_{\adfGH}$,

$(45,38,34,5,30,19,53,47,46,60,13,2,\adfsplit 32,25,59,44,26,28,29,41,40,24,23,14)_{\adfGH}$,

\adfLgap
$(\infty,50,43,48,59,61,21,5,22,6,49,45,\adfsplit 24,15,30,9,4,17,41,28,18,2,25,42)_{\adfGI}$,

$(0,25,1,50,6,30,34,60,39,45,8,57,\adfsplit 15,9,47,29,44,40,17,19,23,61,55,5)_{\adfGI}$,

\adfSgap
$(60,23,53,14,3,19,37,20,29,47,46,62,\adfsplit 35,6,9,4,22,40,52,26,25,61,41,59)_{\adfGI}$,

$(22,58,25,32,53,55,52,38,39,18,16,62,\adfsplit 34,31,56,11,37,2,59,40,10,46,5,44)_{\adfGI}$,

\adfLgap
$(\infty,53,51,25,37,46,58,62,21,45,50,43,\adfsplit 19,26,24,1,32,60,11,16,38,27,23,61)_{\adfGJ}$,

$(17,29,59,14,4,54,58,51,43,31,44,1,\adfsplit 19,52,8,3,15,34,56,27,55,33,53,61)_{\adfGJ}$,

\adfSgap
$(41,19,9,42,7,31,21,8,30,12,54,43,\adfsplit 38,57,36,39,26,45,2,27,60,33,18,3)_{\adfGJ}$,

$(38,5,62,51,54,33,40,8,24,42,57,48,\adfsplit 6,41,0,59,50,46,2,17,9,11,60,30)_{\adfGJ}$,

\adfLgap
$(\infty,34,50,51,37,42,18,40,25,59,4,28,\adfsplit 5,23,56,38,0,54,7,33,12,39,3,16)_{\adfGK}$,

$(24,62,60,12,29,59,2,40,53,30,58,47,\adfsplit 28,9,42,50,45,6,25,18,10,20,52,14)_{\adfGK}$,

\adfSgap
$(60,18,59,13,53,39,36,50,38,37,56,11,\adfsplit 1,43,35,2,21,49,12,55,34,5,22,33)_{\adfGK}$,

$(13,8,41,49,40,53,59,24,25,55,62,34,\adfsplit 1,52,47,19,16,0,4,14,26,46,29,61)_{\adfGK}$

\noindent and

$(\infty,55,5,21,24,54,52,9,38,59,53,16,\adfsplit 30,62,58,7,41,46,17,40,36,50,3,47)_{\adfGL}$,

$(16,47,34,18,39,5,8,38,20,31,53,7,\adfsplit 43,15,2,40,10,13,24,41,33,30,22,59)_{\adfGL}$,

\adfSgap
$(47,0,55,56,26,36,17,49,6,31,27,54,\adfsplit 30,1,34,51,58,39,57,10,42,16,46,52)_{\adfGL}$,

$(48,12,30,24,54,36,33,51,42,50,29,40,\adfsplit 6,35,7,3,59,9,32,0,46,43,58,39)_{\adfGL}$

\noindent under the action of the mapping $\infty \mapsto \infty$, $x \mapsto x + 3$ (mod 63)
for the first two graphs in each design, and $x \mapsto x + 9$ (mod 63) for the last two.

\adfVfy{64, \{\{63,7,9\},\{1,1,1\}\}, 63, -1, \{2,\{\{63,3,3\},\{1,1,1\}\}\}} 

Let the vertex set of $K_{73}$ be $Z_{73}$. The decompositions consist of

$(0,1,2,3,5,6,11,13,15,16,24,30,\adfsplit 28,44,4,47,49,10,35,70,8,21,61,45)_{\adfGa}$,

\adfLgap
$(0,1,2,3,5,6,11,13,15,16,24,25,\adfsplit 28,44,42,4,47,12,70,39,20,8,45,63)_{\adfGb}$,

\adfLgap
$(0,1,2,3,5,6,11,13,15,16,24,25,\adfsplit 28,31,42,4,9,50,54,7,10,48,43,23)_{\adfGc}$,

\adfLgap
$(0,1,2,3,5,6,11,13,15,16,24,25,\adfsplit 28,31,42,43,4,8,55,53,7,72,45,22)_{\adfGd}$,

\adfLgap
$(0,1,2,3,5,6,11,13,15,16,14,30,\adfsplit 28,45,34,54,4,7,60,18,31,52,29,59)_{\adfGe}$,

\adfLgap
$(0,1,2,3,5,6,11,13,15,16,14,30,\adfsplit 28,44,34,47,4,7,64,8,32,50,29,55)_{\adfGf}$,

\adfLgap
$(0,1,2,3,5,6,11,13,15,16,14,30,\adfsplit 28,31,34,47,50,7,8,55,63,20,18,46)_{\adfGg}$,

\adfLgap
$(0,1,2,3,5,6,11,13,15,16,14,30,\adfsplit 28,31,34,48,8,52,7,4,54,64,17,36)_{\adfGh}$,

\adfLgap
$(0,1,2,3,5,6,11,13,15,16,14,30,\adfsplit 28,31,33,48,54,52,8,4,53,71,10,45)_{\adfGi}$,

\adfLgap
$(0,1,2,3,5,6,11,13,15,16,14,25,\adfsplit 28,31,42,34,53,61,51,4,71,52,9,39)_{\adfGj}$,

\adfLgap
$(0,1,2,3,5,6,11,13,15,16,14,25,\adfsplit 28,31,42,39,51,71,50,70,7,49,10,32)_{\adfGk}$,

\adfLgap
$(0,1,2,3,5,6,11,13,15,16,14,25,\adfsplit 28,31,33,45,63,72,56,7,8,58,20,39)_{\adfGl}$,

\adfLgap
$(0,1,2,3,5,6,11,13,15,16,14,25,\adfsplit 28,31,42,37,60,67,52,72,70,50,9,30)_{\adfGm}$,

\adfLgap
$(0,1,2,3,5,6,11,9,13,22,21,19,\adfsplit 27,39,37,40,59,4,61,15,66,23,33,54)_{\adfGn}$,

\adfLgap
$(0,1,2,3,5,6,11,9,13,22,21,19,\adfsplit 23,37,39,40,55,50,69,12,62,15,25,46)_{\adfGo}$,

\adfLgap
$(0,1,2,3,5,6,11,9,13,22,21,19,\adfsplit 23,39,40,4,41,56,8,34,10,70,36,58)_{\adfGp}$,

\adfLgap
$(0,1,2,3,5,6,11,9,13,22,21,19,\adfsplit 23,29,39,57,46,49,72,4,51,14,7,37)_{\adfGq}$,

\adfLgap
$(0,1,2,3,5,6,11,9,13,22,21,19,\adfsplit 23,29,39,43,66,51,71,8,53,7,20,42)_{\adfGr}$,

\adfLgap
$(0,1,2,3,5,6,11,9,13,22,21,19,\adfsplit 23,29,4,49,56,69,62,71,25,39,15,42)_{\adfGs}$,

\adfLgap
$(0,1,2,3,5,6,11,9,13,22,18,19,\adfsplit 27,33,38,55,46,62,66,10,12,14,40,41)_{\adfGt}$,

\adfLgap
$(0,1,2,3,5,6,11,9,13,22,18,19,\adfsplit 27,33,38,40,62,49,72,64,4,8,37,31)_{\adfGu}$,

\adfLgap
$(0,1,2,3,5,6,11,9,13,22,18,19,\adfsplit 27,33,39,34,40,53,10,8,7,56,12,46)_{\adfGv}$,

\adfLgap
$(0,1,2,3,5,6,11,9,13,22,18,19,\adfsplit 27,33,39,40,48,4,49,14,10,61,15,46)_{\adfGw}$,

\adfLgap
$(0,1,2,3,5,6,11,9,13,22,18,19,\adfsplit 27,33,38,44,50,39,58,15,72,17,21,60)_{\adfGx}$,

\adfLgap
$(0,1,2,3,5,6,11,9,13,22,18,19,\adfsplit 27,24,38,41,43,65,70,69,63,12,33,30)_{\adfGy}$,

\adfLgap
$(0,1,2,3,5,6,11,9,13,22,18,19,\adfsplit 27,24,38,45,49,44,10,52,64,4,26,15)_{\adfGz}$,

\adfLgap
$(0,1,2,3,5,6,11,9,13,22,18,19,\adfsplit 27,24,38,35,45,42,69,63,68,71,17,36)_{\adfGA}$,

\adfLgap
$(0,1,2,3,5,6,11,9,13,22,18,19,\adfsplit 27,24,38,35,53,46,51,72,61,68,23,20)_{\adfGB}$,

\adfLgap
$(0,1,2,3,5,6,11,9,13,22,18,19,\adfsplit 32,24,29,47,60,43,4,42,7,46,20,72)_{\adfGC}$,

\adfLgap
$(0,1,2,3,5,6,11,9,13,22,18,19,\adfsplit 32,24,40,39,37,56,65,69,8,57,27,30)_{\adfGD}$,

\adfLgap
$(0,1,2,3,5,6,11,9,13,22,18,19,\adfsplit 32,24,40,39,59,54,68,72,64,20,23,42)_{\adfGE}$,

\adfLgap
$(0,1,2,3,5,6,11,9,13,22,18,19,\adfsplit 32,24,40,39,42,37,4,62,60,71,20,30)_{\adfGF}$,

\adfLgap
$(0,1,2,3,5,6,11,9,13,22,18,19,\adfsplit 32,24,40,39,36,47,69,55,17,15,60,21)_{\adfGG}$,

\adfLgap
$(0,1,2,3,5,6,11,9,13,22,18,19,\adfsplit 32,24,29,4,43,56,45,7,52,72,51,10)_{\adfGH}$,

\adfLgap
$(0,1,2,3,5,6,11,9,13,22,18,19,\adfsplit 23,24,39,40,43,48,53,62,68,10,27,26)_{\adfGI}$,

\adfLgap
$(0,1,2,3,5,6,11,9,13,22,18,19,\adfsplit 23,24,39,40,41,51,55,66,61,7,31,26)_{\adfGJ}$,

\adfLgap
$(0,1,2,3,5,6,11,9,13,14,18,19,\adfsplit 25,28,26,34,32,53,51,45,56,68,69,22)_{\adfGK}$

\noindent and

$(0,1,2,3,5,6,11,9,13,14,18,19,\adfsplit 25,28,26,38,43,66,56,35,59,63,4,16)_{\adfGL}$

\noindent under the action of the mapping $x \mapsto x + 1$ (mod 73).

\adfVfy{73, \{\{73,73,1\}\}, -1, -1, -1} 

Let the vertex set of $K_{145}$ be $Z_{145}$. The decompositions consist of

$(0,1,2,3,5,6,11,13,15,16,24,30,\adfsplit 28,44,4,47,7,8,26,9,49,33,58,76)_{\adfGa}$,

$(0,28,30,31,61,62,96,99,101,70,5,6,\adfsplit 1,42,49,95,63,124,98,23,11,114,53,71)_{\adfGa}$,

\adfLgap
$(0,1,2,3,5,6,11,13,15,16,24,25,\adfsplit 28,44,42,4,9,7,8,26,32,51,78,62)_{\adfGb}$,

$(0,28,31,32,61,63,99,103,108,71,4,10,\adfsplit 2,47,112,56,62,129,29,106,120,14,78,41)_{\adfGb}$,

\adfLgap
$(122,41,16,65,93,90,83,99,25,137,55,9,\adfsplit 8,105,95,47,1,2,3,5,0,11,17,24)_{\adfGc}$,

$(0,1,4,5,15,16,33,35,38,25,50,3,\adfsplit 7,6,49,87,77,108,104,82,90,138,24,48)_{\adfGc}$,

\adfLgap
$(0,1,2,3,5,6,11,13,15,16,24,25,\adfsplit 28,31,42,4,44,7,50,54,74,67,101,8)_{\adfGd}$,

$(0,28,29,30,59,60,93,95,92,66,5,1,\adfsplit 2,3,105,100,55,89,39,41,88,31,11,111)_{\adfGd}$,

\adfLgap
$(5,17,106,0,110,72,81,117,83,136,85,61,\adfsplit 69,10,142,4,12,2,86,134,8,11,21,25)_{\adfGe}$,

$(0,2,3,14,18,20,44,41,46,45,50,10,\adfsplit 4,88,6,83,51,134,91,25,142,101,38,68)_{\adfGe}$,

\adfLgap
$(114,9,69,140,130,39,47,40,88,99,135,117,\adfsplit 67,75,53,23,35,78,49,58,93,15,126,132)_{\adfGf}$,

$(0,2,5,6,9,11,19,22,21,25,27,50,\adfsplit 1,57,80,87,24,107,4,15,102,103,42,43)_{\adfGf}$,

\adfLgap
$(0,1,2,3,5,6,11,13,15,16,14,30,\adfsplit 28,31,34,47,4,7,53,54,8,72,81,101)_{\adfGg}$,

$(0,28,30,31,60,62,95,99,98,67,1,4,\adfsplit 6,3,106,45,48,108,103,49,44,101,123,16)_{\adfGg}$,

\adfLgap
$(107,12,20,21,57,10,81,102,22,138,84,2,\adfsplit 14,121,6,41,113,4,143,18,95,120,56,140)_{\adfGh}$,

$(0,4,5,6,12,15,24,32,42,27,33,2,\adfsplit 3,1,65,102,51,101,43,45,76,142,21,105)_{\adfGh}$,

\adfLgap
$(0,1,2,3,5,6,11,13,15,16,14,30,\adfsplit 28,31,33,48,53,4,7,54,56,70,82,97)_{\adfGi}$,

$(0,28,30,31,60,61,95,99,101,67,1,3,\adfsplit 2,7,4,48,62,44,93,100,47,92,130,21)_{\adfGi}$,

\adfLgap
$(51,125,128,43,67,63,38,9,48,122,56,113,\adfsplit 101,33,44,124,89,20,138,90,42,93,121,65)_{\adfGj}$,

$(0,1,4,6,11,14,25,34,33,22,32,49,\adfsplit 3,67,95,58,96,10,108,45,20,2,52,104)_{\adfGj}$,

\adfLgap
$(25,122,34,11,4,6,49,99,75,53,67,29,\adfsplit 23,22,142,18,86,114,127,69,44,112,91,74)_{\adfGk}$,

$(0,1,2,3,5,7,12,15,21,14,24,35,\adfsplit 39,46,61,45,58,102,6,78,97,141,83,137)_{\adfGk}$,

\adfLgap
$(117,64,104,116,42,114,113,36,99,80,126,88,\adfsplit 63,55,29,43,128,50,14,9,127,89,95,85)_{\adfGl}$,

$(0,3,4,5,10,11,20,22,24,25,26,43,\adfsplit 46,50,54,83,124,109,108,132,139,93,9,57)_{\adfGl}$,

\adfLgap
$(91,61,10,84,43,137,19,6,33,123,44,92,\adfsplit 51,120,100,28,8,88,94,93,15,16,113,114)_{\adfGm}$,

$(0,1,2,3,6,7,17,24,19,28,34,45,\adfsplit 53,50,77,63,101,132,86,141,120,87,4,67)_{\adfGm}$,

\adfLgap
$(134,32,14,58,85,129,23,139,130,133,75,141,\adfsplit 107,86,41,17,40,124,8,125,33,26,37,119)_{\adfGn}$,

$(0,2,3,8,12,13,25,18,38,44,57,30,\adfsplit 69,71,1,54,87,117,83,133,111,26,43,68)_{\adfGn}$,

\adfLgap
$(66,125,73,74,21,18,25,0,24,130,79,6,\adfsplit 45,133,141,12,135,98,30,83,108,78,131,100)_{\adfGo}$,

$(0,1,5,9,14,15,30,23,33,51,43,41,\adfsplit 50,4,88,6,82,132,8,49,99,31,68,107)_{\adfGo}$,

\adfLgap
$(0,1,2,3,5,6,11,9,13,22,21,19,\adfsplit 23,39,40,4,41,43,7,27,12,65,90,55)_{\adfGp}$,

$(0,27,30,32,61,63,98,71,103,137,2,3,\adfsplit 1,45,48,94,55,80,121,16,143,106,26,62)_{\adfGp}$,

\adfLgap
$(0,1,2,3,5,6,11,9,13,22,21,19,\adfsplit 23,29,39,57,40,45,46,4,72,7,8,70)_{\adfGq}$,

$(0,27,28,29,57,58,95,62,97,130,1,2,\adfsplit 3,4,47,86,56,54,100,107,53,132,128,42)_{\adfGq}$,

\adfLgap
$(142,130,88,24,71,113,89,6,105,109,36,98,\adfsplit 102,54,51,59,19,99,91,86,134,112,90,119)_{\adfGr}$,

$(0,2,4,5,8,10,18,20,27,42,41,35,\adfsplit 46,56,77,79,135,136,82,97,139,125,66,30)_{\adfGr}$,

\adfLgap
$(89,100,36,52,43,14,51,83,55,65,9,18,\adfsplit 81,11,143,40,91,39,131,74,7,121,73,6)_{\adfGs}$,

$(0,6,10,16,23,25,46,30,51,70,57,55,\adfsplit 4,2,1,108,96,9,85,42,41,61,141,84)_{\adfGs}$,

\adfLgap
$(112,139,142,61,34,76,99,101,107,67,91,132,\adfsplit 4,135,74,56,134,78,10,124,133,2,33,95)_{\adfGt}$,

$(0,3,4,5,12,13,23,16,26,58,29,35,\adfsplit 49,53,1,60,61,2,115,141,111,22,80,45)_{\adfGt}$,

\adfLgap
$(72,20,125,26,54,1,50,41,133,86,128,8,\adfsplit 47,131,4,109,75,79,126,73,53,101,112,30)_{\adfGu}$,

$(0,1,2,9,11,13,24,17,36,40,30,33,\adfsplit 65,71,70,57,107,16,128,134,122,92,34,39)_{\adfGu}$,

\adfLgap
$(103,45,19,86,121,44,47,9,42,27,142,139,\adfsplit 77,123,136,82,99,54,134,110,74,124,97,10)_{\adfGv}$,

$(0,4,5,6,11,13,23,18,35,38,29,43,\adfsplit 57,66,76,63,84,109,20,25,14,117,27,87)_{\adfGv}$,

\adfLgap
$(0,1,2,3,5,6,11,9,13,22,18,19,\adfsplit 27,33,39,34,40,4,55,7,65,60,10,92)_{\adfGw}$,

$(0,28,29,33,62,63,98,66,100,137,4,1,\adfsplit 2,106,47,56,49,91,19,96,117,14,93,26)_{\adfGw}$,

\adfLgap
$(132,119,38,85,58,10,67,64,78,66,133,98,\adfsplit 118,88,49,62,87,137,16,131,6,44,39,56)_{\adfGx}$,

$(0,2,3,4,7,8,28,11,34,38,23,24,\adfsplit 66,60,69,101,103,76,125,27,16,39,57,134)_{\adfGx}$,

\adfLgap
$(0,1,2,3,5,6,11,9,13,22,18,19,\adfsplit 27,24,38,41,40,43,63,65,4,66,8,92)_{\adfGy}$,

$(0,28,29,30,59,60,92,66,97,131,1,4,\adfsplit 2,5,61,99,106,80,47,58,115,41,120,12)_{\adfGy}$,

\adfLgap
$(0,1,2,3,5,6,11,9,13,22,18,19,\adfsplit 27,24,38,41,40,43,62,87,4,10,66,117)_{\adfGz}$,

$(0,27,28,29,63,64,106,68,104,67,9,3,\adfsplit 8,2,109,60,51,72,103,59,37,7,105,118)_{\adfGz}$,

\adfLgap
$(38,59,6,110,55,3,98,16,50,63,27,130,\adfsplit 53,54,12,101,34,111,89,113,105,24,71,103)_{\adfGA}$,

$(0,1,6,7,9,10,22,17,32,48,27,37,\adfsplit 61,39,81,55,85,67,130,119,120,137,36,42)_{\adfGA}$,

\adfLgap
$(0,1,2,3,5,6,11,9,13,22,18,19,\adfsplit 27,24,38,35,45,41,4,59,7,8,70,71)_{\adfGB}$,

$(0,28,29,32,62,64,99,68,103,137,3,4,\adfsplit 2,9,45,46,81,50,53,102,92,100,12,5)_{\adfGB}$,

\adfLgap
$(8,18,14,110,87,35,138,67,48,66,31,126,\adfsplit 102,68,95,88,84,6,39,136,29,75,38,127)_{\adfGC}$,

$(0,2,3,5,9,10,22,14,21,36,25,34,\adfsplit 51,41,49,84,116,134,102,1,79,136,20,61)_{\adfGC}$,

\adfLgap
$(140,124,62,10,102,24,53,98,33,85,79,30,\adfsplit 7,55,50,109,41,59,11,21,67,1,32,115)_{\adfGD}$,

$(0,1,2,3,8,9,20,12,22,53,34,35,\adfsplit 47,40,81,105,73,143,139,141,21,120,58,78)_{\adfGD}$,

\adfLgap
$(0,1,2,3,5,6,11,9,13,22,18,19,\adfsplit 32,24,40,39,57,41,4,61,7,31,66,91)_{\adfGE}$,

$(0,29,31,32,63,64,107,71,104,68,13,1,\adfsplit 113,117,37,115,67,78,49,28,86,15,139,103)_{\adfGE}$,

\adfLgap
$(0,1,2,3,5,6,11,9,13,22,18,19,\adfsplit 32,24,40,39,36,41,4,62,7,10,86,44)_{\adfGF}$,

$(0,25,29,30,58,60,110,71,101,69,104,105,\adfsplit 6,1,47,8,44,16,100,83,81,99,140,23)_{\adfGF}$,

\adfLgap
$(115,106,13,130,77,79,140,55,38,109,33,63,\adfsplit 52,116,129,51,10,22,104,35,21,15,54,26)_{\adfGG}$,

$(0,1,2,6,9,11,26,13,41,62,33,37,\adfsplit 70,47,92,102,71,110,140,107,31,24,121,44)_{\adfGG}$,

\adfLgap
$(51,30,107,25,55,15,63,77,84,39,96,33,\adfsplit 83,81,126,99,135,94,86,59,138,36,14,61)_{\adfGH}$,

$(0,1,2,3,8,10,18,14,31,35,29,30,\adfsplit 60,48,82,76,134,138,79,117,113,16,7,61)_{\adfGH}$,

\adfLgap
$(56,129,21,79,81,142,86,92,29,39,14,40,\adfsplit 99,120,35,59,13,9,58,118,20,25,112,143)_{\adfGI}$,

$(0,3,4,8,12,13,23,18,38,44,28,29,\adfsplit 42,45,1,82,91,113,119,102,6,25,81,73)_{\adfGI}$,

\adfLgap
$(0,1,2,3,5,6,11,9,13,22,18,19,\adfsplit 23,24,39,40,41,43,44,60,7,4,70,88)_{\adfGJ}$,

$(0,29,30,31,61,63,99,69,102,134,1,2,\adfsplit 4,9,47,46,55,96,75,95,93,133,44,25)_{\adfGJ}$,

\adfLgap
$(6,53,130,28,121,115,9,63,76,23,43,11,\adfsplit 37,132,12,104,0,1,2,3,5,7,14,20)_{\adfGK}$,

$(0,3,6,7,15,19,32,25,42,36,44,49,\adfsplit 1,79,64,110,126,88,134,114,51,22,21,85)_{\adfGK}$

\noindent and

$(69,86,13,120,46,112,76,133,82,141,66,94,\adfsplit 84,80,129,72,1,0,3,4,7,16,17,23)_{\adfGL}$,

$(0,3,9,11,17,18,33,32,45,38,47,49,\adfsplit 1,80,69,90,138,89,136,93,51,140,10,65)_{\adfGL}$

\noindent under the action of the mapping $x \mapsto x + 1$ (mod 145).\eproof

\adfVfy{145, \{\{145,145,1\}\}, -1, -1, -1} 


\section{Proof of Lemma~\ref{lem:Snark24 multipartite}}
\label{app:Snark24 multipartite}


\noindent\textbf{Proof of Lemma~\ref{lem:Snark24 multipartite}.}
Recall that the lemma asserts the existence of decompositions of the complete multipartite graphs
$K_{12,12,12}$, $K_{24,24,15}$, $K_{72,72,63}$, $K_{24,24,24,24}$ and $K_{24,24,24,21}$
for each of the thirty-eight 24-vertex non-trivial snarks.

Let the vertex set of $K_{12,12,12}$ be $Z_{36}$ partitioned according to residue classes modulo 3.
The decompositions consist of

$(0,1,2,4,3,5,10,6,11,9,19,20,\adfsplit 29,13,27,7,15,30,17,18,28,35,16,23)_{\adfGa}$,

\adfLgap
$(0,1,2,4,3,5,10,6,11,9,19,18,\adfsplit 17,16,31,35,27,33,20,12,7,32,22,26)_{\adfGb}$,

\adfLgap
$(0,1,2,4,3,5,10,6,11,9,19,8,\adfsplit 32,24,34,33,15,12,7,17,29,13,28,26)_{\adfGc}$,

\adfLgap
$(0,1,2,4,3,5,10,6,11,9,19,18,\adfsplit 17,30,23,8,31,27,33,16,13,28,14,26)_{\adfGd}$,

\adfLgap
$(0,1,2,4,3,5,10,6,11,9,7,20,\adfsplit 23,33,31,15,30,22,32,28,13,17,12,8)_{\adfGe}$,

\adfLgap
$(0,1,2,4,3,5,10,6,11,9,12,20,\adfsplit 31,21,28,25,8,35,29,14,24,27,13,16)_{\adfGf}$,

\adfLgap
$(0,1,2,4,3,5,10,6,11,9,7,20,\adfsplit 34,16,23,27,12,15,17,35,33,14,22,28)_{\adfGg}$,

\adfLgap
$(0,1,2,4,3,5,10,6,11,9,12,20,\adfsplit 16,30,31,28,15,27,35,14,22,29,17,13)_{\adfGh}$,

\adfLgap
$(0,1,2,4,3,5,10,6,11,9,7,20,\adfsplit 29,16,19,30,27,24,12,17,32,8,13,31)_{\adfGi}$,

\adfLgap
$(0,1,2,4,3,5,10,6,11,9,7,30,\adfsplit 26,34,23,35,18,16,21,25,8,20,31,33)_{\adfGj}$,

\adfLgap
$(0,1,2,4,3,5,10,6,11,9,7,33,\adfsplit 17,27,23,25,15,20,31,24,8,13,14,34)_{\adfGk}$,

\adfLgap
$(0,1,2,4,3,5,10,6,11,9,7,23,\adfsplit 17,33,32,28,18,22,12,24,25,14,34,8)_{\adfGl}$,

\adfLgap
$(0,1,2,4,3,5,10,6,11,9,7,35,\adfsplit 25,16,19,31,15,20,12,8,26,27,17,33)_{\adfGm}$,

\adfLgap
$(0,1,2,4,3,5,10,6,11,15,16,7,\adfsplit 25,8,30,27,35,21,14,19,17,18,34,32)_{\adfGn}$,

\adfLgap
$(0,1,2,4,3,5,10,6,11,15,16,18,\adfsplit 17,30,8,32,22,31,9,12,34,35,14,25)_{\adfGo}$,

\adfLgap
$(0,1,2,4,3,5,10,6,11,15,16,12,\adfsplit 25,33,35,14,17,26,9,34,7,18,31,20)_{\adfGp}$,

\adfLgap
$(0,1,2,4,3,5,10,6,11,15,16,12,\adfsplit 22,25,34,29,32,20,9,33,26,13,21,23)_{\adfGq}$,

\adfLgap
$(0,1,2,4,3,5,10,6,11,15,16,7,\adfsplit 17,34,8,33,23,31,30,18,29,25,9,35)_{\adfGr}$,

\adfLgap
$(0,1,2,4,3,5,10,6,11,15,16,12,\adfsplit 22,34,26,17,32,29,8,30,27,31,18,13)_{\adfGs}$,

\adfLgap
$(0,1,2,4,3,5,10,6,11,15,12,7,\adfsplit 34,26,28,17,33,35,9,27,16,8,13,32)_{\adfGt}$,

\adfLgap
$(0,1,2,4,3,5,10,6,11,15,7,12,\adfsplit 22,29,26,34,32,9,28,27,14,8,18,19)_{\adfGu}$,

\adfLgap
$(0,1,2,4,3,5,10,6,11,15,7,24,\adfsplit 22,29,8,21,23,17,16,18,13,31,12,32)_{\adfGv}$,

\adfLgap
$(0,1,2,4,3,5,10,6,11,15,7,12,\adfsplit 34,26,28,23,31,8,21,14,9,33,32,19)_{\adfGw}$,

\adfLgap
$(0,1,2,4,3,5,10,6,11,15,7,12,\adfsplit 25,32,31,33,34,17,13,8,18,14,29,30)_{\adfGx}$,

\adfLgap
$(0,1,2,4,3,5,10,6,11,15,7,12,\adfsplit 16,19,26,27,34,29,32,8,13,35,9,30)_{\adfGy}$,

\adfLgap
$(0,1,2,4,3,5,10,6,11,15,7,16,\adfsplit 28,20,8,30,32,24,9,25,22,29,12,35)_{\adfGz}$,

\adfLgap
$(0,1,2,4,3,5,10,6,11,15,7,12,\adfsplit 22,28,20,24,29,26,9,14,21,16,31,32)_{\adfGA}$,

\adfLgap
$(0,1,2,4,3,5,10,6,11,15,7,19,\adfsplit 22,25,8,21,35,9,30,20,31,32,18,14)_{\adfGB}$,

\adfLgap
$(0,1,2,4,3,5,10,6,11,15,12,7,\adfsplit 22,19,31,26,32,8,27,14,9,17,13,30)_{\adfGC}$,

\adfLgap
$(0,1,2,4,3,5,10,6,11,15,7,12,\adfsplit 26,16,33,32,23,34,13,24,9,22,8,14)_{\adfGD}$,

\adfLgap
$(0,1,2,4,3,5,10,6,11,15,7,16,\adfsplit 23,34,12,20,33,17,9,24,32,13,22,35)_{\adfGE}$,

\adfLgap
$(0,1,2,4,3,5,10,6,11,15,7,18,\adfsplit 23,25,33,8,21,17,13,12,19,31,32,14)_{\adfGF}$,

\adfLgap
$(0,1,2,4,3,5,10,6,11,15,7,12,\adfsplit 23,29,34,25,21,31,18,9,35,20,13,32)_{\adfGG}$,

\adfLgap
$(0,1,2,4,3,5,10,6,11,15,7,12,\adfsplit 17,32,27,29,30,28,31,9,13,8,19,35)_{\adfGH}$,

\adfLgap
$(0,1,2,4,3,5,10,6,11,15,7,19,\adfsplit 17,25,28,8,33,35,27,12,31,14,9,26)_{\adfGI}$,

\adfLgap
$(0,1,2,4,3,5,10,6,11,15,7,12,\adfsplit 23,22,33,20,30,34,32,8,19,13,29,9)_{\adfGJ}$,

\adfLgap
$(0,1,2,4,3,5,10,6,11,9,7,30,\adfsplit 23,25,19,32,26,33,14,16,18,35,27,34)_{\adfGK}$

\noindent and

$(0,1,2,4,3,5,10,6,11,9,7,30,\adfsplit 23,16,34,32,31,18,8,13,33,35,14,15)_{\adfGL}$

\noindent under the action of the mapping $x \mapsto x + 3$ (mod 36).

\adfVfy{36, \{\{36,12,3\}\}, -1, \{\{12,\{0,1,2\}\}\}, -1} 

Let the vertex set of $K_{24,24,15}$ be $\{0, 1, \dots, 62\}$ partitioned into
$\{2j + i: j = 0, 1, \dots, 23\}$, $i = 0, 1$, and $\{48, 49, \dots, 62\}$.
The decompositions consist of

$(27,22,53,48,9,33,26,12,0,16,51,23,\adfsplit 50,37,20,54,34,10,7,4,60,5,59,47)_{\adfGa}$,

$(46,15,55,21,42,4,35,37,53,26,59,61,\adfsplit 40,5,1,8,9,14,49,52,10,6,39,23)_{\adfGa}$,

$(13,60,36,53,16,0,7,3,47,46,5,61,\adfsplit 4,57,52,9,33,44,38,62,56,18,19,23)_{\adfGa}$,

\adfLgap
$(29,30,53,12,57,5,45,18,39,21,38,26,\adfsplit 58,24,56,23,14,51,37,55,59,3,44,16)_{\adfGb}$,

$(38,15,55,43,34,53,25,0,59,2,23,60,\adfsplit 48,21,11,12,39,61,40,41,52,4,35,46)_{\adfGb}$,

$(40,39,25,54,38,49,50,22,1,20,43,18,\adfsplit 29,28,45,51,57,62,0,47,13,34,3,8)_{\adfGb}$,

\adfLgap
$(51,45,15,8,16,6,57,53,3,47,56,29,\adfsplit 10,0,34,40,21,13,54,58,7,52,28,38)_{\adfGc}$,

$(12,57,17,31,19,21,20,8,42,55,36,58,\adfsplit 23,27,9,5,16,51,54,10,59,22,39,37)_{\adfGc}$,

$(0,61,60,59,13,38,34,9,26,35,53,50,\adfsplit 8,23,24,7,39,29,30,22,1,2,55,62)_{\adfGc}$,

\adfLgap
$(2,45,59,3,54,18,20,35,14,16,1,25,\adfsplit 60,23,52,48,38,40,19,61,44,17,43,6)_{\adfGd}$,

$(47,26,36,49,19,29,57,54,38,43,50,27,\adfsplit 37,15,58,9,20,12,30,48,56,28,33,45)_{\adfGd}$,

$(37,38,50,46,52,58,6,36,43,17,42,39,\adfsplit 19,0,61,55,20,4,28,1,51,2,13,57)_{\adfGd}$,

\adfLgap
$(32,48,37,21,18,13,53,26,50,52,34,11,\adfsplit 54,0,4,27,43,51,5,42,38,58,44,39)_{\adfGe}$,

$(49,0,17,31,13,47,44,18,16,36,30,55,\adfsplit 7,51,52,3,60,6,46,24,21,27,54,57)_{\adfGe}$,

$(54,5,44,37,12,20,31,3,58,18,1,52,\adfsplit 38,47,9,56,21,8,55,30,26,61,19,33)_{\adfGe}$,

\adfLgap
$(33,34,52,50,23,35,16,10,26,4,49,31,\adfsplit 48,45,27,12,1,30,58,54,44,18,47,41)_{\adfGf}$,

$(6,45,33,15,60,2,20,56,28,52,5,37,\adfsplit 13,51,18,4,36,39,49,47,29,57,40,8)_{\adfGf}$,

$(29,52,59,20,47,2,42,15,51,9,61,53,\adfsplit 50,7,34,8,11,4,55,19,38,26,33,58)_{\adfGf}$,

\adfLgap
$(31,58,36,0,28,18,13,7,57,51,45,44,\adfsplit 55,5,4,49,25,16,30,48,19,10,42,29)_{\adfGg}$,

$(57,9,46,27,22,48,3,39,56,10,37,52,\adfsplit 16,38,59,28,53,17,35,55,19,60,6,2)_{\adfGg}$,

$(14,52,5,25,21,19,51,8,18,4,56,43,\adfsplit 7,49,41,42,20,0,1,3,36,15,59,60)_{\adfGg}$,

\adfLgap
$(29,42,18,26,55,48,27,47,25,52,32,24,\adfsplit 12,0,5,37,54,49,1,7,61,38,28,51)_{\adfGh}$,

$(18,23,49,37,32,40,17,39,53,56,21,16,\adfsplit 20,26,46,28,43,57,31,27,48,52,25,42)_{\adfGh}$,

$(61,43,23,42,30,59,60,46,25,11,33,38,\adfsplit 5,34,16,62,35,50,48,13,27,55,4,0)_{\adfGh}$,

\adfLgap
$(46,19,33,54,44,12,48,61,1,47,3,0,\adfsplit 18,10,22,60,53,52,55,21,15,5,32,2)_{\adfGi}$,

$(4,17,60,33,53,8,9,11,18,36,51,59,\adfsplit 58,50,52,42,43,30,3,5,24,27,54,10)_{\adfGi}$,

$(52,29,40,34,12,49,27,11,33,41,24,16,\adfsplit 10,56,61,9,19,39,62,42,43,20,30,21)_{\adfGi}$,

\adfLgap
$(39,50,28,46,9,42,59,21,53,54,37,25,\adfsplit 58,15,14,23,24,6,8,22,60,48,43,19)_{\adfGj}$,

$(33,32,55,61,7,29,2,12,4,38,40,51,\adfsplit 39,47,21,50,1,28,52,14,10,49,24,43)_{\adfGj}$,

$(59,35,24,41,28,4,58,9,52,26,43,17,\adfsplit 36,7,62,13,61,14,51,57,46,6,27,5)_{\adfGj}$,

\adfLgap
$(53,29,35,9,60,6,14,49,44,28,21,25,\adfsplit 41,59,56,47,34,52,58,55,40,3,20,33)_{\adfGk}$,

$(45,12,6,61,58,56,47,3,46,36,34,32,\adfsplit 51,53,41,59,60,35,22,31,39,26,5,30)_{\adfGk}$,

$(29,54,34,28,30,2,15,5,57,3,27,62,\adfsplit 8,13,14,12,16,11,55,52,24,40,47,41)_{\adfGk}$,

\adfLgap
$(45,61,36,58,41,28,55,39,6,11,13,23,\adfsplit 50,19,32,4,53,57,16,35,51,2,47,12)_{\adfGl}$,

$(49,16,29,43,11,13,22,34,44,4,8,54,\adfsplit 27,21,45,25,60,28,52,26,18,56,10,48)_{\adfGl}$,

$(12,47,25,53,26,51,55,2,3,9,15,34,\adfsplit 21,54,10,59,6,1,62,52,19,36,37,38)_{\adfGl}$,

\adfLgap
$(52,34,25,39,35,23,16,42,60,12,0,51,\adfsplit 15,43,7,55,5,49,58,6,4,59,22,33)_{\adfGm}$,

$(11,44,12,46,37,54,56,19,51,61,41,4,\adfsplit 14,6,7,23,48,52,33,39,55,57,36,28)_{\adfGm}$,

$(49,20,17,30,55,9,14,32,25,29,44,47,\adfsplit 62,10,6,42,58,45,38,41,7,53,1,24)_{\adfGm}$,

\adfLgap
$(35,22,8,18,48,61,13,56,33,54,5,4,\adfsplit 16,29,42,27,25,44,17,36,12,55,57,9)_{\adfGn}$,

$(0,52,1,54,43,22,12,8,10,15,58,29,\adfsplit 5,51,31,49,2,47,24,19,3,61,60,42)_{\adfGn}$,

$(16,54,27,35,20,38,48,18,53,42,37,39,\adfsplit 9,55,50,2,30,14,24,23,52,17,11,62)_{\adfGn}$,

\adfLgap
$(56,6,22,21,59,35,47,58,36,53,4,18,\adfsplit 20,7,3,61,60,17,24,8,15,10,52,13)_{\adfGo}$,

$(20,50,7,37,29,22,14,12,57,59,23,33,\adfsplit 27,24,8,34,56,4,1,43,9,44,48,55)_{\adfGo}$,

$(46,59,25,52,17,31,56,44,16,43,49,10,\adfsplit 3,37,53,29,6,60,14,38,22,57,9,27)_{\adfGo}$,

\adfLgap
$(28,43,31,19,32,8,61,48,55,40,41,33,\adfsplit 20,17,51,36,52,50,30,11,0,38,29,49)_{\adfGp}$,

$(36,31,43,9,18,58,60,49,59,46,29,24,\adfsplit 8,17,57,16,13,7,41,47,38,55,10,56)_{\adfGp}$,

$(45,26,18,4,57,35,23,53,21,52,34,6,\adfsplit 7,14,39,29,56,30,20,42,33,31,54,58)_{\adfGp}$,

\adfLgap
$(16,17,31,7,53,18,32,50,46,61,5,37,\adfsplit 35,43,34,51,40,56,4,0,58,9,29,59)_{\adfGq}$,

$(37,12,14,46,31,21,52,45,56,15,16,54,\adfsplit 4,20,28,39,22,53,11,47,60,48,33,26)_{\adfGq}$,

$(45,24,49,12,37,52,22,7,57,55,6,35,\adfsplit 54,41,42,1,30,18,3,60,14,62,43,0)_{\adfGq}$,

\adfLgap
$(45,38,48,16,7,51,24,10,55,9,25,22,\adfsplit 21,14,54,52,41,26,34,27,43,36,50,56)_{\adfGr}$,

$(39,42,20,52,33,48,56,50,16,43,13,4,\adfsplit 30,19,49,36,3,35,15,28,54,25,14,58)_{\adfGr}$,

$(23,28,46,12,7,54,58,25,60,47,42,13,\adfsplit 2,29,56,45,0,52,43,30,20,40,21,62)_{\adfGr}$,

\adfLgap
$(50,43,36,13,24,12,31,9,59,10,56,57,\adfsplit 48,20,5,29,11,26,25,2,30,54,34,35)_{\adfGs}$,

$(23,58,48,59,28,44,43,29,3,46,57,61,\adfsplit 40,16,5,27,19,22,41,50,18,60,30,35)_{\adfGs}$,

$(27,34,24,61,55,5,21,3,12,6,52,20,\adfsplit 54,29,48,30,62,19,0,13,10,28,41,56)_{\adfGs}$,

\adfLgap
$(31,53,8,46,3,12,58,33,61,5,47,55,\adfsplit 0,4,22,21,2,49,17,57,41,56,37,6)_{\adfGt}$,

$(22,41,25,35,18,60,50,61,59,24,23,34,\adfsplit 38,16,54,39,15,47,9,46,44,42,52,53)_{\adfGt}$,

$(26,53,23,5,17,42,22,32,16,59,31,11,\adfsplit 37,57,0,51,55,40,29,27,39,20,62,12)_{\adfGt}$,

\adfLgap
$(11,44,58,40,5,1,32,17,56,9,59,54,\adfsplit 26,46,12,15,0,19,55,41,2,51,47,61)_{\adfGu}$,

$(37,34,36,26,21,50,55,43,59,15,33,38,\adfsplit 42,18,40,47,53,25,57,52,48,8,35,30)_{\adfGu}$,

$(53,44,33,36,31,29,26,34,49,47,57,52,\adfsplit 15,19,4,46,11,32,51,55,21,41,14,8)_{\adfGu}$,

\adfLgap
$(13,55,6,61,40,47,33,43,31,8,22,50,\adfsplit 2,52,27,25,17,15,20,48,44,30,53,19)_{\adfGv}$,

$(29,46,58,28,7,17,40,1,59,37,52,0,\adfsplit 36,38,49,33,48,39,51,42,34,60,4,9)_{\adfGv}$,

$(54,33,14,1,0,44,56,19,52,10,3,57,\adfsplit 18,49,43,61,15,35,38,22,36,9,4,60)_{\adfGv}$,

\adfLgap
$(31,12,48,16,56,17,29,47,43,60,58,59,\adfsplit 14,38,10,40,11,52,9,49,50,1,8,18)_{\adfGw}$,

$(58,42,44,11,56,33,25,31,34,46,8,26,\adfsplit 51,61,9,49,1,4,37,40,12,5,2,35)_{\adfGw}$,

$(21,55,60,0,38,28,11,3,7,57,39,19,\adfsplit 10,2,20,52,49,25,48,62,35,1,4,30)_{\adfGw}$,

\adfLgap
$(10,27,17,45,53,48,32,61,18,56,19,34,\adfsplit 47,28,5,6,50,37,60,30,24,15,41,57)_{\adfGx}$,

$(3,46,40,59,52,45,19,11,18,24,14,49,\adfsplit 56,50,54,39,5,30,15,16,27,8,28,31)_{\adfGx}$,

$(54,19,29,18,46,34,53,36,9,33,23,51,\adfsplit 48,62,28,20,6,0,16,1,50,17,25,52)_{\adfGx}$,

\adfLgap
$(52,41,40,34,14,22,47,60,9,53,43,11,\adfsplit 20,30,8,57,55,2,3,29,1,44,61,51)_{\adfGy}$,

$(42,45,7,58,10,44,54,49,21,13,41,47,\adfsplit 14,35,57,12,18,30,11,61,16,33,50,31)_{\adfGy}$,

$(57,35,34,12,42,10,25,60,53,56,15,49,\adfsplit 1,11,20,48,5,8,36,32,54,28,47,17)_{\adfGy}$,

\adfLgap
$(32,58,55,3,2,6,21,25,51,42,47,29,\adfsplit 22,14,53,59,30,8,57,35,5,15,54,56)_{\adfGz}$,

$(7,20,18,30,5,13,55,31,51,35,22,52,\adfsplit 0,40,50,37,34,39,49,45,36,6,3,57)_{\adfGz}$,

$(39,50,10,48,20,18,45,11,23,44,25,13,\adfsplit 0,12,57,58,32,24,27,56,28,41,49,37)_{\adfGz}$,

\adfLgap
$(9,4,28,40,50,49,37,7,15,59,2,31,\adfsplit 56,26,11,60,58,12,23,41,44,19,61,57)_{\adfGA}$,

$(52,0,22,47,31,39,50,54,14,36,44,2,\adfsplit 17,26,58,51,27,7,24,1,13,34,42,57)_{\adfGA}$,

$(11,24,49,12,21,13,0,33,55,25,26,30,\adfsplit 6,34,48,37,27,56,58,62,8,46,9,41)_{\adfGA}$,

\adfLgap
$(17,48,34,56,2,37,41,9,13,42,61,54,\adfsplit 8,6,15,4,33,35,58,7,60,57,44,36)_{\adfGB}$,

$(47,2,49,30,61,27,20,34,19,39,40,51,\adfsplit 60,1,58,35,28,22,41,8,38,59,9,55)_{\adfGB}$,

$(21,24,4,12,37,31,60,35,57,23,18,40,\adfsplit 6,30,62,5,15,54,47,52,9,10,48,32)_{\adfGB}$,

\adfLgap
$(4,35,43,31,50,34,46,26,24,59,9,19,\adfsplit 61,57,41,44,30,55,23,27,32,0,1,53)_{\adfGC}$,

$(13,55,50,34,18,36,5,30,33,61,21,23,\adfsplit 49,3,26,32,46,4,48,22,53,31,9,57)_{\adfGC}$,

$(48,1,38,29,6,46,47,21,57,32,11,37,\adfsplit 0,40,4,51,61,49,39,62,5,59,15,12)_{\adfGC}$,

\adfLgap
$(17,40,59,30,60,19,3,25,31,48,61,36,\adfsplit 32,10,56,22,57,27,29,47,39,46,58,28)_{\adfGD}$,

$(26,60,33,9,23,46,28,50,14,53,31,11,\adfsplit 16,10,39,8,6,49,57,29,51,45,18,40)_{\adfGD}$,

$(47,56,26,12,14,21,25,54,49,50,4,36,\adfsplit 7,0,13,5,53,9,10,29,44,3,52,57)_{\adfGD}$,

\adfLgap
$(2,45,19,31,54,59,30,32,20,53,18,3,\adfsplit 1,60,33,11,28,52,26,23,55,27,13,61)_{\adfGE}$,

$(52,22,45,30,54,50,36,26,39,55,40,42,\adfsplit 47,27,12,11,57,41,59,14,13,24,8,5)_{\adfGE}$,

$(16,7,56,52,22,30,33,14,29,36,45,58,\adfsplit 51,11,8,61,19,9,20,6,1,24,53,48)_{\adfGE}$,

\adfLgap
$(16,54,61,33,15,34,38,27,10,58,5,60,\adfsplit 44,8,19,41,20,25,57,49,0,50,45,35)_{\adfGF}$,

$(23,6,38,44,27,19,57,56,60,17,20,24,\adfsplit 33,28,35,61,1,58,42,3,34,18,59,45)_{\adfGF}$,

$(18,52,7,57,8,33,58,4,45,16,46,26,\adfsplit 43,19,53,21,49,41,42,56,22,40,3,50)_{\adfGF}$,

\adfLgap
$(17,26,28,48,11,53,57,13,44,39,10,1,\adfsplit 55,59,45,38,31,46,0,41,52,51,34,27)_{\adfGG}$,

$(26,31,1,50,36,55,57,40,17,19,5,14,\adfsplit 56,51,2,48,6,27,37,35,4,59,43,18)_{\adfGG}$,

$(52,13,40,18,34,6,3,15,54,51,1,37,\adfsplit 28,55,8,32,24,20,21,4,59,23,31,58)_{\adfGG}$,

\adfLgap
$(14,1,53,13,49,16,29,20,22,0,55,27,\adfsplit 15,9,5,60,33,26,51,2,52,36,56,32)_{\adfGH}$,

$(14,43,55,11,4,50,1,0,53,46,26,39,\adfsplit 41,56,42,31,58,61,36,45,34,54,57,5)_{\adfGH}$,

$(29,6,4,59,1,33,58,39,44,47,12,22,\adfsplit 62,34,3,16,31,43,52,27,11,42,28,49)_{\adfGH}$,

\adfLgap
$(15,40,59,50,17,58,28,21,16,11,14,13,\adfsplit 38,36,31,49,5,56,53,52,24,44,23,29)_{\adfGI}$,

$(40,21,5,47,49,22,46,30,28,61,17,43,\adfsplit 60,50,25,38,44,8,32,14,9,35,52,56)_{\adfGI}$,

$(33,60,54,61,27,17,10,31,32,19,14,0,\adfsplit 46,18,58,53,43,62,41,7,11,57,34,6)_{\adfGI}$,

\adfLgap
$(30,45,23,53,59,2,18,49,10,3,60,47,\adfsplit 8,7,35,28,6,58,4,44,51,5,0,33)_{\adfGJ}$,

$(8,54,33,55,40,21,57,46,23,37,53,56,\adfsplit 61,48,2,24,25,31,17,19,27,34,52,20)_{\adfGJ}$,

$(30,61,49,9,16,3,33,27,55,24,32,42,\adfsplit 40,14,7,57,25,13,23,62,41,50,38,34)_{\adfGJ}$,

\adfLgap
$(20,52,19,39,41,11,60,53,16,22,6,49,\adfsplit 0,21,18,3,45,55,29,24,50,61,48,30)_{\adfGK}$,

$(15,38,12,60,57,47,25,45,42,9,26,4,\adfsplit 28,59,58,49,30,17,39,56,7,11,43,0)_{\adfGK}$,

$(28,52,11,48,27,31,26,22,34,20,61,16,\adfsplit 21,5,56,29,49,40,9,33,14,32,30,62)_{\adfGK}$

\noindent and

$(3,22,58,36,49,57,32,26,31,25,13,24,\adfsplit 29,50,43,12,60,4,53,41,56,42,46,21)_{\adfGL}$,

$(1,46,61,56,50,17,47,42,16,23,12,59,\adfsplit 8,15,3,49,20,43,10,37,27,2,26,41)_{\adfGL}$,

$(21,14,50,36,53,48,39,24,17,15,33,40,\adfsplit 57,52,62,10,56,18,3,49,43,37,6,44)_{\adfGL}$

\noindent under the action of the mapping $x \mapsto x + 4$ (mod 48) for $x < 48$,
$x \mapsto 48 + (x - 48 + 5 \mathrm{~(mod~15)})$ for $x \ge 48$.

\adfVfy{63, \{\{48,12,4\},\{15,3,5\}\}, -1, \{\{24,\{0,1\}\},\{15,\{2\}\}\}, -1} 

Let the vertex set of $K_{72,72,63}$ be $\{0, 1, \dots, 206\}$ partitioned into
$\{2j + i: j = 0, 1, \dots, 71\}$, $i = 0, 1$, and $\{144, 145, \dots, 206\}$.
The decompositions consist of

$(42,173,182,189,138,48,139,101,31,135,146,50,\adfsplit 180,52,108,123,113,193,119,12,26,199,76,95)_{\adfGa}$,

$(43,174,185,192,139,49,140,102,32,136,147,51,\adfsplit 183,53,109,124,114,196,120,13,27,202,77,96)_{\adfGa}$,

$(44,175,188,195,140,50,141,103,33,137,148,52,\adfsplit 186,54,110,125,115,199,121,14,28,205,78,97)_{\adfGa}$,

$(45,176,191,198,141,51,142,104,34,138,149,53,\adfsplit 189,55,111,126,116,202,122,15,29,181,79,98)_{\adfGa}$,

$(121,148,200,155,34,22,101,87,37,125,180,199,\adfsplit 96,169,63,28,197,93,70,67,53,94,178,158)_{\adfGa}$,

$(122,149,203,156,35,23,102,88,38,126,183,202,\adfsplit 97,170,64,29,200,94,71,68,54,95,179,159)_{\adfGa}$,

$(123,150,206,157,36,24,103,89,39,127,186,205,\adfsplit 98,171,65,30,203,95,72,69,55,96,144,160)_{\adfGa}$,

$(124,151,182,158,37,25,104,90,40,128,189,181,\adfsplit 99,172,66,31,206,96,73,70,56,97,145,161)_{\adfGa}$,

$(30,143,146,133,164,80,96,85,84,92,167,173,\adfsplit 40,35,67,44,176,168,38,91,135,46,87,116)_{\adfGa}$,

$(164,125,77,25,10,30,170,177,149,72,53,103,\adfsplit 61,142,78,146,32,156,91,118,171,113,136,5)_{\adfGa}$,

$(77,169,154,145,0,118,115,73,131,15,174,134,\adfsplit 48,163,8,111,61,39,121,106,74,147,102,167)_{\adfGa}$,

\adfLgap
$(14,154,39,111,64,3,181,46,148,183,134,33,\adfsplit 176,113,56,155,55,197,73,124,94,171,12,17)_{\adfGb}$,

$(15,155,40,112,65,4,184,47,149,186,135,34,\adfsplit 177,114,57,156,56,200,74,125,95,172,13,18)_{\adfGb}$,

$(16,156,41,113,66,5,187,48,150,189,136,35,\adfsplit 178,115,58,157,57,203,75,126,96,173,14,19)_{\adfGb}$,

$(17,157,42,114,67,6,190,49,151,192,137,36,\adfsplit 179,116,59,158,58,206,76,127,97,174,15,20)_{\adfGb}$,

$(157,132,91,30,129,166,60,2,162,164,128,148,\adfsplit 59,14,119,95,32,198,161,202,61,105,12,56)_{\adfGb}$,

$(158,133,92,31,130,167,61,3,163,165,129,149,\adfsplit 60,15,120,96,33,201,162,205,62,106,13,57)_{\adfGb}$,

$(159,134,93,32,131,168,62,4,164,166,130,150,\adfsplit 61,16,121,97,34,204,163,181,63,107,14,58)_{\adfGb}$,

$(160,135,94,33,132,169,63,5,165,167,131,151,\adfsplit 62,17,122,98,35,180,164,184,64,108,15,59)_{\adfGb}$,

$(161,50,66,43,183,196,85,7,32,201,97,120,\adfsplit 174,135,55,179,128,86,14,104,79,200,5,191)_{\adfGb}$,

$(25,156,36,134,9,2,197,55,203,119,199,47,\adfsplit 20,21,110,132,73,6,181,189,32,138,5,97)_{\adfGb}$,

$(186,86,92,15,182,200,33,11,24,194,10,118,\adfsplit 64,73,77,119,49,122,96,191,193,112,111,103)_{\adfGb}$,

\adfLgap
$(41,138,70,164,171,166,33,95,63,25,48,152,\adfsplit 147,169,56,27,7,38,4,23,121,180,18,109)_{\adfGc}$,

$(42,139,71,165,172,167,34,96,64,26,49,153,\adfsplit 148,170,57,28,8,39,5,24,122,183,19,110)_{\adfGc}$,

$(43,140,72,166,173,168,35,97,65,27,50,154,\adfsplit 149,171,58,29,9,40,6,25,123,186,20,111)_{\adfGc}$,

$(44,141,73,167,174,169,36,98,66,28,51,155,\adfsplit 150,172,59,30,10,41,7,26,124,189,21,112)_{\adfGc}$,

$(92,99,125,3,14,176,148,98,184,66,101,56,\adfsplit 174,108,147,155,52,53,193,169,140,71,39,4)_{\adfGc}$,

$(93,100,126,4,15,177,149,99,187,67,102,57,\adfsplit 175,109,148,156,53,54,196,170,141,72,40,5)_{\adfGc}$,

$(94,101,127,5,16,178,150,100,190,68,103,58,\adfsplit 176,110,149,157,54,55,199,171,142,73,41,6)_{\adfGc}$,

$(95,102,128,6,17,179,151,101,193,69,104,59,\adfsplit 177,111,150,158,55,56,202,172,143,74,42,7)_{\adfGc}$,

$(180,15,65,128,181,191,118,110,109,111,92,192,\adfsplit 193,194,102,47,68,28,11,85,45,138,198,182)_{\adfGc}$,

$(13,202,12,30,71,104,191,21,195,113,204,3,\adfsplit 112,8,96,2,59,103,189,197,141,194,142,122)_{\adfGc}$,

$(187,99,98,9,6,197,188,11,183,206,49,75,\adfsplit 138,51,76,30,129,203,201,112,25,191,182,32)_{\adfGc}$,

\adfLgap
$(187,125,23,1,12,76,144,186,150,182,85,94,\adfsplit 65,55,71,90,3,82,153,154,56,171,176,114)_{\adfGd}$,

$(190,126,24,2,13,77,145,189,151,185,86,95,\adfsplit 66,56,72,91,4,83,154,155,57,172,177,115)_{\adfGd}$,

$(193,127,25,3,14,78,146,192,152,188,87,96,\adfsplit 67,57,73,92,5,84,155,156,58,173,178,116)_{\adfGd}$,

$(196,128,26,4,15,79,147,195,153,191,88,97,\adfsplit 68,58,74,93,6,85,156,157,59,174,179,117)_{\adfGd}$,

$(148,3,122,23,58,108,143,71,165,164,133,161,\adfsplit 175,114,45,153,79,61,12,180,41,38,76,113)_{\adfGd}$,

$(149,4,123,24,59,109,0,72,166,165,134,162,\adfsplit 176,115,46,154,80,62,13,183,42,39,77,114)_{\adfGd}$,

$(150,5,124,25,60,110,1,73,167,166,135,163,\adfsplit 177,116,47,155,81,63,14,186,43,40,78,115)_{\adfGd}$,

$(151,6,125,26,61,111,2,74,168,167,136,164,\adfsplit 178,117,48,156,82,64,15,189,44,41,79,116)_{\adfGd}$,

$(77,10,202,36,203,89,45,116,135,63,205,86,\adfsplit 127,94,182,111,190,180,66,31,33,34,100,185)_{\adfGd}$,

$(120,101,139,184,36,201,203,74,69,100,17,121,\adfsplit 93,24,194,4,195,187,138,135,91,137,78,202)_{\adfGd}$,

$(186,55,46,132,22,96,33,184,85,71,191,185,\adfsplit 38,4,82,64,143,65,193,206,135,127,18,80)_{\adfGd}$,

\adfLgap
$(117,161,193,36,55,98,96,53,158,125,41,167,\adfsplit 191,119,166,78,93,88,133,150,60,66,162,21)_{\adfGe}$,

$(118,162,196,37,56,99,97,54,159,126,42,168,\adfsplit 194,120,167,79,94,89,134,151,61,67,163,22)_{\adfGe}$,

$(119,163,199,38,57,100,98,55,160,127,43,169,\adfsplit 197,121,168,80,95,90,135,152,62,68,164,23)_{\adfGe}$,

$(120,164,202,39,58,101,99,56,161,128,44,170,\adfsplit 200,122,169,81,96,91,136,153,63,69,165,24)_{\adfGe}$,

$(39,130,28,10,139,129,158,176,148,179,175,126,\adfsplit 89,92,48,41,200,197,97,198,112,73,116,15)_{\adfGe}$,

$(40,131,29,11,140,130,159,177,149,144,176,127,\adfsplit 90,93,49,42,203,200,98,201,113,74,117,16)_{\adfGe}$,

$(41,132,30,12,141,131,160,178,150,145,177,128,\adfsplit 91,94,50,43,206,203,99,204,114,75,118,17)_{\adfGe}$,

$(42,133,31,13,142,132,161,179,151,146,178,129,\adfsplit 92,95,51,44,182,206,100,180,115,76,119,18)_{\adfGe}$,

$(151,108,42,15,85,205,202,25,136,66,124,89,\adfsplit 50,71,195,137,115,181,34,54,192,199,51,133)_{\adfGe}$,

$(166,123,66,57,110,160,196,117,195,22,96,17,\adfsplit 130,47,190,91,95,140,129,204,70,112,186,1)_{\adfGe}$,

$(47,68,192,60,204,7,18,92,135,190,2,87,\adfsplit 117,100,64,125,96,193,173,189,37,71,94,112)_{\adfGe}$,

\adfLgap
$(39,56,140,158,37,59,148,144,70,8,118,71,\adfsplit 98,31,150,205,107,195,6,142,172,161,5,63)_{\adfGf}$,

$(40,57,141,159,38,60,149,145,71,9,119,72,\adfsplit 99,32,151,181,108,198,7,143,173,162,6,64)_{\adfGf}$,

$(41,58,142,160,39,61,150,146,72,10,120,73,\adfsplit 100,33,152,184,109,201,8,0,174,163,7,65)_{\adfGf}$,

$(42,59,143,161,40,62,151,147,73,11,121,74,\adfsplit 101,34,153,187,110,204,9,1,175,164,8,66)_{\adfGf}$,

$(58,79,150,7,48,72,19,45,183,192,47,90,\adfsplit 191,53,104,196,94,18,171,26,105,190,73,169)_{\adfGf}$,

$(59,80,151,8,49,73,20,46,186,195,48,91,\adfsplit 194,54,105,199,95,19,172,27,106,193,74,170)_{\adfGf}$,

$(60,81,152,9,50,74,21,47,189,198,49,92,\adfsplit 197,55,106,202,96,20,173,28,107,196,75,171)_{\adfGf}$,

$(61,82,153,10,51,75,22,48,192,201,50,93,\adfsplit 200,56,107,205,97,21,174,29,108,199,76,172)_{\adfGf}$,

$(9,4,144,155,53,17,58,70,98,20,182,111,\adfsplit 21,179,115,74,191,92,40,51,3,10,139,194)_{\adfGf}$,

$(163,67,138,38,128,54,143,17,203,25,101,157,\adfsplit 78,127,158,4,139,152,134,133,182,148,59,92)_{\adfGf}$,

$(135,32,38,4,27,197,173,53,150,73,134,12,\adfsplit 16,61,206,35,137,130,52,170,194,84,145,89)_{\adfGf}$,

\adfLgap
$(105,118,70,54,21,31,172,157,161,154,197,130,\adfsplit 75,66,2,23,148,142,109,186,57,72,36,201)_{\adfGg}$,

$(106,119,71,55,22,32,173,158,162,155,200,131,\adfsplit 76,67,3,24,149,143,110,189,58,73,37,204)_{\adfGg}$,

$(107,120,72,56,23,33,174,159,163,156,203,132,\adfsplit 77,68,4,25,150,0,111,192,59,74,38,180)_{\adfGg}$,

$(108,121,73,57,24,34,175,160,164,157,206,133,\adfsplit 78,69,5,26,151,1,112,195,60,75,39,183)_{\adfGg}$,

$(154,37,106,53,100,110,77,101,158,114,146,170,\adfsplit 150,16,49,84,40,95,43,183,148,59,9,26)_{\adfGg}$,

$(155,38,107,54,101,111,78,102,159,115,147,171,\adfsplit 151,17,50,85,41,96,44,186,149,60,10,27)_{\adfGg}$,

$(156,39,108,55,102,112,79,103,160,116,148,172,\adfsplit 152,18,51,86,42,97,45,189,150,61,11,28)_{\adfGg}$,

$(157,40,109,56,103,113,80,104,161,117,149,173,\adfsplit 153,19,52,87,43,98,46,192,151,62,12,29)_{\adfGg}$,

$(18,167,25,205,43,84,184,202,32,132,185,121,\adfsplit 46,187,93,30,45,182,130,85,54,196,200,39)_{\adfGg}$,

$(60,148,113,19,137,143,199,102,184,174,52,86,\adfsplit 43,111,97,182,50,190,181,104,96,37,185,18)_{\adfGg}$,

$(65,84,200,72,182,125,31,22,91,196,193,76,\adfsplit 121,197,140,75,173,199,60,43,35,74,54,119)_{\adfGg}$,

\adfLgap
$(183,18,29,33,101,180,28,78,204,169,49,125,\adfsplit 99,37,108,52,144,70,178,158,140,135,1,112)_{\adfGh}$,

$(186,19,30,34,102,183,29,79,180,170,50,126,\adfsplit 100,38,109,53,145,71,179,159,141,136,2,113)_{\adfGh}$,

$(189,20,31,35,103,186,30,80,183,171,51,127,\adfsplit 101,39,110,54,146,72,144,160,142,137,3,114)_{\adfGh}$,

$(192,21,32,36,104,189,31,81,186,172,52,128,\adfsplit 102,40,111,55,147,73,145,161,143,138,4,115)_{\adfGh}$,

$(139,150,171,124,2,83,13,64,146,172,170,66,\adfsplit 101,11,131,133,194,158,20,98,161,18,24,184)_{\adfGh}$,

$(140,151,172,125,3,84,14,65,147,173,171,67,\adfsplit 102,12,132,134,197,159,21,99,162,19,25,187)_{\adfGh}$,

$(141,152,173,126,4,85,15,66,148,174,172,68,\adfsplit 103,13,133,135,200,160,22,100,163,20,26,190)_{\adfGh}$,

$(142,153,174,127,5,86,16,67,149,175,173,69,\adfsplit 104,14,134,136,203,161,23,101,164,21,27,193)_{\adfGh}$,

$(178,124,119,17,23,203,197,8,190,0,70,77,\adfsplit 15,132,137,5,116,184,199,13,89,185,50,30)_{\adfGh}$,

$(191,61,127,45,18,86,144,145,181,188,93,85,\adfsplit 56,68,126,147,94,65,31,196,184,27,52,70)_{\adfGh}$,

$(125,193,138,80,90,70,11,200,1,35,15,2,\adfsplit 7,56,22,92,203,96,16,135,197,47,79,196)_{\adfGh}$,

\adfLgap
$(35,146,160,52,17,19,122,54,176,187,0,59,\adfsplit 185,99,138,203,113,140,83,171,92,109,153,106)_{\adfGi}$,

$(36,147,161,53,18,20,123,55,177,190,1,60,\adfsplit 188,100,139,206,114,141,84,172,93,110,154,107)_{\adfGi}$,

$(37,148,162,54,19,21,124,56,178,193,2,61,\adfsplit 191,101,140,182,115,142,85,173,94,111,155,108)_{\adfGi}$,

$(38,149,163,55,20,22,125,57,179,196,3,62,\adfsplit 194,102,141,185,116,143,86,174,95,112,156,109)_{\adfGi}$,

$(120,63,146,135,184,187,100,42,84,194,62,121,\adfsplit 53,13,109,158,178,102,145,180,46,108,143,7)_{\adfGi}$,

$(121,64,147,136,187,190,101,43,85,197,63,122,\adfsplit 54,14,110,159,179,103,146,183,47,109,0,8)_{\adfGi}$,

$(122,65,148,137,190,193,102,44,86,200,64,123,\adfsplit 55,15,111,160,144,104,147,186,48,110,1,9)_{\adfGi}$,

$(123,66,149,138,193,196,103,45,87,203,65,124,\adfsplit 56,16,112,161,145,105,148,189,49,111,2,10)_{\adfGi}$,

$(44,159,123,39,112,59,166,180,170,186,14,85,\adfsplit 36,142,26,192,103,63,119,153,164,109,116,2)_{\adfGi}$,

$(130,178,186,153,100,131,55,70,23,33,179,140,\adfsplit 75,102,116,49,89,6,189,43,13,144,168,30)_{\adfGi}$,

$(153,108,3,125,33,201,92,40,46,112,135,13,\adfsplit 204,1,195,82,60,142,25,151,78,109,137,20)_{\adfGi}$,

\adfLgap
$(48,81,139,203,12,152,193,10,95,119,20,27,\adfsplit 164,171,161,78,178,205,18,69,141,52,41,56)_{\adfGj}$,

$(49,82,140,206,13,153,196,11,96,120,21,28,\adfsplit 165,172,162,79,179,181,19,70,142,53,42,57)_{\adfGj}$,

$(50,83,141,182,14,154,199,12,97,121,22,29,\adfsplit 166,173,163,80,144,184,20,71,143,54,43,58)_{\adfGj}$,

$(51,84,142,185,15,155,202,13,98,122,23,30,\adfsplit 167,174,164,81,145,187,21,72,0,55,44,59)_{\adfGj}$,

$(168,125,33,21,144,132,76,149,72,174,83,173,\adfsplit 26,81,44,53,180,28,146,193,18,82,79,105)_{\adfGj}$,

$(169,126,34,22,145,133,77,150,73,175,84,174,\adfsplit 27,82,45,54,183,29,147,196,19,83,80,106)_{\adfGj}$,

$(170,127,35,23,146,134,78,151,74,176,85,175,\adfsplit 28,83,46,55,186,30,148,199,20,84,81,107)_{\adfGj}$,

$(171,128,36,24,147,135,79,152,75,177,86,176,\adfsplit 29,84,47,56,189,31,149,202,21,85,82,108)_{\adfGj}$,

$(90,61,43,125,56,151,185,74,188,192,27,49,\adfsplit 55,111,22,126,189,204,180,110,140,182,120,75)_{\adfGj}$,

$(91,56,20,60,169,103,117,144,115,150,183,136,\adfsplit 90,8,37,98,25,188,109,119,201,186,114,82)_{\adfGj}$,

$(89,124,42,58,189,93,107,191,13,131,14,88,\adfsplit 113,128,197,104,59,33,194,76,123,49,86,203)_{\adfGj}$,

\adfLgap
$(9,102,44,114,163,13,27,182,57,127,159,152,\adfsplit 119,175,28,191,36,117,151,146,132,54,49,7)_{\adfGk}$,

$(10,103,45,115,164,14,28,185,58,128,160,153,\adfsplit 120,176,29,194,37,118,152,147,133,55,50,8)_{\adfGk}$,

$(11,104,46,116,165,15,29,188,59,129,161,154,\adfsplit 121,177,30,197,38,119,153,148,134,56,51,9)_{\adfGk}$,

$(12,105,47,117,166,16,30,191,60,130,162,155,\adfsplit 122,178,31,200,39,120,154,149,135,57,52,10)_{\adfGk}$,

$(1,163,185,151,41,79,42,112,32,70,132,158,\adfsplit 153,169,47,23,117,96,9,120,76,93,205,184)_{\adfGk}$,

$(2,164,188,152,42,80,43,113,33,71,133,159,\adfsplit 154,170,48,24,118,97,10,121,77,94,181,187)_{\adfGk}$,

$(3,165,191,153,43,81,44,114,34,72,134,160,\adfsplit 155,171,49,25,119,98,11,122,78,95,184,190)_{\adfGk}$,

$(4,166,194,154,44,82,45,115,35,73,135,161,\adfsplit 156,172,50,26,120,99,12,123,79,96,187,193)_{\adfGk}$,

$(159,122,99,65,133,81,157,74,108,180,193,120,\adfsplit 183,204,83,123,112,201,85,102,6,53,189,71)_{\adfGk}$,

$(126,187,123,181,11,137,86,62,124,97,30,129,\adfsplit 158,127,204,36,131,195,64,85,196,116,105,82)_{\adfGk}$,

$(140,3,33,61,186,201,128,58,64,68,136,79,\adfsplit 204,160,189,43,75,0,134,113,123,30,180,82)_{\adfGk}$,

\adfLgap
$(34,13,111,153,164,0,8,189,98,27,51,205,\adfsplit 35,190,16,46,147,75,188,183,63,103,132,128)_{\adfGl}$,

$(35,14,112,154,165,1,9,192,99,28,52,181,\adfsplit 36,193,17,47,148,76,191,186,64,104,133,129)_{\adfGl}$,

$(36,15,113,155,166,2,10,195,100,29,53,184,\adfsplit 37,196,18,48,149,77,194,189,65,105,134,130)_{\adfGl}$,

$(37,16,114,156,167,3,11,198,101,30,54,187,\adfsplit 38,199,19,49,150,78,197,192,66,106,135,131)_{\adfGl}$,

$(65,56,172,178,57,189,72,68,26,128,133,115,\adfsplit 161,181,188,173,76,32,11,4,16,111,145,49)_{\adfGl}$,

$(66,57,173,179,58,192,73,69,27,129,134,116,\adfsplit 162,184,191,174,77,33,12,5,17,112,146,50)_{\adfGl}$,

$(67,58,174,144,59,195,74,70,28,130,135,117,\adfsplit 163,187,194,175,78,34,13,6,18,113,147,51)_{\adfGl}$,

$(68,59,175,145,60,198,75,71,29,131,136,118,\adfsplit 164,190,197,176,79,35,14,7,19,114,148,52)_{\adfGl}$,

$(185,47,63,72,146,70,110,158,79,99,139,149,\adfsplit 112,6,50,4,179,75,168,109,101,145,69,28)_{\adfGl}$,

$(71,16,144,153,19,39,92,108,58,64,147,155,\adfsplit 3,51,127,37,110,8,185,30,54,176,137,77)_{\adfGl}$,

$(15,175,78,162,14,44,53,154,17,37,133,145,\adfsplit 26,56,126,118,206,57,203,144,9,49,0,32)_{\adfGl}$,

\adfLgap
$(0,196,150,93,138,20,119,87,161,180,169,168,\adfsplit 149,23,121,71,82,88,65,130,70,152,56,143)_{\adfGm}$,

$(1,199,151,94,139,21,120,88,162,183,170,169,\adfsplit 150,24,122,72,83,89,66,131,71,153,57,0)_{\adfGm}$,

$(2,202,152,95,140,22,121,89,163,186,171,170,\adfsplit 151,25,123,73,84,90,67,132,72,154,58,1)_{\adfGm}$,

$(3,205,153,96,141,23,122,90,164,189,172,171,\adfsplit 152,26,124,74,85,91,68,133,73,155,59,2)_{\adfGm}$,

$(114,69,91,158,52,177,167,112,11,84,95,97,\adfsplit 192,180,188,41,70,98,54,127,115,46,146,152)_{\adfGm}$,

$(115,70,92,159,53,178,168,113,12,85,96,98,\adfsplit 195,183,191,42,71,99,55,128,116,47,147,153)_{\adfGm}$,

$(116,71,93,160,54,179,169,114,13,86,97,99,\adfsplit 198,186,194,43,72,100,56,129,117,48,148,154)_{\adfGm}$,

$(117,72,94,161,55,144,170,115,14,87,98,100,\adfsplit 201,189,197,44,73,101,57,130,118,49,149,155)_{\adfGm}$,

$(161,73,18,0,126,24,77,131,188,7,103,142,\adfsplit 148,22,194,6,193,12,60,41,61,197,59,191)_{\adfGm}$,

$(150,134,7,97,133,75,199,114,28,182,184,124,\adfsplit 77,81,17,94,37,96,190,181,74,50,15,131)_{\adfGm}$,

$(142,202,47,135,16,64,191,171,197,60,5,30,\adfsplit 82,75,99,91,187,40,117,13,124,185,4,196)_{\adfGm}$,

\adfLgap
$(22,165,11,151,16,142,201,18,41,39,170,73,\adfsplit 191,74,25,10,4,160,150,135,148,81,92,78)_{\adfGn}$,

$(23,166,12,152,17,143,204,19,42,40,171,74,\adfsplit 194,75,26,11,5,161,151,136,149,82,93,79)_{\adfGn}$,

$(24,167,13,153,18,0,180,20,43,41,172,75,\adfsplit 197,76,27,12,6,162,152,137,150,83,94,80)_{\adfGn}$,

$(25,168,14,154,19,1,183,21,44,42,173,76,\adfsplit 200,77,28,13,7,163,153,138,151,84,95,81)_{\adfGn}$,

$(165,8,9,74,97,176,102,151,188,194,79,3,\adfsplit 140,98,187,20,41,150,83,146,11,64,124,107)_{\adfGn}$,

$(166,9,10,75,98,177,103,152,191,197,80,4,\adfsplit 141,99,190,21,42,151,84,147,12,65,125,108)_{\adfGn}$,

$(167,10,11,76,99,178,104,153,194,200,81,5,\adfsplit 142,100,193,22,43,152,85,148,13,66,126,109)_{\adfGn}$,

$(168,11,12,77,100,179,105,154,197,203,82,6,\adfsplit 143,101,196,23,44,153,86,149,14,67,127,110)_{\adfGn}$,

$(161,116,11,128,129,19,192,195,193,115,56,172,\adfsplit 30,130,180,103,91,71,44,187,74,204,181,132)_{\adfGn}$,

$(62,179,101,186,137,2,40,134,30,125,49,119,\adfsplit 95,8,166,4,16,35,133,190,33,192,106,140)_{\adfGn}$,

$(69,54,36,82,183,133,101,21,53,204,196,199,\adfsplit 92,7,66,127,19,138,5,112,90,190,205,27)_{\adfGn}$,

\adfLgap
$(98,180,133,198,127,15,4,154,108,1,40,205,\adfsplit 91,185,134,19,60,124,17,193,32,115,160,12)_{\adfGo}$,

$(99,183,134,201,128,16,5,155,109,2,41,181,\adfsplit 92,188,135,20,61,125,18,196,33,116,161,13)_{\adfGo}$,

$(100,186,135,204,129,17,6,156,110,3,42,184,\adfsplit 93,191,136,21,62,126,19,199,34,117,162,14)_{\adfGo}$,

$(101,189,136,180,130,18,7,157,111,4,43,187,\adfsplit 94,194,137,22,63,127,20,202,35,118,163,15)_{\adfGo}$,

$(9,108,163,186,153,174,8,22,110,61,15,11,\adfsplit 205,149,10,90,63,39,194,188,41,122,116,49)_{\adfGo}$,

$(10,109,164,189,154,175,9,23,111,62,16,12,\adfsplit 181,150,11,91,64,40,197,191,42,123,117,50)_{\adfGo}$,

$(11,110,165,192,155,176,10,24,112,63,17,13,\adfsplit 184,151,12,92,65,41,200,194,43,124,118,51)_{\adfGo}$,

$(12,111,166,195,156,177,11,25,113,64,18,14,\adfsplit 187,152,13,93,66,42,203,197,44,125,119,52)_{\adfGo}$,

$(168,42,80,26,166,172,71,127,17,159,98,15,\adfsplit 90,147,58,151,129,174,162,120,5,157,62,73)_{\adfGo}$,

$(163,67,93,104,124,171,158,153,168,61,130,23,\adfsplit 59,101,14,179,169,96,141,135,126,79,157,118)_{\adfGo}$,

$(142,105,41,150,68,48,153,24,29,11,149,158,\adfsplit 156,0,28,62,115,16,5,39,119,144,172,140)_{\adfGo}$,

\adfLgap
$(96,77,61,175,22,136,189,204,143,6,67,156,\adfsplit 37,190,69,78,4,181,112,135,11,171,9,196)_{\adfGp}$,

$(97,78,62,176,23,137,192,180,0,7,68,157,\adfsplit 38,193,70,79,5,184,113,136,12,172,10,199)_{\adfGp}$,

$(98,79,63,177,24,138,195,183,1,8,69,158,\adfsplit 39,196,71,80,6,187,114,137,13,173,11,202)_{\adfGp}$,

$(99,80,64,178,25,139,198,186,2,9,70,159,\adfsplit 40,199,72,81,7,190,115,138,14,174,12,205)_{\adfGp}$,

$(76,176,182,3,41,73,50,101,146,195,161,94,\adfsplit 149,35,119,96,63,36,168,125,154,46,49,20)_{\adfGp}$,

$(77,177,185,4,42,74,51,102,147,198,162,95,\adfsplit 150,36,120,97,64,37,169,126,155,47,50,21)_{\adfGp}$,

$(78,178,188,5,43,75,52,103,148,201,163,96,\adfsplit 151,37,121,98,65,38,170,127,156,48,51,22)_{\adfGp}$,

$(79,179,191,6,44,76,53,104,149,204,164,97,\adfsplit 152,38,122,99,66,39,171,128,157,49,52,23)_{\adfGp}$,

$(168,66,23,62,152,133,22,176,99,203,138,84,\adfsplit 119,2,79,101,51,104,12,191,161,102,11,103)_{\adfGp}$,

$(108,109,71,145,136,32,188,66,39,59,197,65,\adfsplit 159,106,14,139,154,29,13,44,24,81,151,166)_{\adfGp}$,

$(127,78,100,157,25,151,203,9,40,32,106,113,\adfsplit 42,191,170,91,86,57,132,136,200,37,182,85)_{\adfGp}$,

\adfLgap
$(37,146,60,188,45,25,155,121,110,79,54,186,\adfsplit 96,187,136,81,98,27,149,12,7,198,177,105)_{\adfGq}$,

$(38,147,61,191,46,26,156,122,111,80,55,189,\adfsplit 97,190,137,82,99,28,150,13,8,201,178,106)_{\adfGq}$,

$(39,148,62,194,47,27,157,123,112,81,56,192,\adfsplit 98,193,138,83,100,29,151,14,9,204,179,107)_{\adfGq}$,

$(40,149,63,197,48,28,158,124,113,82,57,195,\adfsplit 99,196,139,84,101,30,152,15,10,180,144,108)_{\adfGq}$,

$(203,20,9,5,51,181,48,30,86,188,121,78,\adfsplit 49,174,125,196,33,160,180,74,127,66,128,79)_{\adfGq}$,

$(206,21,10,6,52,184,49,31,87,191,122,79,\adfsplit 50,175,126,199,34,161,183,75,128,67,129,80)_{\adfGq}$,

$(182,22,11,7,53,187,50,32,88,194,123,80,\adfsplit 51,176,127,202,35,162,186,76,129,68,130,81)_{\adfGq}$,

$(185,23,12,8,54,190,51,33,89,197,124,81,\adfsplit 52,177,128,205,36,163,189,77,130,69,131,82)_{\adfGq}$,

$(47,144,10,154,73,84,176,27,82,113,165,93,\adfsplit 145,11,174,42,50,54,112,167,160,149,43,41)_{\adfGq}$,

$(77,88,18,161,146,3,125,157,100,84,136,169,\adfsplit 103,156,41,171,73,158,166,110,99,66,56,7)_{\adfGq}$,

$(46,117,179,174,14,147,25,44,141,84,51,8,\adfsplit 7,104,171,42,111,50,24,176,95,144,167,17)_{\adfGq}$,

\adfLgap
$(149,131,125,1,14,116,188,164,145,46,97,186,\adfsplit 58,20,57,179,135,129,194,187,95,24,126,107)_{\adfGr}$,

$(150,132,126,2,15,117,191,165,146,47,98,189,\adfsplit 59,21,58,144,136,130,197,190,96,25,127,108)_{\adfGr}$,

$(151,133,127,3,16,118,194,166,147,48,99,192,\adfsplit 60,22,59,145,137,131,200,193,97,26,128,109)_{\adfGr}$,

$(152,134,128,4,17,119,197,167,148,49,100,195,\adfsplit 61,23,60,146,138,132,203,196,98,27,129,110)_{\adfGr}$,

$(85,114,188,120,55,135,108,88,172,192,173,140,\adfsplit 39,95,51,8,195,202,118,14,199,166,0,141)_{\adfGr}$,

$(86,115,191,121,56,136,109,89,173,195,174,141,\adfsplit 40,96,52,9,198,205,119,15,202,167,1,142)_{\adfGr}$,

$(87,116,194,122,57,137,110,90,174,198,175,142,\adfsplit 41,97,53,10,201,181,120,16,205,168,2,143)_{\adfGr}$,

$(88,117,197,123,58,138,111,91,175,201,176,143,\adfsplit 42,98,54,11,204,184,121,17,181,169,3,0)_{\adfGr}$,

$(70,31,174,168,138,72,87,29,113,30,163,73,\adfsplit 16,177,157,107,172,156,135,11,34,56,121,161)_{\adfGr}$,

$(124,152,57,5,94,38,159,64,145,52,173,1,\adfsplit 71,125,151,35,102,171,58,138,158,136,103,97)_{\adfGr}$,

$(176,84,20,143,123,121,150,61,56,12,54,179,\adfsplit 62,5,59,95,70,55,146,102,48,96,162,9)_{\adfGr}$,

\adfLgap
$(128,97,75,19,162,108,14,30,16,153,107,83,\adfsplit 172,188,134,122,64,189,53,113,93,27,201,106)_{\adfGs}$,

$(129,98,76,20,163,109,15,31,17,154,108,84,\adfsplit 173,191,135,123,65,192,54,114,94,28,204,107)_{\adfGs}$,

$(130,99,77,21,164,110,16,32,18,155,109,85,\adfsplit 174,194,136,124,66,195,55,115,95,29,180,108)_{\adfGs}$,

$(131,100,78,22,165,111,17,33,19,156,110,86,\adfsplit 175,197,137,125,67,198,56,116,96,30,183,109)_{\adfGs}$,

$(28,77,1,152,6,128,165,102,45,36,205,199,\adfsplit 193,134,172,37,49,111,139,122,106,32,146,201)_{\adfGs}$,

$(29,78,2,153,7,129,166,103,46,37,181,202,\adfsplit 196,135,173,38,50,112,140,123,107,33,147,204)_{\adfGs}$,

$(30,79,3,154,8,130,167,104,47,38,184,205,\adfsplit 199,136,174,39,51,113,141,124,108,34,148,180)_{\adfGs}$,

$(31,80,4,155,9,131,168,105,48,39,187,181,\adfsplit 202,137,175,40,52,114,142,125,109,35,149,183)_{\adfGs}$,

$(79,191,165,10,85,19,52,143,145,157,176,188,\adfsplit 16,63,0,86,20,82,174,170,203,185,89,65)_{\adfGs}$,

$(112,157,165,172,73,103,56,21,136,23,175,80,\adfsplit 194,158,176,14,59,143,49,95,128,126,148,155)_{\adfGs}$,

$(162,72,1,9,206,93,24,122,42,171,185,151,\adfsplit 11,105,22,30,29,10,146,188,176,167,50,31)_{\adfGs}$,

\adfLgap
$(165,126,34,102,67,99,170,113,155,85,132,26,\adfsplit 104,193,179,47,153,172,43,161,107,112,63,100)_{\adfGt}$,

$(166,127,35,103,68,100,171,114,156,86,133,27,\adfsplit 105,196,144,48,154,173,44,162,108,113,64,101)_{\adfGt}$,

$(167,128,36,104,69,101,172,115,157,87,134,28,\adfsplit 106,199,145,49,155,174,45,163,109,114,65,102)_{\adfGt}$,

$(168,129,37,105,70,102,173,116,158,88,135,29,\adfsplit 107,202,146,50,156,175,46,164,110,115,66,103)_{\adfGt}$,

$(54,137,170,176,175,86,68,14,34,25,162,81,\adfsplit 139,105,36,153,46,195,192,97,69,135,20,94)_{\adfGt}$,

$(55,138,171,177,176,87,69,15,35,26,163,82,\adfsplit 140,106,37,154,47,198,195,98,70,136,21,95)_{\adfGt}$,

$(56,139,172,178,177,88,70,16,36,27,164,83,\adfsplit 141,107,38,155,48,201,198,99,71,137,22,96)_{\adfGt}$,

$(57,140,173,179,178,89,71,17,37,28,165,84,\adfsplit 142,108,39,156,49,204,201,100,72,138,23,97)_{\adfGt}$,

$(95,136,197,201,39,200,40,14,38,55,79,13,\adfsplit 182,185,190,126,183,3,129,45,134,32,194,195)_{\adfGt}$,

$(192,16,107,92,205,197,62,185,61,193,6,26,\adfsplit 98,87,137,188,63,184,124,73,32,203,28,182)_{\adfGt}$,

$(80,190,65,67,11,99,187,203,114,82,0,36,\adfsplit 83,1,81,191,49,199,193,38,12,90,137,53)_{\adfGt}$,

\adfLgap
$(173,43,28,92,94,196,89,21,160,169,122,100,\adfsplit 112,153,59,45,101,154,67,130,72,198,52,155)_{\adfGu}$,

$(174,44,29,93,95,199,90,22,161,170,123,101,\adfsplit 113,154,60,46,102,155,68,131,73,201,53,156)_{\adfGu}$,

$(175,45,30,94,96,202,91,23,162,171,124,102,\adfsplit 114,155,61,47,103,156,69,132,74,204,54,157)_{\adfGu}$,

$(176,46,31,95,97,205,92,24,163,172,125,103,\adfsplit 115,156,62,48,104,157,70,133,75,180,55,158)_{\adfGu}$,

$(45,156,112,144,41,50,168,165,140,8,125,155,\adfsplit 21,116,151,134,126,80,127,92,121,115,203,184)_{\adfGu}$,

$(46,157,113,145,42,51,169,166,141,9,126,156,\adfsplit 22,117,152,135,127,81,128,93,122,116,206,187)_{\adfGu}$,

$(47,158,114,146,43,52,170,167,142,10,127,157,\adfsplit 23,118,153,136,128,82,129,94,123,117,182,190)_{\adfGu}$,

$(48,159,115,147,44,53,171,168,143,11,128,158,\adfsplit 24,119,154,137,129,83,130,95,124,118,185,193)_{\adfGu}$,

$(0,31,39,35,64,70,199,22,193,66,190,182,\adfsplit 137,133,9,103,28,194,20,136,138,200,129,25)_{\adfGu}$,

$(126,187,17,15,0,63,191,104,201,94,60,6,\adfsplit 39,183,97,49,192,82,111,136,186,197,90,35)_{\adfGu}$,

$(5,128,195,185,111,191,4,20,114,119,11,135,\adfsplit 49,186,32,194,98,180,42,97,41,118,81,204)_{\adfGu}$,

\adfLgap
$(184,125,61,31,179,190,52,156,36,204,112,124,\adfsplit 73,143,19,131,205,149,0,8,118,164,35,75)_{\adfGv}$,

$(187,126,62,32,144,193,53,157,37,180,113,125,\adfsplit 74,0,20,132,181,150,1,9,119,165,36,76)_{\adfGv}$,

$(190,127,63,33,145,196,54,158,38,183,114,126,\adfsplit 75,1,21,133,184,151,2,10,120,166,37,77)_{\adfGv}$,

$(193,128,64,34,146,199,55,159,39,186,115,127,\adfsplit 76,2,22,134,187,152,3,11,121,167,38,78)_{\adfGv}$,

$(84,19,151,144,92,154,117,4,65,28,31,10,\adfsplit 157,149,73,155,39,95,12,90,130,141,54,158)_{\adfGv}$,

$(85,20,152,145,93,155,118,5,66,29,32,11,\adfsplit 158,150,74,156,40,96,13,91,131,142,55,159)_{\adfGv}$,

$(86,21,153,146,94,156,119,6,67,30,33,12,\adfsplit 159,151,75,157,41,97,14,92,132,143,56,160)_{\adfGv}$,

$(87,22,154,147,95,157,120,7,68,31,34,13,\adfsplit 160,152,76,158,42,98,15,93,133,0,57,161)_{\adfGv}$,

$(26,19,151,29,12,144,121,57,188,36,122,6,\adfsplit 142,204,195,200,91,191,69,60,126,13,124,182)_{\adfGv}$,

$(166,136,113,72,189,153,106,192,69,185,95,123,\adfsplit 22,74,100,62,186,188,25,77,59,195,32,191)_{\adfGv}$,

$(76,125,135,27,104,90,191,38,182,136,113,139,\adfsplit 17,194,192,2,197,186,99,23,123,185,35,120)_{\adfGv}$,

\adfLgap
$(137,171,191,56,50,92,33,62,167,160,37,61,\adfsplit 91,201,40,22,42,148,84,174,168,29,172,143)_{\adfGw}$,

$(138,172,194,57,51,93,34,63,168,161,38,62,\adfsplit 92,204,41,23,43,149,85,175,169,30,173,0)_{\adfGw}$,

$(139,173,197,58,52,94,35,64,169,162,39,63,\adfsplit 93,180,42,24,44,150,86,176,170,31,174,1)_{\adfGw}$,

$(140,174,200,59,53,95,36,65,170,163,40,64,\adfsplit 94,183,43,25,45,151,87,177,171,32,175,2)_{\adfGw}$,

$(142,63,172,169,66,179,103,130,77,32,11,54,\adfsplit 196,203,9,165,3,52,113,50,78,7,56,197)_{\adfGw}$,

$(143,64,173,170,67,144,104,131,78,33,12,55,\adfsplit 199,206,10,166,4,53,114,51,79,8,57,200)_{\adfGw}$,

$(0,65,174,171,68,145,105,132,79,34,13,56,\adfsplit 202,182,11,167,5,54,115,52,80,9,58,203)_{\adfGw}$,

$(1,66,175,172,69,146,106,133,80,35,14,57,\adfsplit 205,185,12,168,6,55,116,53,81,10,59,206)_{\adfGw}$,

$(154,118,71,95,198,127,80,12,48,187,42,122,\adfsplit 47,204,45,173,87,202,23,78,126,65,201,125)_{\adfGw}$,

$(72,180,155,67,88,81,11,96,190,186,4,108,\adfsplit 71,109,49,9,61,138,192,40,34,199,22,7)_{\adfGw}$,

$(168,49,60,1,58,24,184,190,186,6,11,25,\adfsplit 117,130,187,116,110,82,83,48,143,81,196,192)_{\adfGw}$,

\adfLgap
$(59,136,161,110,173,97,78,90,31,159,205,148,\adfsplit 203,143,73,2,94,121,130,166,181,152,87,34)_{\adfGx}$,

$(60,137,162,111,174,98,79,91,32,160,181,149,\adfsplit 206,0,74,3,95,122,131,167,184,153,88,35)_{\adfGx}$,

$(61,138,163,112,175,99,80,92,33,161,184,150,\adfsplit 182,1,75,4,96,123,132,168,187,154,89,36)_{\adfGx}$,

$(62,139,164,113,176,100,81,93,34,162,187,151,\adfsplit 185,2,76,5,97,124,133,169,190,155,90,37)_{\adfGx}$,

$(69,68,118,52,83,137,144,156,163,31,24,154,\adfsplit 87,113,167,152,27,122,186,136,119,64,124,37)_{\adfGx}$,

$(70,69,119,53,84,138,145,157,164,32,25,155,\adfsplit 88,114,168,153,28,123,189,137,120,65,125,38)_{\adfGx}$,

$(71,70,120,54,85,139,146,158,165,33,26,156,\adfsplit 89,115,169,154,29,124,192,138,121,66,126,39)_{\adfGx}$,

$(72,71,121,55,86,140,147,159,166,34,27,157,\adfsplit 90,116,170,155,30,125,195,139,122,67,127,40)_{\adfGx}$,

$(80,203,99,57,89,77,182,195,118,138,140,201,\adfsplit 45,5,27,204,128,13,116,186,38,91,191,120)_{\adfGx}$,

$(140,7,25,57,90,16,187,36,182,94,123,119,\adfsplit 107,193,201,184,142,186,92,71,54,21,188,30)_{\adfGx}$,

$(81,201,182,114,110,140,3,95,139,198,188,21,\adfsplit 14,72,8,25,2,62,11,199,71,44,185,85)_{\adfGx}$,

\adfLgap
$(106,196,169,79,108,29,97,135,175,0,156,30,\adfsplit 42,197,15,10,85,178,51,88,58,177,132,27)_{\adfGy}$,

$(107,199,170,80,109,30,98,136,176,1,157,31,\adfsplit 43,200,16,11,86,179,52,89,59,178,133,28)_{\adfGy}$,

$(108,202,171,81,110,31,99,137,177,2,158,32,\adfsplit 44,203,17,12,87,144,53,90,60,179,134,29)_{\adfGy}$,

$(109,205,172,82,111,32,100,138,178,3,159,33,\adfsplit 45,206,18,13,88,145,54,91,61,144,135,30)_{\adfGy}$,

$(140,200,152,163,15,57,120,5,134,63,126,144,\adfsplit 155,194,179,97,109,90,6,38,49,42,55,192)_{\adfGy}$,

$(141,203,153,164,16,58,121,6,135,64,127,145,\adfsplit 156,197,144,98,110,91,7,39,50,43,56,195)_{\adfGy}$,

$(142,206,154,165,17,59,122,7,136,65,128,146,\adfsplit 157,200,145,99,111,92,8,40,51,44,57,198)_{\adfGy}$,

$(143,182,155,166,18,60,123,8,137,66,129,147,\adfsplit 158,203,146,100,112,93,9,41,52,45,58,201)_{\adfGy}$,

$(193,25,43,128,159,48,38,152,1,189,146,170,\adfsplit 122,126,94,109,104,49,187,19,35,156,70,201)_{\adfGy}$,

$(9,195,4,120,33,134,196,49,157,139,181,202,\adfsplit 108,193,32,27,122,77,205,117,13,183,56,58)_{\adfGy}$,

$(157,107,95,19,124,130,202,140,192,142,167,186,\adfsplit 33,189,43,29,31,24,66,36,204,14,47,3)_{\adfGy}$,

\adfLgap
$(190,101,30,5,194,186,91,123,58,36,31,140,\adfsplit 153,110,185,179,3,11,26,137,113,126,161,145)_{\adfGz}$,

$(193,102,31,6,197,189,92,124,59,37,32,141,\adfsplit 154,111,188,144,4,12,27,138,114,127,162,146)_{\adfGz}$,

$(196,103,32,7,200,192,93,125,60,38,33,142,\adfsplit 155,112,191,145,5,13,28,139,115,128,163,147)_{\adfGz}$,

$(199,104,33,8,203,195,94,126,61,39,34,143,\adfsplit 156,113,194,146,6,14,29,140,116,129,164,148)_{\adfGz}$,

$(136,193,149,178,99,98,32,50,60,95,186,79,\adfsplit 3,5,96,106,200,54,20,192,196,147,53,57)_{\adfGz}$,

$(137,196,150,179,100,99,33,51,61,96,189,80,\adfsplit 4,6,97,107,203,55,21,195,199,148,54,58)_{\adfGz}$,

$(138,199,151,144,101,100,34,52,62,97,192,81,\adfsplit 5,7,98,108,206,56,22,198,202,149,55,59)_{\adfGz}$,

$(139,202,152,145,102,101,35,53,63,98,195,82,\adfsplit 6,8,99,109,182,57,23,201,205,150,56,60)_{\adfGz}$,

$(104,174,31,83,84,6,176,10,144,88,11,91,\adfsplit 79,147,155,60,42,16,175,119,8,65,152,153)_{\adfGz}$,

$(34,99,173,93,88,40,47,76,167,82,109,61,\adfsplit 113,85,174,120,169,130,18,160,136,148,33,57)_{\adfGz}$,

$(155,121,31,118,177,48,138,154,81,97,39,127,\adfsplit 46,9,149,178,92,166,55,14,98,42,33,87)_{\adfGz}$,

\adfLgap
$(54,91,166,205,182,187,124,6,98,59,140,52,\adfsplit 73,139,178,169,174,167,128,108,58,25,35,148)_{\adfGA}$,

$(55,92,167,181,185,190,125,7,99,60,141,53,\adfsplit 74,140,179,170,175,168,129,109,59,26,36,149)_{\adfGA}$,

$(56,93,168,184,188,193,126,8,100,61,142,54,\adfsplit 75,141,144,171,176,169,130,110,60,27,37,150)_{\adfGA}$,

$(57,94,169,187,191,196,127,9,101,62,143,55,\adfsplit 76,142,145,172,177,170,131,111,61,28,38,151)_{\adfGA}$,

$(35,6,134,26,158,101,81,193,31,166,32,142,\adfsplit 92,96,13,45,163,143,65,192,149,170,2,112)_{\adfGA}$,

$(36,7,135,27,159,102,82,196,32,167,33,143,\adfsplit 93,97,14,46,164,0,66,195,150,171,3,113)_{\adfGA}$,

$(37,8,136,28,160,103,83,199,33,168,34,0,\adfsplit 94,98,15,47,165,1,67,198,151,172,4,114)_{\adfGA}$,

$(38,9,137,29,161,104,84,202,34,169,35,1,\adfsplit 95,99,16,48,166,2,68,201,152,173,5,115)_{\adfGA}$,

$(153,123,22,48,180,146,109,85,111,186,44,5,\adfsplit 78,151,125,200,56,26,63,118,14,107,197,204)_{\adfGA}$,

$(80,59,204,101,0,44,39,79,197,182,95,129,\adfsplit 136,122,38,52,152,24,119,37,3,17,200,126)_{\adfGA}$,

$(85,28,46,34,45,201,194,186,198,139,137,70,\adfsplit 64,44,122,188,11,113,203,79,183,80,130,37)_{\adfGA}$,

\adfLgap
$(140,137,13,115,4,151,147,56,183,52,83,0,\adfsplit 143,177,171,155,59,153,161,32,44,113,169,23)_{\adfGB}$,

$(141,138,14,116,5,152,148,57,186,53,84,1,\adfsplit 0,178,172,156,60,154,162,33,45,114,170,24)_{\adfGB}$,

$(142,139,15,117,6,153,149,58,189,54,85,2,\adfsplit 1,179,173,157,61,155,163,34,46,115,171,25)_{\adfGB}$,

$(143,140,16,118,7,154,150,59,192,55,86,3,\adfsplit 2,144,174,158,62,156,164,35,47,116,172,26)_{\adfGB}$,

$(51,120,24,64,144,115,69,31,119,170,2,148,\adfsplit 202,92,93,173,45,100,98,57,20,109,182,190)_{\adfGB}$,

$(52,121,25,65,145,116,70,32,120,171,3,149,\adfsplit 205,93,94,174,46,101,99,58,21,110,185,193)_{\adfGB}$,

$(53,122,26,66,146,117,71,33,121,172,4,150,\adfsplit 181,94,95,175,47,102,100,59,22,111,188,196)_{\adfGB}$,

$(54,123,27,67,147,118,72,34,122,173,5,151,\adfsplit 184,95,96,176,48,103,101,60,23,112,191,199)_{\adfGB}$,

$(81,42,180,118,75,91,112,110,197,184,130,24,\adfsplit 129,203,37,23,47,182,183,56,36,102,181,113)_{\adfGB}$,

$(116,197,21,65,35,48,128,44,205,181,19,29,\adfsplit 83,73,38,20,58,106,64,186,201,198,7,129)_{\adfGB}$,

$(36,69,21,51,70,118,186,200,182,74,53,3,\adfsplit 32,95,197,86,18,183,127,198,9,87,185,48)_{\adfGB}$,

\adfLgap
$(157,26,40,103,61,173,167,135,22,124,39,51,\adfsplit 14,80,15,164,174,204,19,81,92,159,143,90)_{\adfGC}$,

$(158,27,41,104,62,174,168,136,23,125,40,52,\adfsplit 15,81,16,165,175,180,20,82,93,160,0,91)_{\adfGC}$,

$(159,28,42,105,63,175,169,137,24,126,41,53,\adfsplit 16,82,17,166,176,183,21,83,94,161,1,92)_{\adfGC}$,

$(160,29,43,106,64,176,170,138,25,127,42,54,\adfsplit 17,83,18,167,177,186,22,84,95,162,2,93)_{\adfGC}$,

$(95,173,76,177,136,24,160,150,5,3,157,154,\adfsplit 94,34,198,50,89,35,185,7,100,123,190,108)_{\adfGC}$,

$(96,174,77,178,137,25,161,151,6,4,158,155,\adfsplit 95,35,201,51,90,36,188,8,101,124,193,109)_{\adfGC}$,

$(97,175,78,179,138,26,162,152,7,5,159,156,\adfsplit 96,36,204,52,91,37,191,9,102,125,196,110)_{\adfGC}$,

$(98,176,79,144,139,27,163,153,8,6,160,157,\adfsplit 97,37,180,53,92,38,194,10,103,126,199,111)_{\adfGC}$,

$(135,138,182,102,65,25,98,75,193,143,62,140,\adfsplit 205,6,43,200,107,41,180,91,128,60,104,187)_{\adfGC}$,

$(78,95,7,81,58,124,181,200,182,8,121,21,\adfsplit 34,59,18,201,186,106,193,132,75,93,107,194)_{\adfGC}$,

$(139,108,204,36,1,33,18,4,206,205,74,190,\adfsplit 103,194,9,45,28,61,106,89,120,181,200,72)_{\adfGC}$,

\adfLgap
$(119,178,54,166,44,62,135,61,43,94,51,81,\adfsplit 106,70,158,151,76,175,123,154,113,190,18,42)_{\adfGD}$,

$(120,179,55,167,45,63,136,62,44,95,52,82,\adfsplit 107,71,159,152,77,176,124,155,114,193,19,43)_{\adfGD}$,

$(121,144,56,168,46,64,137,63,45,96,53,83,\adfsplit 108,72,160,153,78,177,125,156,115,196,20,44)_{\adfGD}$,

$(122,145,57,169,47,65,138,64,46,97,54,84,\adfsplit 109,73,161,154,79,178,126,157,116,199,21,45)_{\adfGD}$,

$(127,154,190,158,109,32,34,12,140,67,199,174,\adfsplit 9,173,81,38,61,114,46,0,65,198,139,147)_{\adfGD}$,

$(128,155,193,159,110,33,35,13,141,68,202,175,\adfsplit 10,174,82,39,62,115,47,1,66,201,140,148)_{\adfGD}$,

$(129,156,196,160,111,34,36,14,142,69,205,176,\adfsplit 11,175,83,40,63,116,48,2,67,204,141,149)_{\adfGD}$,

$(130,157,199,161,112,35,37,15,143,70,181,177,\adfsplit 12,176,84,41,64,117,49,3,68,180,142,150)_{\adfGD}$,

$(158,49,92,126,34,200,191,45,201,77,122,142,\adfsplit 183,182,36,186,85,55,57,42,44,180,78,185)_{\adfGD}$,

$(160,20,87,58,67,97,44,197,188,200,120,165,\adfsplit 63,56,41,25,99,92,185,201,122,198,105,79)_{\adfGD}$,

$(74,61,59,206,185,138,46,10,20,123,15,17,\adfsplit 155,107,186,0,60,182,89,188,198,127,32,23)_{\adfGD}$,

\adfLgap
$(205,2,121,86,87,178,30,60,156,9,18,80,\adfsplit 195,21,100,145,109,53,142,180,37,12,150,141)_{\adfGE}$,

$(181,3,122,87,88,179,31,61,157,10,19,81,\adfsplit 198,22,101,146,110,54,143,183,38,13,151,142)_{\adfGE}$,

$(184,4,123,88,89,144,32,62,158,11,20,82,\adfsplit 201,23,102,147,111,55,0,186,39,14,152,143)_{\adfGE}$,

$(187,5,124,89,90,145,33,63,159,12,21,83,\adfsplit 204,24,103,148,112,56,1,189,40,15,153,0)_{\adfGE}$,

$(181,10,108,110,1,158,174,143,147,55,127,109,\adfsplit 130,30,99,6,183,166,165,49,157,35,34,86)_{\adfGE}$,

$(184,11,109,111,2,159,175,0,148,56,128,110,\adfsplit 131,31,100,7,186,167,166,50,158,36,35,87)_{\adfGE}$,

$(187,12,110,112,3,160,176,1,149,57,129,111,\adfsplit 132,32,101,8,189,168,167,51,159,37,36,88)_{\adfGE}$,

$(190,13,111,113,4,161,177,2,150,58,130,112,\adfsplit 133,33,102,9,192,169,168,52,160,38,37,89)_{\adfGE}$,

$(141,70,158,42,37,63,140,131,188,203,145,30,\adfsplit 20,182,15,91,144,205,191,43,88,105,93,32)_{\adfGE}$,

$(165,66,51,126,148,155,94,56,121,194,101,77,\adfsplit 185,45,188,0,182,32,22,184,91,132,199,125)_{\adfGE}$,

$(110,190,11,7,91,64,182,158,108,2,175,97,\adfsplit 40,8,197,129,86,90,47,107,17,194,188,118)_{\adfGE}$,

\adfLgap
$(94,151,7,193,78,133,166,187,35,16,161,108,\adfsplit 91,90,195,89,139,19,82,186,202,174,40,135)_{\adfGF}$,

$(95,152,8,196,79,134,167,190,36,17,162,109,\adfsplit 92,91,198,90,140,20,83,189,205,175,41,136)_{\adfGF}$,

$(96,153,9,199,80,135,168,193,37,18,163,110,\adfsplit 93,92,201,91,141,21,84,192,181,176,42,137)_{\adfGF}$,

$(97,154,10,202,81,136,169,196,38,19,164,111,\adfsplit 94,93,204,92,142,22,85,195,184,177,43,138)_{\adfGF}$,

$(28,5,172,166,154,82,122,62,135,15,103,180,\adfsplit 179,137,147,78,168,45,130,76,170,52,73,65)_{\adfGF}$,

$(29,6,173,167,155,83,123,63,136,16,104,183,\adfsplit 144,138,148,79,169,46,131,77,171,53,74,66)_{\adfGF}$,

$(30,7,174,168,156,84,124,64,137,17,105,186,\adfsplit 145,139,149,80,170,47,132,78,172,54,75,67)_{\adfGF}$,

$(31,8,175,169,157,85,125,65,138,18,106,189,\adfsplit 146,140,150,81,171,48,133,79,173,55,76,68)_{\adfGF}$,

$(153,91,21,90,108,6,194,203,119,89,39,174,\adfsplit 0,73,84,4,132,75,49,102,99,188,182,100)_{\adfGF}$,

$(167,143,32,141,96,94,49,61,200,2,133,79,\adfsplit 28,102,67,185,106,52,182,53,99,111,203,110)_{\adfGF}$,

$(48,144,182,7,9,118,42,29,36,188,135,79,\adfsplit 24,18,77,126,200,64,206,59,17,49,110,34)_{\adfGF}$,

\adfLgap
$(43,160,163,171,133,23,54,8,38,77,62,150,\adfsplit 81,101,154,0,165,76,102,189,17,135,18,166)_{\adfGG}$,

$(44,161,164,172,134,24,55,9,39,78,63,151,\adfsplit 82,102,155,1,166,77,103,192,18,136,19,167)_{\adfGG}$,

$(45,162,165,173,135,25,56,10,40,79,64,152,\adfsplit 83,103,156,2,167,78,104,195,19,137,20,168)_{\adfGG}$,

$(46,163,166,174,136,26,57,11,41,80,65,153,\adfsplit 84,104,157,3,168,79,105,198,20,138,21,169)_{\adfGG}$,

$(67,205,175,36,90,14,35,93,29,188,89,172,\adfsplit 134,178,116,129,149,23,105,55,150,8,115,189)_{\adfGG}$,

$(68,181,176,37,91,15,36,94,30,191,90,173,\adfsplit 135,179,117,130,150,24,106,56,151,9,116,192)_{\adfGG}$,

$(69,184,177,38,92,16,37,95,31,194,91,174,\adfsplit 136,144,118,131,151,25,107,57,152,10,117,195)_{\adfGG}$,

$(70,187,178,39,93,17,38,96,32,197,92,175,\adfsplit 137,145,119,132,152,26,108,58,153,11,118,198)_{\adfGG}$,

$(89,122,201,205,87,59,52,143,78,43,84,32,\adfsplit 184,191,111,197,190,139,58,133,14,96,115,193)_{\adfGG}$,

$(114,101,33,7,0,74,194,205,196,34,202,190,\adfsplit 137,52,85,77,32,122,112,115,24,191,27,185)_{\adfGG}$,

$(1,36,6,38,197,204,41,198,35,200,143,60,\adfsplit 130,33,30,65,183,27,121,191,56,2,53,40)_{\adfGG}$,

\adfLgap
$(74,97,180,131,30,157,133,54,198,187,14,73,\adfsplit 3,175,46,39,149,34,24,107,99,197,171,92)_{\adfGH}$,

$(75,98,183,132,31,158,134,55,201,190,15,74,\adfsplit 4,176,47,40,150,35,25,108,100,200,172,93)_{\adfGH}$,

$(76,99,186,133,32,159,135,56,204,193,16,75,\adfsplit 5,177,48,41,151,36,26,109,101,203,173,94)_{\adfGH}$,

$(77,100,189,134,33,160,136,57,180,196,17,76,\adfsplit 6,178,49,42,152,37,27,110,102,206,174,95)_{\adfGH}$,

$(97,140,186,178,179,131,129,104,55,40,68,52,\adfsplit 67,171,28,25,205,158,69,66,115,177,144,62)_{\adfGH}$,

$(98,141,189,179,144,132,130,105,56,41,69,53,\adfsplit 68,172,29,26,181,159,70,67,116,178,145,63)_{\adfGH}$,

$(99,142,192,144,145,133,131,106,57,42,70,54,\adfsplit 69,173,30,27,184,160,71,68,117,179,146,64)_{\adfGH}$,

$(100,143,195,145,146,134,132,107,58,43,71,55,\adfsplit 70,174,31,28,187,161,72,69,118,144,147,65)_{\adfGH}$,

$(144,13,79,71,26,42,197,10,182,104,101,55,\adfsplit 191,172,38,19,88,128,114,17,78,5,188,59)_{\adfGH}$,

$(163,33,98,68,78,178,123,63,151,177,131,11,\adfsplit 34,84,61,22,203,181,197,59,174,80,130,109)_{\adfGH}$,

$(128,194,5,29,69,137,116,64,104,188,193,118,\adfsplit 83,199,33,112,141,177,139,40,55,8,187,22)_{\adfGH}$,

\adfLgap
$(23,20,153,147,157,13,58,76,19,5,140,178,\adfsplit 127,167,86,152,164,24,55,17,48,42,3,158)_{\adfGI}$,

$(24,21,154,148,158,14,59,77,20,6,141,179,\adfsplit 128,168,87,153,165,25,56,18,49,43,4,159)_{\adfGI}$,

$(25,22,155,149,159,15,60,78,21,7,142,144,\adfsplit 129,169,88,154,166,26,57,19,50,44,5,160)_{\adfGI}$,

$(26,23,156,150,160,16,61,79,22,8,143,145,\adfsplit 130,170,89,155,167,27,58,20,51,45,6,161)_{\adfGI}$,

$(131,183,34,181,88,84,168,15,45,85,187,184,\adfsplit 166,136,130,116,117,81,40,73,105,198,32,182)_{\adfGI}$,

$(132,186,35,184,89,85,169,16,46,86,190,187,\adfsplit 167,137,131,117,118,82,41,74,106,201,33,185)_{\adfGI}$,

$(133,189,36,187,90,86,170,17,47,87,193,190,\adfsplit 168,138,132,118,119,83,42,75,107,204,34,188)_{\adfGI}$,

$(134,192,37,190,91,87,171,18,48,88,196,193,\adfsplit 169,139,133,119,120,84,43,76,108,180,35,191)_{\adfGI}$,

$(163,102,123,109,131,77,191,108,104,48,78,76,\adfsplit 1,148,195,13,154,67,101,23,50,18,200,182)_{\adfGI}$,

$(49,54,14,78,91,115,180,175,195,45,0,4,\adfsplit 3,22,122,204,172,13,37,131,40,10,185,197)_{\adfGI}$,

$(35,60,72,30,194,55,57,101,3,197,120,22,\adfsplit 154,86,64,41,93,87,7,191,177,169,52,8)_{\adfGI}$,

\adfLgap
$(146,13,26,43,76,187,69,83,145,175,92,3,\adfsplit 62,157,37,57,9,160,185,122,148,40,1,105)_{\adfGJ}$,

$(147,14,27,44,77,190,70,84,146,176,93,4,\adfsplit 63,158,38,58,10,161,188,123,149,41,2,106)_{\adfGJ}$,

$(148,15,28,45,78,193,71,85,147,177,94,5,\adfsplit 64,159,39,59,11,162,191,124,150,42,3,107)_{\adfGJ}$,

$(149,16,29,46,79,196,72,86,148,178,95,6,\adfsplit 65,160,40,60,12,163,194,125,151,43,4,108)_{\adfGJ}$,

$(205,110,23,52,155,107,102,184,123,57,98,148,\adfsplit 114,71,177,145,163,127,68,118,31,27,16,186)_{\adfGJ}$,

$(181,111,24,53,156,108,103,187,124,58,99,149,\adfsplit 115,72,178,146,164,128,69,119,32,28,17,189)_{\adfGJ}$,

$(184,112,25,54,157,109,104,190,125,59,100,150,\adfsplit 116,73,179,147,165,129,70,120,33,29,18,192)_{\adfGJ}$,

$(187,113,26,55,158,110,105,193,126,60,101,151,\adfsplit 117,74,144,148,166,130,71,121,34,30,19,195)_{\adfGJ}$,

$(8,166,185,101,0,21,89,65,186,192,112,201,\adfsplit 132,28,17,19,194,131,142,182,200,53,130,141)_{\adfGJ}$,

$(164,34,114,91,198,53,61,95,2,38,76,82,\adfsplit 124,153,201,1,103,55,185,8,20,54,191,115)_{\adfGJ}$,

$(101,159,26,74,122,109,194,186,7,23,108,14,\adfsplit 39,60,128,197,127,201,51,33,103,198,12,126)_{\adfGJ}$,

\adfLgap
$(180,122,128,28,85,73,179,123,148,162,178,74,\adfsplit 64,52,88,113,142,75,127,81,183,145,147,70)_{\adfGK}$,

$(183,123,129,29,86,74,144,124,149,163,179,75,\adfsplit 65,53,89,114,143,76,128,82,186,146,148,71)_{\adfGK}$,

$(186,124,130,30,87,75,145,125,150,164,144,76,\adfsplit 66,54,90,115,0,77,129,83,189,147,149,72)_{\adfGK}$,

$(189,125,131,31,88,76,146,126,151,165,145,77,\adfsplit 67,55,91,116,1,78,130,84,192,148,150,73)_{\adfGK}$,

$(133,136,171,180,146,150,93,139,111,98,44,134,\adfsplit 68,199,28,125,174,129,83,152,197,115,19,130)_{\adfGK}$,

$(134,137,172,183,147,151,94,140,112,99,45,135,\adfsplit 69,202,29,126,175,130,84,153,200,116,20,131)_{\adfGK}$,

$(135,138,173,186,148,152,95,141,113,100,46,136,\adfsplit 70,205,30,127,176,131,85,154,203,117,21,132)_{\adfGK}$,

$(136,139,174,189,149,153,96,142,114,101,47,137,\adfsplit 71,181,31,128,177,132,86,155,206,118,22,133)_{\adfGK}$,

$(13,44,116,157,185,89,41,205,138,5,112,110,\adfsplit 196,35,21,187,126,143,197,123,46,80,188,15)_{\adfGK}$,

$(174,47,22,138,26,76,125,9,203,25,205,167,\adfsplit 199,38,28,12,104,91,51,131,113,181,182,130)_{\adfGK}$,

$(152,61,36,44,18,185,139,15,63,196,52,3,\adfsplit 64,34,58,73,62,81,193,200,206,103,137,120)_{\adfGK}$

\noindent and

$(22,158,198,144,115,24,140,134,88,108,155,125,\adfsplit 169,189,121,13,171,10,184,96,20,81,170,132)_{\adfGL}$,

$(23,159,201,145,116,25,141,135,89,109,156,126,\adfsplit 170,192,122,14,172,11,187,97,21,82,171,133)_{\adfGL}$,

$(24,160,204,146,117,26,142,136,90,110,157,127,\adfsplit 171,195,123,15,173,12,190,98,22,83,172,134)_{\adfGL}$,

$(25,161,180,147,118,27,143,137,91,111,158,128,\adfsplit 172,198,124,16,174,13,193,99,23,84,173,135)_{\adfGL}$,

$(152,47,55,64,196,12,14,72,135,11,151,168,\adfsplit 186,179,162,94,8,3,66,77,83,134,116,75)_{\adfGL}$,

$(153,48,56,65,199,13,15,73,136,12,152,169,\adfsplit 189,144,163,95,9,4,67,78,84,135,117,76)_{\adfGL}$,

$(154,49,57,66,202,14,16,74,137,13,153,170,\adfsplit 192,145,164,96,10,5,68,79,85,136,118,77)_{\adfGL}$,

$(155,50,58,67,205,15,17,75,138,14,154,171,\adfsplit 195,146,165,97,11,6,69,80,86,137,119,78)_{\adfGL}$,

$(80,57,203,199,182,78,105,68,70,95,194,33,\adfsplit 193,91,37,58,32,87,14,200,59,110,106,13)_{\adfGL}$,

$(145,51,5,134,32,116,199,104,97,166,194,53,\adfsplit 52,188,67,83,0,89,185,172,197,78,68,47)_{\adfGL}$,

$(93,38,24,72,103,127,182,179,194,185,28,197,\adfsplit 30,39,13,29,70,7,94,61,184,108,126,139)_{\adfGL}$

\noindent under the action of the mapping $x \mapsto x + 4$ (mod 144) for $x < 144$,
$x \mapsto 144 + (x + 4 \mathrm{~(mod~36)})$ for $144 \le x < 180$,
$x \mapsto 180 + (x - 180 + 12 \mathrm{~(mod~27)})$ for $x \ge 180$.

\adfVfy{207, \{\{144,36,4\},\{36,9,4\},\{27,9,12\}\}, -1, \{\{72,\{0,1\}\},\{63,\{2\}\}\}, -1} 

Let the vertex set of $K_{24,24,24,24}$ be $Z_{96}$ partitioned according to residue classes modulo 4.
The decompositions consist of

$(0,1,2,3,6,7,13,16,21,20,30,39,\adfsplit 29,50,5,66,51,91,62,19,8,82,32,53)_{\adfGa}$,

\adfLgap
$(0,1,2,3,6,7,13,16,21,20,30,26,\adfsplit 35,56,51,4,64,14,90,45,27,11,60,72)_{\adfGb}$,

\adfLgap
$(0,1,2,3,6,7,13,16,21,20,30,26,\adfsplit 35,42,51,4,8,69,71,12,5,63,59,28)_{\adfGc}$,

\adfLgap
$(0,1,2,3,6,7,13,16,21,20,30,26,\adfsplit 35,42,51,52,81,77,5,8,75,89,24,54)_{\adfGd}$,

\adfLgap
$(0,1,2,3,6,7,13,16,21,20,17,39,\adfsplit 29,50,43,68,56,4,77,26,18,91,34,67)_{\adfGe}$,

\adfLgap
$(0,1,2,3,6,7,13,16,21,20,17,39,\adfsplit 29,50,42,62,56,4,80,27,15,90,41,57)_{\adfGf}$,

\adfLgap
$(0,1,2,3,6,7,13,16,21,20,17,39,\adfsplit 29,42,47,64,63,94,4,71,69,19,28,61)_{\adfGg}$,

\adfLgap
$(0,1,2,3,6,7,13,16,21,20,17,39,\adfsplit 29,42,47,46,84,62,90,8,65,10,15,51)_{\adfGh}$,

\adfLgap
$(0,1,2,3,6,7,13,16,21,20,17,39,\adfsplit 29,42,43,64,70,63,94,4,72,14,9,55)_{\adfGi}$,

\adfLgap
$(0,1,2,3,6,7,13,16,21,20,17,26,\adfsplit 35,42,51,43,72,85,82,78,12,4,37,39)_{\adfGj}$,

\adfLgap
$(0,1,2,3,6,7,13,16,21,20,17,26,\adfsplit 35,42,51,57,66,95,65,90,88,68,14,27)_{\adfGk}$,

\adfLgap
$(0,1,2,3,6,7,13,16,21,20,17,26,\adfsplit 35,42,43,51,86,77,94,76,72,9,48,23)_{\adfGl}$,

\adfLgap
$(0,1,2,3,6,7,13,16,21,20,17,26,\adfsplit 35,42,49,45,75,82,62,90,5,64,11,40)_{\adfGm}$,

\adfLgap
$(0,1,2,3,6,7,13,11,16,28,25,24,\adfsplit 37,5,47,54,86,56,84,23,87,44,29,90)_{\adfGn}$,

\adfLgap
$(0,1,2,3,6,7,13,11,16,28,25,24,\adfsplit 30,46,50,47,73,63,88,8,76,22,26,57)_{\adfGo}$,

\adfLgap
$(0,1,2,3,6,7,13,11,16,28,25,24,\adfsplit 30,51,50,74,53,68,93,15,14,84,38,48)_{\adfGp}$,

\adfLgap
$(0,1,2,3,6,7,13,11,16,28,25,24,\adfsplit 30,37,50,79,71,57,95,87,72,20,14,54)_{\adfGq}$,

\adfLgap
$(0,1,2,3,6,7,13,11,16,28,25,24,\adfsplit 30,37,50,52,87,64,89,94,60,5,22,48)_{\adfGr}$,

\adfLgap
$(0,1,2,3,6,7,13,11,16,28,25,24,\adfsplit 30,37,5,63,66,91,80,93,32,50,19,54)_{\adfGs}$,

\adfLgap
$(0,1,2,3,6,7,13,11,16,28,21,24,\adfsplit 34,40,54,75,58,81,94,14,12,20,51,55)_{\adfGt}$,

\adfLgap
$(0,1,2,3,6,7,13,11,16,28,21,24,\adfsplit 34,40,54,45,79,56,86,15,90,80,37,52)_{\adfGu}$,

\adfLgap
$(0,1,2,3,6,7,13,11,16,28,21,24,\adfsplit 34,40,50,42,54,67,94,12,95,77,20,49)_{\adfGv}$,

\adfLgap
$(0,1,2,3,6,7,13,11,16,28,21,24,\adfsplit 34,40,50,42,67,4,73,19,8,77,26,58)_{\adfGw}$,

\adfLgap
$(0,1,2,3,6,7,13,11,16,28,21,24,\adfsplit 34,40,49,56,65,50,85,18,88,23,35,77)_{\adfGx}$,

\adfLgap
$(0,1,2,3,6,7,13,11,16,28,21,24,\adfsplit 34,30,50,51,57,86,89,83,81,20,47,36)_{\adfGy}$,

\adfLgap
$(0,1,2,3,6,7,13,11,16,28,21,24,\adfsplit 34,30,50,51,73,61,8,65,77,94,32,19)_{\adfGz}$,

\adfLgap
$(0,1,2,3,6,7,13,11,16,28,21,24,\adfsplit 34,30,50,42,57,55,95,72,91,5,29,44)_{\adfGA}$,

\adfLgap
$(0,1,2,3,6,7,13,11,16,28,21,24,\adfsplit 34,30,50,42,67,57,61,87,80,88,33,19)_{\adfGB}$,

\adfLgap
$(0,1,2,3,6,7,13,11,16,28,21,24,\adfsplit 40,29,37,58,73,93,52,90,59,95,19,45)_{\adfGC}$,

\adfLgap
$(0,1,2,3,6,7,13,11,16,28,21,24,\adfsplit 40,29,50,51,46,78,5,84,8,77,31,45)_{\adfGD}$,

\adfLgap
$(0,1,2,3,6,7,13,11,16,28,21,24,\adfsplit 40,29,50,51,77,87,79,56,91,42,33,5)_{\adfGE}$,

\adfLgap
$(0,1,2,3,6,7,13,11,16,28,21,24,\adfsplit 40,29,50,51,52,54,89,80,88,87,30,34)_{\adfGF}$,

\adfLgap
$(0,1,2,3,6,7,13,11,16,28,21,24,\adfsplit 40,29,50,51,43,61,81,67,25,10,78,31)_{\adfGG}$,

\adfLgap
$(0,1,2,3,6,7,13,11,16,28,21,24,\adfsplit 40,29,37,5,55,75,57,94,72,15,64,17)_{\adfGH}$,

\adfLgap
$(0,1,2,3,6,7,13,11,16,28,21,24,\adfsplit 29,30,50,51,56,69,73,93,82,18,40,23)_{\adfGI}$,

\adfLgap
$(0,1,2,3,6,7,13,11,16,28,21,24,\adfsplit 29,30,50,51,52,67,74,91,77,5,38,32)_{\adfGJ}$,

\adfLgap
$(0,1,2,3,6,7,13,11,16,17,22,24,\adfsplit 31,37,30,44,40,65,69,54,94,92,89,27)_{\adfGK}$

\noindent and

$(0,1,2,3,6,7,13,11,16,17,22,24,\adfsplit 31,37,30,51,48,73,71,87,90,64,52,18)_{\adfGL}$

\noindent under the action of the mapping $x \mapsto x + 1$ (mod 96).

\adfVfy{96, \{\{96,96,1\}\}, -1, \{\{24,\{0,1,2,3\}\}\}, -1} 

Let the vertex set of $K_{24,24,24,21}$ be $\{0, 1, \dots, 92\}$ partitioned into
$\{2j + i: j = 0, 1, \dots, 23\}$, $i = 0, 1, 2$, and $\{72, 73, \dots, 92\}$.
The decompositions consist of

$(64,66,24,42,61,77,65,53,35,49,69,18,\adfsplit 13,0,80,56,2,17,50,10,55,27,79,36)_{\adfGa}$,

$(85,5,22,30,16,87,90,2,72,53,4,45,\adfsplit 21,19,84,23,49,82,42,66,56,62,34,52)_{\adfGa}$,

$(75,31,46,29,3,87,76,59,90,89,10,40,\adfsplit 86,13,12,78,56,42,73,5,67,63,25,81)_{\adfGa}$,

$(52,86,11,83,10,41,77,43,51,58,42,20,\adfsplit 29,82,32,48,55,44,31,38,35,85,39,81)_{\adfGa}$,

$(36,16,74,23,24,88,71,5,61,39,38,49,\adfsplit 4,26,15,44,29,66,32,78,84,45,55,31)_{\adfGa}$,

\adfLgap
$(63,62,44,26,79,72,61,25,15,13,27,77,\adfsplit 18,28,42,17,59,24,84,81,7,33,35,8)_{\adfGb}$,

$(9,65,56,68,60,33,75,22,80,43,14,63,\adfsplit 36,16,5,70,46,86,90,81,24,27,44,1)_{\adfGb}$,

$(76,11,2,44,6,90,80,1,52,18,38,10,\adfsplit 87,53,41,74,39,57,48,8,85,49,78,22)_{\adfGb}$,

$(74,67,24,27,45,11,89,13,65,68,39,31,\adfsplit 21,75,36,80,90,38,59,5,51,30,85,55)_{\adfGb}$,

$(53,24,9,58,68,59,85,77,66,62,34,7,\adfsplit 10,32,30,47,75,84,72,1,43,18,16,36)_{\adfGb}$,

\adfLgap
$(39,86,52,87,7,25,54,36,65,43,38,53,\adfsplit 46,16,48,85,6,63,44,14,89,10,64,90)_{\adfGc}$,

$(19,0,3,56,78,25,65,89,15,36,91,12,\adfsplit 43,50,4,41,2,32,6,7,47,45,84,70)_{\adfGc}$,

$(48,88,65,26,5,61,43,9,90,49,35,15,\adfsplit 7,28,11,81,18,82,79,41,70,24,80,14)_{\adfGc}$,

$(39,74,40,58,1,54,63,84,82,21,28,23,\adfsplit 5,64,86,76,45,83,80,50,16,17,78,43)_{\adfGc}$,

$(42,22,87,26,53,41,16,0,86,21,49,81,\adfsplit 68,12,38,6,51,3,84,85,39,17,59,10)_{\adfGc}$,

\adfLgap
$(59,75,83,46,32,5,60,21,78,47,77,8,\adfsplit 81,39,82,4,74,25,67,50,71,38,66,15)_{\adfGd}$,

$(11,48,16,61,26,74,76,59,6,9,62,58,\adfsplit 78,32,79,36,72,60,69,87,73,29,2,7)_{\adfGd}$,

$(16,47,17,20,9,31,46,3,43,52,33,24,\adfsplit 56,82,41,87,14,63,85,1,29,49,22,45)_{\adfGd}$,

$(89,14,68,22,7,10,77,54,83,17,63,2,\adfsplit 78,15,19,52,84,50,43,1,79,62,49,60)_{\adfGd}$,

$(84,70,17,55,23,14,48,81,80,56,54,1,\adfsplit 47,26,33,65,36,39,79,13,60,90,52,8)_{\adfGd}$,

\adfLgap
$(37,78,45,41,35,16,74,72,81,48,69,22,\adfsplit 26,60,85,34,39,90,4,76,73,32,51,14)_{\adfGe}$,

$(12,56,35,14,13,25,76,39,69,31,5,33,\adfsplit 84,7,15,67,1,50,79,8,23,66,19,28)_{\adfGe}$,

$(6,70,75,85,9,3,44,67,62,11,22,89,\adfsplit 35,46,40,81,42,47,30,69,2,8,73,64)_{\adfGe}$,

$(4,29,30,74,61,84,72,26,15,50,40,0,\adfsplit 13,59,79,36,85,80,33,73,67,45,2,17)_{\adfGe}$,

$(2,39,1,22,4,62,3,91,17,75,48,64,\adfsplit 24,16,35,84,76,33,49,54,15,5,82,56)_{\adfGe}$,

\adfLgap
$(23,84,82,72,61,27,39,34,2,4,71,62,\adfsplit 78,60,90,86,38,5,48,66,28,67,87,74)_{\adfGf}$,

$(13,30,59,89,40,50,79,28,32,47,18,45,\adfsplit 44,37,80,29,15,69,48,33,81,88,35,62)_{\adfGf}$,

$(22,35,59,84,55,51,78,40,20,18,90,44,\adfsplit 63,1,53,42,29,77,52,61,27,79,12,17)_{\adfGf}$,

$(88,8,10,41,77,45,15,32,46,49,22,73,\adfsplit 58,62,3,30,24,61,67,64,53,20,9,78)_{\adfGf}$,

$(26,51,70,76,52,11,75,15,47,0,42,49,\adfsplit 65,92,82,87,34,50,1,56,32,9,90,39)_{\adfGf}$,

\adfLgap
$(44,66,72,16,61,38,53,9,63,54,85,22,\adfsplit 88,65,67,32,4,50,3,89,76,12,69,84)_{\adfGg}$,

$(15,7,83,41,30,68,38,0,76,12,45,73,\adfsplit 50,61,35,13,87,3,42,77,48,89,17,47)_{\adfGg}$,

$(2,77,6,18,48,42,64,8,88,44,65,73,\adfsplit 89,71,7,31,58,54,51,37,0,72,85,11)_{\adfGg}$,

$(82,0,15,71,43,58,80,32,21,52,47,57,\adfsplit 64,34,85,10,17,2,3,86,65,88,45,19)_{\adfGg}$,

$(29,91,30,25,15,0,11,37,88,8,31,10,\adfsplit 24,68,79,92,5,19,4,83,14,49,66,65)_{\adfGg}$,

\adfLgap
$(39,29,78,88,67,69,8,46,17,14,75,83,\adfsplit 0,85,66,61,41,62,23,27,16,52,33,74)_{\adfGh}$,

$(84,67,45,46,21,54,34,32,5,15,58,90,\adfsplit 82,80,14,88,25,63,60,38,20,79,16,43)_{\adfGh}$,

$(79,44,33,10,13,76,8,23,17,74,30,24,\adfsplit 51,73,49,84,32,62,29,43,66,12,64,75)_{\adfGh}$,

$(12,47,72,7,87,90,19,56,8,82,17,2,\adfsplit 15,66,18,25,84,53,68,78,40,11,51,69)_{\adfGh}$,

$(54,49,62,22,80,56,6,37,60,12,33,55,\adfsplit 77,92,50,19,23,47,32,63,10,72,39,61)_{\adfGh}$,

\adfLgap
$(63,40,83,65,0,62,41,7,74,37,19,79,\adfsplit 76,39,89,50,6,21,56,61,77,24,16,81)_{\adfGi}$,

$(33,35,83,44,43,3,54,50,85,21,10,29,\adfsplit 61,65,14,87,48,24,51,70,31,12,32,41)_{\adfGi}$,

$(77,45,6,1,88,79,46,22,18,69,64,89,\adfsplit 23,32,87,20,11,27,49,72,30,55,84,71)_{\adfGi}$,

$(4,59,54,83,1,16,50,85,20,22,11,9,\adfsplit 36,40,8,78,39,23,61,80,44,55,75,48)_{\adfGi}$,

$(79,32,35,30,43,25,9,82,38,91,87,8,\adfsplit 10,66,37,18,2,81,85,47,14,17,51,33)_{\adfGi}$,

\adfLgap
$(71,16,90,37,48,82,29,56,84,47,43,74,\adfsplit 52,51,36,34,45,69,89,30,67,59,70,78)_{\adfGj}$,

$(21,62,13,73,34,30,74,90,17,54,23,26,\adfsplit 19,83,80,49,51,79,50,44,58,18,85,0)_{\adfGj}$,

$(62,84,72,39,40,23,63,37,73,52,22,88,\adfsplit 60,32,2,65,24,64,86,53,51,21,14,80)_{\adfGj}$,

$(20,81,25,49,43,27,86,78,29,24,5,32,\adfsplit 56,13,67,4,18,76,66,48,75,91,71,7)_{\adfGj}$,

$(64,85,48,29,8,69,10,43,7,89,52,84,\adfsplit 23,42,6,61,18,65,57,63,38,59,34,37)_{\adfGj}$,

\adfLgap
$(42,72,34,10,8,29,9,85,60,53,21,67,\adfsplit 0,81,86,64,65,3,41,70,19,47,48,77)_{\adfGk}$,

$(8,31,4,37,51,65,29,84,26,75,73,74,\adfsplit 68,22,45,12,28,2,33,85,78,66,47,59)_{\adfGk}$,

$(33,90,59,87,27,36,82,31,52,18,17,68,\adfsplit 85,79,13,55,66,88,69,4,71,42,62,89)_{\adfGk}$,

$(60,20,65,14,3,69,13,10,87,73,34,53,\adfsplit 26,44,51,11,89,37,67,48,79,76,27,62)_{\adfGk}$,

$(12,4,25,46,23,11,66,6,91,84,22,64,\adfsplit 74,41,36,39,76,35,16,1,3,90,26,2)_{\adfGk}$,

\adfLgap
$(1,6,3,26,32,49,81,5,73,15,56,35,\adfsplit 67,48,80,90,42,9,47,36,86,28,82,25)_{\adfGl}$,

$(53,85,76,46,23,28,9,62,42,63,71,90,\adfsplit 86,50,84,31,21,12,17,68,80,10,88,25)_{\adfGl}$,

$(53,91,4,0,66,8,78,27,85,26,25,65,\adfsplit 41,31,88,60,47,40,54,30,2,35,70,77)_{\adfGl}$,

$(11,55,88,4,14,39,42,44,90,6,79,86,\adfsplit 54,32,5,28,18,51,67,89,33,43,68,13)_{\adfGl}$,

$(67,66,73,32,19,29,27,6,36,0,10,75,\adfsplit 92,52,82,28,65,49,43,17,50,21,84,9)_{\adfGl}$,

\adfLgap
$(9,8,20,64,48,46,88,73,84,36,5,54,\adfsplit 7,66,35,61,83,87,53,30,89,22,69,67)_{\adfGm}$,

$(36,13,76,44,48,20,71,7,78,51,84,70,\adfsplit 62,4,86,16,59,55,91,40,6,14,57,82)_{\adfGm}$,

$(37,75,9,30,43,8,71,80,81,40,82,60,\adfsplit 58,47,86,54,22,69,72,21,0,31,53,16)_{\adfGm}$,

$(43,79,30,91,1,22,71,85,41,37,2,33,\adfsplit 39,88,50,23,4,7,5,65,9,67,27,78)_{\adfGm}$,

$(28,24,80,14,4,55,54,41,30,83,3,65,\adfsplit 74,92,46,68,76,50,8,6,21,71,63,1)_{\adfGm}$,

\adfLgap
$(57,61,14,86,56,84,66,54,63,19,71,24,\adfsplit 1,47,58,37,88,48,62,60,68,18,76,32)_{\adfGn}$,

$(83,23,17,70,58,54,85,61,32,47,35,53,\adfsplit 69,67,25,28,2,89,41,72,11,49,27,68)_{\adfGn}$,

$(78,29,59,20,55,18,81,74,84,40,1,72,\adfsplit 4,48,38,39,80,11,25,87,31,53,90,24)_{\adfGn}$,

$(90,54,41,24,86,11,10,77,17,82,6,52,\adfsplit 49,53,74,20,32,3,21,73,16,66,75,2)_{\adfGn}$,

$(41,21,7,75,23,32,87,91,34,2,22,51,\adfsplit 6,4,38,53,64,71,78,46,92,35,60,39)_{\adfGn}$,

\adfLgap
$(45,87,49,55,11,41,72,42,63,50,16,27,\adfsplit 34,82,40,20,57,38,51,64,3,59,4,76)_{\adfGo}$,

$(28,23,79,32,34,51,50,61,69,74,75,20,\adfsplit 8,31,37,63,57,88,86,62,85,58,60,39)_{\adfGo}$,

$(63,88,25,80,61,19,53,36,2,42,59,66,\adfsplit 73,24,11,89,56,16,84,14,22,91,17,0)_{\adfGo}$,

$(86,16,14,23,6,8,43,34,4,91,76,60,\adfsplit 5,9,17,78,31,37,85,51,62,18,21,1)_{\adfGo}$,

$(77,22,66,29,80,56,59,26,9,88,85,49,\adfsplit 90,71,30,76,3,25,8,50,36,5,82,64)_{\adfGo}$,

\adfLgap
$(40,23,5,62,87,70,42,39,13,52,71,77,\adfsplit 79,24,38,88,30,64,31,33,17,81,36,85)_{\adfGp}$,

$(11,48,55,27,19,8,14,29,20,90,31,28,\adfsplit 15,1,30,91,82,37,79,57,44,25,33,10)_{\adfGp}$,

$(79,32,23,14,6,84,78,91,58,43,40,71,\adfsplit 29,3,12,53,51,90,4,19,70,82,0,87)_{\adfGp}$,

$(71,74,9,80,69,67,64,38,53,8,6,63,\adfsplit 3,83,55,19,65,26,75,79,10,13,18,41)_{\adfGp}$,

$(87,62,8,65,81,75,34,57,64,36,46,37,\adfsplit 13,85,47,12,18,44,28,26,22,25,90,78)_{\adfGp}$,

\adfLgap
$(9,65,73,32,13,24,71,37,86,88,8,56,\adfsplit 0,50,58,18,28,70,91,77,25,83,21,57)_{\adfGq}$,

$(61,53,79,27,46,34,0,56,20,67,75,23,\adfsplit 39,40,30,50,9,83,8,2,91,81,59,51)_{\adfGq}$,

$(84,3,47,70,86,75,19,78,57,11,6,23,\adfsplit 55,62,27,73,43,79,44,54,45,24,14,10)_{\adfGq}$,

$(21,75,43,80,5,13,83,38,6,20,39,74,\adfsplit 60,31,42,87,33,78,85,57,50,52,61,44)_{\adfGq}$,

$(27,65,17,4,33,90,28,15,8,3,52,38,\adfsplit 70,74,34,69,9,54,56,64,43,71,85,82)_{\adfGq}$,

\adfLgap
$(85,15,18,12,89,47,25,16,10,2,74,37,\adfsplit 62,68,86,44,72,77,4,52,79,19,7,24)_{\adfGr}$,

$(87,11,55,20,24,70,41,50,82,39,69,63,\adfsplit 37,67,10,34,59,36,9,17,47,66,58,88)_{\adfGr}$,

$(79,61,2,18,50,76,46,73,67,26,66,47,\adfsplit 41,23,63,77,57,15,70,52,81,37,7,60)_{\adfGr}$,

$(16,48,27,82,7,90,76,84,5,14,52,1,\adfsplit 38,89,9,3,41,0,60,26,67,87,78,64)_{\adfGr}$,

$(1,83,42,68,50,52,25,81,85,9,45,35,\adfsplit 28,70,92,49,91,53,69,39,73,29,16,40)_{\adfGr}$,

\adfLgap
$(64,91,84,57,9,24,11,46,4,76,80,19,\adfsplit 17,12,70,14,29,6,49,63,75,44,18,23)_{\adfGs}$,

$(59,18,63,83,52,34,73,56,0,29,12,15,\adfsplit 13,35,76,72,44,14,91,6,28,36,67,85)_{\adfGs}$,

$(47,58,9,25,91,48,29,22,6,80,89,52,\adfsplit 82,20,40,21,5,12,19,35,30,88,33,79)_{\adfGs}$,

$(81,40,39,31,82,29,61,34,48,78,27,16,\adfsplit 14,67,57,87,33,3,79,46,41,85,56,88)_{\adfGs}$,

$(86,61,45,41,38,66,81,7,37,58,68,70,\adfsplit 42,63,44,43,6,23,76,71,78,13,50,15)_{\adfGs}$,

\adfLgap
$(12,5,58,25,60,76,41,85,72,82,70,24,\adfsplit 66,59,62,4,84,26,87,0,3,32,75,38)_{\adfGt}$,

$(67,68,23,5,63,78,16,69,91,0,13,21,\adfsplit 35,81,83,7,79,76,28,59,26,56,6,88)_{\adfGt}$,

$(65,67,18,87,26,76,64,53,54,21,49,28,\adfsplit 13,30,29,51,74,59,55,11,25,71,66,0)_{\adfGt}$,

$(29,60,73,52,89,80,26,11,15,85,56,59,\adfsplit 19,70,3,86,0,74,77,64,10,61,62,66)_{\adfGt}$,

$(63,79,62,75,66,49,1,12,53,71,65,92,\adfsplit 77,43,28,36,35,2,29,68,6,25,69,58)_{\adfGt}$,

\adfLgap
$(23,27,52,82,61,44,42,56,41,43,81,80,\adfsplit 64,9,11,13,59,51,16,57,62,53,45,7)_{\adfGu}$,

$(77,61,71,66,84,73,69,88,7,10,29,23,\adfsplit 33,46,18,74,28,49,60,11,8,30,78,87)_{\adfGu}$,

$(38,87,1,36,10,66,90,82,53,67,28,16,\adfsplit 52,60,42,14,75,17,86,84,58,64,2,23)_{\adfGu}$,

$(25,89,26,20,38,66,58,22,72,83,23,32,\adfsplit 29,48,37,76,13,19,47,71,3,60,85,78)_{\adfGu}$,

$(2,7,81,79,0,71,13,16,19,41,22,61,\adfsplit 68,77,78,59,76,24,8,14,72,9,34,55)_{\adfGu}$,

\adfLgap
$(30,81,35,56,23,33,70,76,77,69,62,79,\adfsplit 40,19,38,34,49,60,5,48,3,83,2,55)_{\adfGv}$,

$(42,59,78,25,51,61,38,3,85,87,33,81,\adfsplit 47,43,63,86,31,91,41,52,56,2,69,80)_{\adfGv}$,

$(0,26,35,70,73,72,51,54,57,50,71,16,\adfsplit 22,84,82,28,75,18,52,24,38,5,67,81)_{\adfGv}$,

$(31,8,12,57,63,90,70,20,25,84,9,58,\adfsplit 45,13,50,62,82,22,44,80,1,51,29,0)_{\adfGv}$,

$(6,4,22,74,65,86,78,60,1,16,15,47,\adfsplit 44,56,48,82,61,9,92,35,45,31,83,40)_{\adfGv}$,

\adfLgap
$(14,22,86,43,33,78,5,66,56,85,32,23,\adfsplit 87,70,25,90,30,27,38,73,17,45,55,16)_{\adfGw}$,

$(30,2,53,82,48,22,16,12,32,26,91,6,\adfsplit 34,67,25,11,78,83,36,33,64,88,31,62)_{\adfGw}$,

$(90,63,1,43,14,7,15,75,0,53,21,82,\adfsplit 65,23,40,77,22,61,86,20,29,55,32,88)_{\adfGw}$,

$(58,84,71,39,65,38,25,3,34,47,76,87,\adfsplit 88,28,10,27,7,54,74,5,59,61,8,79)_{\adfGw}$,

$(43,80,65,2,45,36,10,60,28,24,61,40,\adfsplit 86,82,77,14,29,7,35,8,78,52,9,0)_{\adfGw}$,

\adfLgap
$(46,6,60,38,67,86,87,49,51,4,35,22,\adfsplit 75,90,15,69,57,83,68,47,1,10,78,39)_{\adfGx}$,

$(86,57,14,20,5,73,39,75,77,85,46,25,\adfsplit 43,53,71,26,6,35,51,40,70,84,8,42)_{\adfGx}$,

$(86,15,29,36,41,52,16,46,1,23,85,69,\adfsplit 76,30,43,38,88,65,58,39,57,33,35,34)_{\adfGx}$,

$(40,18,81,69,14,73,60,32,64,16,36,7,\adfsplit 75,28,77,35,26,88,44,37,6,12,89,53)_{\adfGx}$,

$(58,88,85,27,26,25,69,31,73,77,53,32,\adfsplit 5,3,60,4,52,76,1,35,7,50,86,40)_{\adfGx}$,

\adfLgap
$(23,49,81,24,14,86,0,43,90,50,3,20,\adfsplit 33,71,79,37,58,27,85,57,46,2,10,5)_{\adfGy}$,

$(72,45,3,49,84,13,10,56,26,73,71,88,\adfsplit 9,31,42,82,44,58,68,91,32,48,87,60)_{\adfGy}$,

$(73,11,51,68,83,22,28,25,49,36,75,35,\adfsplit 21,8,26,60,67,15,46,30,45,56,76,34)_{\adfGy}$,

$(51,91,83,2,59,69,9,42,78,89,67,64,\adfsplit 37,10,5,8,38,21,39,81,36,85,44,29)_{\adfGy}$,

$(73,57,55,17,70,1,50,63,91,16,30,36,\adfsplit 28,81,71,37,79,78,89,8,67,42,43,51)_{\adfGy}$,

\adfLgap
$(31,59,53,47,80,12,9,33,34,87,82,37,\adfsplit 72,41,28,7,77,36,54,85,63,62,64,58)_{\adfGz}$,

$(91,43,52,65,63,57,44,36,82,9,25,58,\adfsplit 14,16,32,54,88,41,49,83,37,79,56,24)_{\adfGz}$,

$(41,72,31,83,36,3,8,77,43,54,2,32,\adfsplit 14,51,90,58,61,0,80,56,27,11,74,78)_{\adfGz}$,

$(89,17,12,25,66,10,46,77,47,50,18,51,\adfsplit 58,41,88,37,14,49,0,55,36,11,73,90)_{\adfGz}$,

$(35,81,31,48,37,36,63,33,50,92,40,46,\adfsplit 22,43,29,86,78,89,53,12,34,3,11,2)_{\adfGz}$,

\adfLgap
$(9,43,74,4,48,29,55,51,5,91,49,81,\adfsplit 64,85,26,89,14,22,60,31,57,39,80,83)_{\adfGA}$,

$(61,21,12,18,84,49,26,28,37,16,85,56,\adfsplit 45,71,69,64,34,67,11,83,13,78,2,30)_{\adfGA}$,

$(71,27,31,7,23,37,86,47,84,5,11,30,\adfsplit 16,46,74,61,83,79,53,73,44,59,64,89)_{\adfGA}$,

$(0,74,58,19,26,28,36,50,6,68,76,87,\adfsplit 75,78,51,69,8,31,46,1,25,73,11,62)_{\adfGA}$,

$(30,17,25,2,67,48,26,75,42,31,4,37,\adfsplit 8,40,68,3,72,86,6,14,84,28,18,53)_{\adfGA}$,

\adfLgap
$(50,16,15,64,35,87,89,79,78,24,26,6,\adfsplit 9,70,4,74,76,2,65,30,85,28,57,18)_{\adfGB}$,

$(40,88,54,15,59,45,46,64,16,84,73,37,\adfsplit 42,83,71,3,24,35,14,55,18,1,75,13)_{\adfGB}$,

$(40,35,27,36,34,86,11,23,91,25,44,24,\adfsplit 57,9,45,83,59,65,81,71,8,42,80,61)_{\adfGB}$,

$(86,18,51,5,89,91,22,82,49,12,23,16,\adfsplit 0,37,40,67,32,14,8,78,57,3,70,31)_{\adfGB}$,

$(3,83,81,49,17,9,64,38,73,84,26,70,\adfsplit 36,63,33,72,78,18,43,32,1,61,42,11)_{\adfGB}$,

\adfLgap
$(66,88,34,1,36,23,71,59,87,27,30,81,\adfsplit 52,37,54,4,57,39,41,17,25,80,72,70)_{\adfGC}$,

$(34,18,89,83,70,77,60,29,24,56,47,25,\adfsplit 39,15,72,67,80,75,55,37,42,17,71,90)_{\adfGC}$,

$(23,48,64,91,76,84,14,54,37,58,13,32,\adfsplit 41,71,78,15,57,89,4,18,19,36,46,81)_{\adfGC}$,

$(71,16,76,87,17,84,58,55,40,56,26,35,\adfsplit 82,68,11,3,45,51,52,54,79,78,70,53)_{\adfGC}$,

$(72,56,63,13,70,37,2,29,62,85,8,81,\adfsplit 36,64,78,46,28,65,33,26,39,44,48,61)_{\adfGC}$,

\adfLgap
$(21,52,37,67,32,82,90,45,3,53,5,38,\adfsplit 40,91,50,31,19,30,66,18,56,81,65,58)_{\adfGD}$,

$(50,77,72,86,63,21,41,67,44,43,40,70,\adfsplit 23,87,57,2,80,24,61,13,75,71,78,68)_{\adfGD}$,

$(8,84,87,67,59,36,52,10,83,18,81,46,\adfsplit 14,24,17,28,90,57,61,38,27,64,72,43)_{\adfGD}$,

$(39,49,59,83,69,30,85,81,20,26,52,86,\adfsplit 3,37,42,63,4,71,78,8,75,65,33,40)_{\adfGD}$,

$(78,64,44,62,84,0,39,69,33,89,16,87,\adfsplit 14,38,31,17,23,42,21,27,30,10,81,59)_{\adfGD}$,

\adfLgap
$(17,46,67,76,30,33,88,81,62,28,85,61,\adfsplit 71,8,55,36,35,66,57,80,24,43,64,2)_{\adfGE}$,

$(29,90,75,13,69,26,23,18,62,80,6,16,\adfsplit 53,86,70,17,89,57,34,8,27,59,47,79)_{\adfGE}$,

$(87,64,8,2,17,90,60,69,91,1,63,20,\adfsplit 38,49,9,65,88,70,43,48,31,50,42,84)_{\adfGE}$,

$(23,40,13,52,68,6,72,32,66,39,78,43,\adfsplit 22,36,20,62,75,9,18,89,27,10,84,56)_{\adfGE}$,

$(34,48,76,88,78,71,35,27,1,33,84,55,\adfsplit 32,41,36,16,51,86,43,37,92,23,5,87)_{\adfGE}$,

\adfLgap
$(20,64,49,74,14,89,51,65,57,37,25,47,\adfsplit 63,78,40,82,18,86,85,42,58,62,8,66)_{\adfGF}$,

$(5,66,64,22,65,80,6,32,77,89,9,13,\adfsplit 72,76,36,48,56,75,29,23,59,0,43,91)_{\adfGF}$,

$(8,4,70,87,57,71,65,59,52,15,21,0,\adfsplit 25,40,6,7,78,84,90,11,14,27,55,46)_{\adfGF}$,

$(8,63,78,9,50,53,34,61,55,35,64,22,\adfsplit 33,12,87,45,79,57,14,43,70,5,47,74)_{\adfGF}$,

$(21,91,53,87,58,62,83,90,6,11,15,4,\adfsplit 47,60,35,72,74,2,12,88,38,0,55,64)_{\adfGF}$,

\adfLgap
$(83,1,53,71,66,27,61,31,64,12,68,65,\adfsplit 79,30,38,22,52,10,42,74,2,63,76,41)_{\adfGG}$,

$(69,1,16,28,29,44,74,21,54,23,0,58,\adfsplit 70,87,68,18,34,85,15,14,89,3,71,13)_{\adfGG}$,

$(68,60,13,28,84,80,66,29,17,90,35,69,\adfsplit 85,72,49,19,51,32,62,36,89,91,37,12)_{\adfGG}$,

$(31,48,41,57,47,61,81,51,55,37,80,82,\adfsplit 7,13,84,36,38,35,3,50,75,22,49,79)_{\adfGG}$,

$(30,77,76,64,29,10,36,15,78,81,3,79,\adfsplit 40,55,26,14,23,48,73,85,24,52,46,71)_{\adfGG}$,

\adfLgap
$(41,57,63,36,86,55,31,35,11,85,3,73,\adfsplit 49,40,64,8,70,24,2,33,21,89,90,20)_{\adfGH}$,

$(63,71,89,62,69,74,31,13,72,82,22,43,\adfsplit 66,90,46,1,60,14,84,48,41,9,20,37)_{\adfGH}$,

$(90,50,39,10,57,74,62,56,79,84,49,17,\adfsplit 80,4,23,65,73,7,67,36,45,77,88,8)_{\adfGH}$,

$(57,49,13,34,2,85,78,65,45,17,70,71,\adfsplit 66,88,11,27,80,64,55,62,28,42,15,26)_{\adfGH}$,

$(25,68,36,71,86,91,40,47,12,48,10,24,\adfsplit 66,30,58,89,44,34,78,22,90,23,0,43)_{\adfGH}$,

\adfLgap
$(1,35,69,89,67,78,44,38,8,66,58,4,\adfsplit 43,74,81,85,76,30,5,57,60,40,50,64)_{\adfGI}$,

$(3,10,58,78,59,5,24,47,45,43,55,66,\adfsplit 16,87,83,17,14,29,22,30,89,84,69,38)_{\adfGI}$,

$(73,40,32,35,89,81,39,7,22,31,38,29,\adfsplit 88,6,51,90,52,79,3,84,69,20,49,13)_{\adfGI}$,

$(66,64,41,80,91,84,24,22,9,65,67,35,\adfsplit 29,54,1,51,42,86,72,50,14,68,48,0)_{\adfGI}$,

$(37,54,38,32,25,73,86,19,27,3,67,0,\adfsplit 20,9,82,47,69,28,61,31,36,92,44,83)_{\adfGI}$,

\adfLgap
$(34,20,18,27,43,9,29,47,87,52,61,49,\adfsplit 75,1,25,2,42,51,32,26,30,90,53,72)_{\adfGJ}$,

$(49,76,21,27,30,47,26,5,84,88,66,46,\adfsplit 85,91,20,1,32,59,71,4,63,40,72,90)_{\adfGJ}$,

$(79,16,57,71,45,48,32,23,84,15,44,59,\adfsplit 19,87,50,10,70,27,39,18,65,85,28,89)_{\adfGJ}$,

$(14,3,76,78,81,47,28,5,29,8,91,61,\adfsplit 69,1,86,9,38,73,89,30,43,52,24,56)_{\adfGJ}$,

$(42,75,19,59,2,40,85,65,22,49,30,11,\adfsplit 12,15,80,81,52,88,16,8,26,21,28,18)_{\adfGJ}$,

\adfLgap
$(2,51,7,64,26,28,66,9,89,77,68,76,\adfsplit 29,5,31,46,53,4,17,60,78,84,86,67)_{\adfGK}$,

$(8,3,81,80,55,44,39,71,66,41,25,7,\adfsplit 73,52,18,74,57,17,76,75,56,14,28,21)_{\adfGK}$,

$(20,21,36,7,67,2,83,64,68,53,72,80,\adfsplit 63,57,44,70,12,29,46,51,85,74,19,48)_{\adfGK}$,

$(41,87,66,79,56,3,22,88,52,47,34,32,\adfsplit 59,26,10,84,39,75,91,54,15,46,53,5)_{\adfGK}$,

$(19,62,18,36,58,1,32,2,79,75,3,14,\adfsplit 85,67,25,13,17,37,76,78,41,63,12,46)_{\adfGK}$

\noindent and

$(14,76,86,66,49,45,39,52,89,88,79,50,\adfsplit 22,57,15,28,67,65,72,40,3,23,59,69)_{\adfGL}$,

$(49,66,35,82,55,38,80,22,18,68,6,12,\adfsplit 62,53,77,73,78,13,61,2,48,91,16,3)_{\adfGL}$,

$(89,19,32,38,86,39,30,34,60,31,28,68,\adfsplit 55,26,81,90,9,74,4,91,48,52,7,5)_{\adfGL}$,

$(38,22,25,15,88,41,57,91,37,80,33,69,\adfsplit 73,71,46,17,6,53,83,32,28,66,85,35)_{\adfGL}$,

$(58,76,62,36,60,39,1,63,8,40,31,7,\adfsplit 56,64,85,35,27,75,33,78,25,57,26,4)_{\adfGL}$

\noindent under the action of the mapping $x \mapsto x + 4$ (mod 72) for $x < 72$,
$x \mapsto 72 + (x - 72 + 7 \mathrm{~(mod~21)})$ for $x \ge 72$.\eproof

\adfVfy{93, \{\{72,18,4\},\{21,3,7\}\}, -1, \{\{24,\{0,1,2\}\},\{21,\{3\}\}\}, -1} 


\end{document}